\documentclass[aos]{imsart}

\RequirePackage{amsthm,amsmath,amsfonts,amssymb}
\RequirePackage[numbers]{natbib}
\usepackage{fix-cm}
\usepackage{graphicx}
\usepackage{subcaption}
\usepackage{tikz}

\usepackage{lipsum} 

\usepackage{titlesec}

\RequirePackage{algorithm}
\RequirePackage{algorithmic}

\usepackage{breakcites}
\usepackage{xcolor}

\usepackage{color}
\usepackage{hyperref}
\usepackage{ulem}
\usepackage{textcomp}
\usepackage{comment}
\usepackage{xcolor}
\usepackage{amsthm}
\usepackage{bbm}
\usepackage[makeroom]{cancel}

\RequirePackage[shortlabels]{enumitem}

\startlocaldefs 

\newcommand{\pa}[1]{\left(#1\right)}
\newcommand{\cro}[1]{\left[#1\right]}
\newcommand{\ac}[1]{\left\{#1\right\}}

\def\cE{\mathcal{E}}

\newcommand{\E}{\operatorname{\mathbb{E}}}
\renewcommand{\P}{\operatorname{\mathbb{P}}}
\newcommand{\R}{\mathbb{R}}

\newcommand{\Gcal}{\mathcal{G}}

\newcommand{\Xcal}{\mathcal{X}}

\newcommand{\de}{\delta}
\newcommand{\ep}{\epsilon}

\newcommand{\om}{\omega}
\newcommand{\si}{\pi}
\newcommand{\al}{\alpha}
\newcommand{\wH}{\widehat H}
\newcommand{\wD}{\widehat D}

\newcommand{\bPi}{\boldsymbol{\Pi}}

\DeclareMathOperator*{\argmax}{arg\,max}
\DeclareMathOperator*{\argmin}{arg\,min}

\newcommand{\cA}{\mathcal{A}}

\newcommand{\lde}{\mathcal{D}}

\newcommand{\ap}{\mathcal{A}}

\theoremstyle{plain}
\newtheorem{thm}{Theorem}[section]

\newtheorem{extended-hyp}{Extension}[section]
\newtheorem{extension-hyp}{Relaxation}[section]
\newtheorem{prop}[thm]{Proposition}
\newtheorem*{prop*}{Proposition}
\newtheorem{lem}[thm]{Lemma}

\newtheorem{cor}[thm]{Corollary}

\theoremstyle{remark}
\newtheorem{remark}{Remark}[section]

\theoremstyle{definition}
\newtheorem{definition}[thm]{Definition}

\newcommand{\va}{\sigma}  
\newcommand{\der}{\gamma}  

\endlocaldefs

\begin{document}

\begin{frontmatter}
\title{Minimax optimal seriation  in polynomial time}
\runtitle{Optimal seriation}

\begin{aug}
\author[A]{\fnms{Yann}~\snm{Issartel}\ead[label=e1]{yann.issartel@telecom-paris.fr}},
\author[B]{\fnms{Christophe}~\snm{Giraud}\ead[label=e2]{christophe.giraud@universite-paris-saclay.fr}}
\and
\author[C]{\fnms{Nicolas}~\snm{Verzelen}\ead[label=e3]{Nicolas.Verzelen@inrae.fr}}
\address[A]{LTCI, Télécom Paris, Institut Polytechnique de Paris\printead[presep={,\ }]{e1}}

\address[B]{Laboratoire de Math\'ematiques d'Orsay, Universit\'e Paris-Saclay\printead[presep={,\ }]{e2}}

\address[C]{MISTEA, INRAE, Institut Agro,
Univ. Montpellier\printead[presep={,\ }]{e3}}

\end{aug}

\begin{abstract}
We consider the seriation problem, whose goal is to recover a hidden ordering from a noisy observation of a permuted Robinson matrix.
We establish sharp minimax rates under average-Lipschitz conditions that strictly extend the bi-Lipschitz framework of~\cite{giraud2021localization}.
We further design a polynomial-time algorithm that attains these optimal rates, thereby resolving two open questions raised in~\cite{giraud2021localization}.
Finally, our analysis extends to a broader class of matrices beyond those generated by exact permutations.
\end{abstract}


\end{frontmatter}


\section{Introduction}
Seriation aims to recover the unknown ordering of $n$ objects from noisy pairwise measurements.
This problem dates back to early applications in archaeology, where researchers reconstructed the chronological order of graves from the observation that graves close in time tend to contain similar sets of artifacts~\cite{Robinson51}. In modern data science, seriation arises in areas such as sparse matrix reordering~\cite{barnard1995spectral}, genome sequencing~\cite{garriga2011banded,bioinfo17}, network synchronization~\cite{Clock-Synchro04,Clock-Synchro06}, and interval graph identification~\cite{fulkerson1965incidence}.

In the statistical formulation of seriation, one observes a symmetric matrix {{\small $A$}} of noisy pairwise similarities between $n$ objects. These similarities are assumed to reflect an underlying permutation {{\small $\pi=(\pi_1,\ldots,\pi_n)$}} of {{\small $[n]$}}, in the sense that {{\small $A_{ij}$}} tends to be large when the positions {{\small $\pi_i$}} and {{\small $\pi_j$}} are close, and small when they are far apart. A common formalization of this structure relies on Robinson matrices~\cite{fogel2013convex,recanati2018reconstructing,janssen2020reconstruction,giraud2021localization,natik2021consistency}; that is, matrices whose rows and columns are unimodal, attaining their maxima on the diagonal. In our setting, we assume that the mean matrix {{\small $\E A$}} is Robinson up to a permutation of its rows and columns by {{\small $\pi$}}.

The objective is to recover the permutation {{\small $\pi$}} from the observed matrix {{\small $A$}}. The performance of an estimator {{\small $\hat \pi$}} is often assessed via the normalized maximum error
{{\small $\max_{i\in[n]} |\hat \pi_i - \pi_i|/n$}}~\cite{janssen2020reconstruction,giraud2021localization,natik2021consistency}; 
a rigorous definition of this loss is given in Section~\ref{section-setting}. Because this criterion is stronger than global metrics such as the Frobenius or Kendall's tau distances, 
bounds obtained under this loss immediately yield corresponding guarantees for these weaker measures.


\paragraph*{Limitations of existing methods.}   
The noiseless seriation problem {{\small $(A = \E A)$}} was solved by Atkins et al.~\cite{atkins1998spectral}, who proposed a spectral algorithm that recovers the permutation {{\small $\pi$}} exactly and efficiently. In the noisy case, however, the situation is less clear. Convex relaxations~\cite{fogel2013convex} and spectral methods~\cite{giraud2021localization,natik2021consistency,cai2022matrix} have been analyzed under strong assumptions on the mean matrix {{\small $\E A$}}, for example when {{\small $\E A$}} is Toeplitz with a large spectral gap, in which case optimal statistical rates can be achieved~\cite{giraud2021localization}. Beyond this specific case, theoretical guarantees for spectral algorithms remain limited, and several works suggest a degradation of performance on more general matrices~\cite{rocha2018recovering,janssen2020reconstruction,giraud2021localization,cai2022matrix}. This gap highlights the need for seriation methods that perform reliably even when {{\small $\E A$}} is not Toeplitz.  

A structural alternative is provided by the class of \textit{bi-Lipschitz matrices}~\cite{giraud2021localization}, which impose Lipschitz-type constraints on variations across rows and columns. This framework has two main advantages. First, it captures heterogeneous structures that cannot be represented within the Toeplitz setting. Second, it comes with explicit regularity parameters {{\small $\alpha,\beta$}}, which allow one to quantify how algorithmic performance depends on structural smoothness.

Over the bi-Lipschitz class, the optimal minimax rate {{\small $\sqrt{\log(n)/n}$}} was established in~\cite{giraud2021localization} by a super-polynomial-time procedure, showing that this rate can be attained in principle. However, no polynomial-time algorithm is currently known to achieve this rate, and the precise dependence of the optimal error on matrix regularity and noise level remains unknown.

\subsection{Contributions} 
This paper addresses the limitations above by introducing {{\small SABRE}} (Seriate by Aggregating Bisections and Re-Evaluating), a new polynomial-time algorithm with provable optimality guarantees. 
Our main contributions are as follows:

\begin{itemize}\setlength{\itemsep}{0.7em}
    \item \textbf{Average-Lipschitz framework.}
    We introduce the class of average-Lipschitz matrices, which strictly extends the bi-Lipschitz framework of~\cite{giraud2021localization} and accommodates natural Robinson structures outside that class.

    \item \textbf{Minimax rates.}
    We establish sharp minimax rates for seriation of average-Lipschitz matrices with respect to the maximum error {{\small $\max_{i\in[n]} |\hat\pi_i - \pi_i|/n$}}.
    These rates describe the statistical difficulty of seriation by making explicit how the optimal error depends on noise level and matrix regularity, thereby addressing a point left open in the bi-Lipschitz analysis of~\cite{giraud2021localization}.
    We also derive error bounds under alternative loss functions such as Frobenius and Kendall's tau distances.

    \item \textbf{Efficient algorithm.}
    {{\small SABRE}} attains these minimax rates while running in {{\small $O(n^3)$}} time. 
    In particular, this yields the first polynomial-time estimator with minimax guarantees in the bi-Lipschitz setting, resolving an open question raised in~\cite{giraud2021localization}.

    \item \textbf{Approximate orderings.}
    Our analysis extends to latent orders represented by approximate permutations ({{\small $\pi\in[n]^n$}}), allowing repeated or missing positions.  
    This relaxation accommodates a broader family of matrices while preserving the minimax guarantees of {{\small SABRE}}.
\end{itemize}



\subsection{Related literature}
\label{subsection:intro:example}

The seriation problem belongs to a broader class of learning tasks where data are disordered by an unknown permutation, 
and the goal is to recover or approximate the latent ordering. 
Examples include ranking~\cite{braverman2009sorting,mao2018minimax,chen2019spectral,chatterjee2019estimation}, 
feature matching~\cite{collier2016minimax,galstyan2022optimal}, 
and matrix estimation under shape constraints~\cite{flammarion2019optimal,chatterjee2019estimation,mao2018towards,hutter-Mao-Rigollet-2020Monge-estimation}. 
Closer to our setting, seriation has been studied for Robinson matrices~\cite{atkins1998spectral,fogel2013convex,FogelNeurIPS,recanati2018reconstructing}, 
with more recent statistical analyses in~\cite{janssen2020reconstruction,giraud2021localization,cai2022matrix}. 
Each problem comes with its own assumptions and goals, so that techniques are not always transferable.

\paragraph*{Noiseless seriation.}
Much of the early work on seriation focused on the noiseless case. 
Efficient algorithms have been proposed via spectral methods~\cite{atkins1998spectral} 
and convex optimization~\cite{fogel2013convex}, 
with exact recovery guarantees for Robinson matrices~\cite{atkins1998spectral} and their toroidal variants~\cite{recanati2018reconstructing}. 
Other related contributions include combinatorial algorithms for recognizing Robinson matrices~\cite{prea2014optimal,carmona2024modules} 
and for best Robinson matrix approximation under $\ell_p$-norms~\cite{chepoi2009approximation}.

\paragraph*{Noisy seriation.}
Statistical analysis in the noisy setting has emerged only recently~\cite{janssen2020reconstruction,giraud2021localization,natik2021consistency,cai2022matrix}. 
Most existing results rely on strong structural assumptions, most notably Toeplitz matrices. 
In this regime, spectral algorithms can succeed under large spectral gaps and, with suitable post-processing, achieve small maximum error~\cite{rocha2018recovering,natik2021consistency}, and even the optimal rate {{\small $\sqrt{\log(n)/n}$}}~\cite{giraud2021localization}.
In contrast, our analysis is carried out within the class of average-Lipschitz matrices.  
These two frameworks (Toeplitz / Lipschitz) are structurally distinct and not nested; see Section~\ref{section-discussion} for further discussion.

Distance-based methods have also been proposed for seriation on Toeplitz matrices~\cite{cai2022matrix}. 
Although {{\small SABRE}} leverages distance information, its construction and analysis are fundamentally different.

Beyond  Toeplitz settings, Janssen et al.~\cite{janssen2020reconstruction} introduce a seriation algorithm for network data, which is based on thresholding the squared adjacency matrix.

\paragraph*{Loss functions.}
Different criteria have been used to evaluate seriation procedures. The normalized maximum error~\cite{janssen2020reconstruction,giraud2021localization,natik2021consistency} is appropriate when accurate recovery of all positions is required, but may be overly stringent in some settings. 
Weaker criteria such as the Frobenius loss have therefore been used as alternatives, for instance in~\cite{cai2022matrix}.

\paragraph*{Related distance-based approaches.}
For reordering Toeplitz matrices, Cai et al.~\cite{cai2022matrix} exploit that consecutive items in the latent order have similar rows,  so that row {{\small $\ell_1$}} distances can be used to track ordering proximity.
Beyond seriation, distance-based algorithms are also studied in matrix estimation problems.  
For example, Hütter et al.~\cite{hutter-Mao-Rigollet-2020Monge-estimation} analyze Monge matrices under noisy observations, building on the same principle with an appropriate distance.  
While their objectives differ from ours, these works reinforce the general principle that latent orderings can be partially reconstructed from suitable distance comparisons.


\subsection{Organization and Notation}
Section~\ref{section-setting+assumption+motivation} introduces the  model, the loss, and the class of average-Lipschitz matrices.
Section~\ref{section-sabre} presents the {{\small SABRE}} algorithm.
Section~\ref{section:bilip-mat} states our main results. 
Section~\ref{section-discussion} concludes with a discussion.
Section~\ref{section-algo-analysis} gives a proof sketch.
Technical proofs and additional details are deferred to the appendix.

\paragraph*{Notation.}
We denote positive numerical  constants by {{\small $c$}} and {{\small $C$}}. The notations {{\small $a \lesssim b$}} and {{\small $a \asymp b$}} mean that there exist  constants {{\small $c,C>0$}} such that {{\small $a \leq Cb$}} and {{\small $cb \leq a \leq Cb$}}.
If a constant depends on parameters {{\small $\alpha,\beta$}}, we denote it by {{\small $C_{\alpha \beta}$}}, and similarly we write {{\small $\asymp_{\alpha,\beta}$}}. 
We denote {{\small $\{1,\ldots,n\}$}} by {{\small $[n]$}}, {{\small $\max(a,b)$}} by {{\small $a \vee b$}}, and {{\small $\min(a,b)$}} by {{\small $a \wedge b$}}. 
For any permutation {{\small $\pi$}}, we shorten {{\small $\pi(i)$}} to {{\small $\pi_i$}}.  
We write {{\small $F_i$}} for the $i$-th row of a matrix {{\small $F$}}, and 
{{\small $\|x\|$}} for the {{\small $\ell_2$}}-norm of a vector~{{\small $x$}}.  
For a sequence {{\small $\zeta = \zeta_n$}}, the notation {{\small $\zeta = o(n)$}} means {{\small $\zeta/n \to 0$}} as {{\small $n\to\infty$}}.
We use {{\small $\mathbf{1}_{\{\cdot\}}$}} for indicator functions, and {{\small $\mathbf{1}_n\in\mathbb{R}^n$}} for the all-ones vector. For sets {{\small $G$}} and {{\small $G'$}}, {{\small $G\setminus G'$}} denotes {{\small $\{k \in G: k \notin G'\}$}}.


\section{Statistical setting and structural assumptions}
\label{section-setting+assumption+motivation}

We introduce the statistical model and the seriation problem,
then describe the structural assumptions on the signal matrix, with illustrative examples.

\subsection{Model}
\label{section-setting}

A symmetric matrix {{\small $F$}} is called a \textit{Robinson matrix} if its rows and columns are unimodal with maxima on the diagonal; that is,
\begin{equation}\label{assumption:robinson-matrix}
\forall \, k \le i<j: \quad F_{jk}<F_{ik}, 
\qquad \qquad
\forall \, i<j \le k: \quad F_{ik}<F_{jk}.
\end{equation}
Equivalently, the entries of {{\small $F$}} decrease monotonically as one moves away from the diagonal. 
For a permutation {{\small $\pi:[n]\to[n]$}}, we denote by {{\small $F_\pi$}} the matrix obtained by permuting the rows and columns of {{\small $F$}} according to~{{\small $\pi$}}, 
that is {{\small $(F_\pi)_{ij} = F_{\pi_i\pi_j}$}}.

We observe a noisy permuted Robinson matrix of the form
\begin{equation}\label{eq:model}
A \;=\; F_\pi + \sigma E ,
\end{equation}
where {{\small $F$}} is an unknown Robinson matrix, {{\small $\pi$}} is an unknown permutation of {{\small $[n]$}}, and {{\small $E$}} is a symmetric noise matrix.  
The entries {{\small $\{E_{ij},\, i \le j\}$}} are independent, mean-zero, sub-Gaussian random variables with variance proxies upper bounded by~{{\small $1$}}, in the sense that {{\small $\E[e^{sE_{ij}}]\le \exp(s^2/2)$}} for all {{\small $s\in\R$}} and  all {{\small $i \le j$}}
(see, e.g., \cite{vershynin2018high,rigollet2023high}). 
Thus {{\small $\sigma$}} provides an upper bound on the noise scale in {{\small $A$}}.
This setting includes, for instance, random graph models where {{\small $A_{ij}$}} follows a Bernoulli distribution with mean {{\small $F_{\pi_i \pi_j}$}}.  
We write {{\small $\P_{(F,\pi)}$}} for the distribution of~\eqref{eq:model}. 

While we present the model in terms of an exact permutation {{\small $\pi$}} of {{\small $[n]$}}, we also consider a relaxed notion of latent order. In Section~\ref{extension-section}, we introduce \textit{approximate permutations}, 
allowing vectors {{\small $\pi\in[n]^n$}} with repeated values or small gaps, provided they remain well spread.


\subsection{Seriation problem and loss}
\label{section-objective}

The goal is to estimate the latent permutation {{\small $\pi$}} from the single observation {{\small $A$}}.  
There is a minor lack of identifiability, since the model~\eqref{eq:model} remains unchanged if both the latent order and the underlying matrix are reversed  
-- that is, when {{\small $(F,\pi)$}} is replaced by {{\small $(F^{\mathrm{rev}},\pi^{\mathrm{rev}})$}},  
where {{\small $\pi^{\mathrm{rev}}_i = n + 1 - \pi_i$}} and {{\small $F^{\mathrm{rev}}$}} denotes the matrix obtained by reversing  the order of its rows and columns.   
Following prior work~\cite{janssen2020reconstruction,giraud2021localization,natik2021consistency}, we evaluate an estimator {{\small $\hat\pi$}} by the normalized maximum error
\begin{equation}\label{eq:max-loss}
L_{\max}(\hat{\pi},\pi)
= \frac{1}{n}\Bigg(
\max_{i\in[n]} |\hat\pi_i-\pi_i| 
\;\wedge\;
\max_{i\in[n]} |\hat\pi_i-\pi^{\mathrm{rev}}_i|
\Bigg).
\end{equation}

The max-error loss~\eqref{eq:max-loss} is a stringent criterion, requiring accurate recovery of all positions. 
We will use it as our primary performance measure throughout the paper.


\subsection{Structural assumptions}
\label{section-assumption}

A natural refinement of the Robinson structure~\eqref{assumption:robinson-matrix} is given by \textit{bi-Lipschitz}  matrices introduced in~\cite{giraud2021localization}, which quantify the rate of decay along each row and column as follows.

\begin{definition}[Bi-Lipschitz matrices]
\label{defi:BL}
Let {{\small $\alpha, \beta > 0$}}.  
A matrix {{\small $F\in[0,1]^{n\times n}$}} is said to be \textit{bi-Lipschitz} if,  
for every {{\small $i<j$}}:

\begin{enumerate}[leftmargin=1.3em,label=--]
\item \textit{Pointwise {{\small $\ell_\infty$}} upper bound:}\;  
{{\small $\max_{k\in[n]} |F_{ik}-F_{jk}| \le \beta\,|i-j|/n.$}}

\item \textit{Pointwise lower bound:} 
\[
F_{ik}-F_{jk} \ \ge \ \alpha\,\frac{|i-j|}{n}
\quad \text{for all }  k\le i,
\qquad
F_{jk}-F_{ik} \ \ge  \alpha\,\frac{|i-j|}{n}
\quad \text{for all }  k\ge j.
\]
\end{enumerate}

We denote this class by {{\small $\mathcal{BL}(\alpha,\beta)$}}.
\end{definition}

Intuitively, the bi-Lipschitz condition strengthens the Robinson shape~\eqref{assumption:robinson-matrix} by enforcing uniform lower and upper control on the decay of each row and column. This refinement rules out degenerate cases where the similarity matrix {{\small $F$}} is nearly flat, which is essential when aiming to control the maximal loss {{\small $L_{\max}$}} in~\eqref{eq:max-loss}.

However, the bi-Lipschitz condition can be overly restrictive in practice,
as it enforces uniform regularity, thereby failing to accommodate matrices with localized irregularities such as abrupt transitions; see Section~\ref{section:examples-BL-AL} for examples.
We therefore introduce a weakened form of the bi-Lipschitz condition that retains local regularity while relaxing its uniform constraints.
We refer to this broader family as \textit{average-Lipschitz}  matrices.

\begin{definition}[average-Lipschitz matrices]
\label{defi:PBL}
Let {{\small $\alpha, \beta, r,r'>0$}}.  
A matrix {{\small $F\in[0,1]^{n\times n}$}} is said to be \textit{average-Lipschitz} if it satisfies:

\begin{enumerate}[leftmargin=1.3em,label=(\roman*)]
\item \textbf{Local conditions.}  
For every {{\small $i<j$}} with {{\small $j-i\le r n$}}:
  \begin{enumerate}[leftmargin=1.3em,label=--]
  \item \textit{Average {{\small $\ell_2$}} upper bound:}\; 
  {{\small $\|F_i - F_j\| \le \beta\,|i-j|/\sqrt{n}.$}}
  \item \textit{Average {{\small $\ell_1$}} lower bound:}
  \[
  \Bigg(
  \sum_{k < i - c_0 n} (F_{ik} - F_{jk})
  \ \vee\
  \sum_{k > j + c_0 n} (F_{jk} - F_{ik})
  \Bigg)
  \ \ge\  \alpha\,|i-j|.
  \]
  \end{enumerate}

\item \textbf{Non-collapse at large distance.}  
If {{\small $|i-j|>r n$}}, then {{\small $\|F_i - F_j\| > r'\sqrt{n}$}}.
\end{enumerate}

We denote this class by {{\small $\mathcal{AL}(\alpha,\beta,r,r')$}},  
where {{\small $c_0 = 1/32$}}.
\end{definition}

The above definition weakens the classical bi-Lipschitz condition:
it replaces the pointwise {{\small $\ell_\infty$}} upper bound by an average {{\small $\ell_2$}} control,
and relaxes the pointwise lower bound to an average {{\small $\ell_1$}} separation. (For technical reasons, a small trimming of {{\small $c_0 n$}} indices is applied in the lower-bound condition to avoid boundary effects.)

Beyond these local constraints, the non-collapse condition acts at larger scales.
It prevents distant indices from having nearly identical rows.

\begin{remark}[On the need for a separation condition]
The {{\small $\ell_1$}} lower bound reflects a separation requirement specific to the {{\small $L_{\max}$}} loss:  
unlike average-Losses such as the Frobenius norm, {{\small $L_{\max}$}} demands accurate recovery of every position, 
which is impossible without sufficient separation between distinct rows.  
At the matrix level, this means that, to recover the correct order under noise, pairs of rows {{\small $F_i$}} and {{\small $F_j$}} must differ sufficiently on at least one subset of indices.  
In the Robinson structure, these subsets naturally lie within {{\small $\{k<i\}$}} or {{\small $\{k>j\}$}}, 
and the {{\small $\ell_1$}} condition simply enforces such a separation.
\end{remark}

The average-Lipschitz condition serves as our standing assumption throughout the paper.  
It generalizes the classical bi-Lipschitz framework by accommodating matrices with more heterogeneous structure, as illustrated in the next subsection.
Formally, {{\small $\mathcal{BL}(\alpha,\beta)\subset\mathcal{AL}(\alpha/4,\beta,r,r')$}} (see Appendix~\ref{append-proof-cor-bi-lip}), and this inclusion is strict.


\subsection{Illustrative examples}
\label{section:examples-BL-AL}

We now provide canonical examples that highlight the limitations of the bi-Lipschitz class and the flexibility of the average-Lipschitz condition.

We consider  matrices of the form {{\small $F_{ij}=a_{|i-j|}/n$}},
where the sequence {{\small $(a_\ell)$}} is non-increasing, so that {{\small $F$}} is Robinson.
The following examples give different ways in which such matrices can violate the bi-Lipschitz condition while still satisfying the average-Lipschitz condition.

\begin{enumerate}[leftmargin=1.3em,label=(\roman*)]
\setlength{\itemsep}{.7em}

\item \textbf{Vanishing long-range similarities.}  
Consider  {{\small $a_\ell = n/2 - \ell$}} for {{\small $\ell \le n/2$}} 
and {{\small $a_\ell = 0$}} otherwise.  
The matrix {{\small $F$}}  is not bi-Lipschitz because {{\small $F_{ij}=0$}} whenever {{\small $|i-j|\ge n/2$}}, 
violating the long-range lower bound required in Definition~\ref{defi:BL}.  
However, this loss of contrast at large distances is compatible with the average-Lipschitz condition, 
which only enforces average row separation.

\item \textbf{Localized jump.}
Consider {{\small $a_\ell = n - \alpha \ell$}} for {{\small $\ell \le \ell_0$}} 
and {{\small $a_\ell = n - \alpha \ell - \Delta$}} for {{\small $\ell \ge \ell_0 +1$}}.  
Under the bi-Lipschitz condition, the pointwise {{\small $\ell_\infty$}} upper bound must pay for the jump:  
the smallest admissible constant {{\small $\beta$}} grows linearly with {{\small $\Delta$}}.  
Under the average-Lipschitz condition, however,
{{\small $\|F_i-F_{i+1}\|\le (\alpha + O(\Delta/\sqrt{n}))/\sqrt{n}$}},  
so the same matrix belongs to {{\small $\mathcal{AL}$}} with  
{{\small $\beta = \alpha + O(\Delta/\sqrt{n})$}}.  
Hence, a localized discontinuity inflates the bi-Lipschitz constant but remains negligible under the average {{\small $\ell_2$}} control.

\item \textbf{Local plateaus.} 
Consider  {{\small $a_\ell = n - \lfloor \ell /2 \rfloor$}}, 
where {{\small $\lfloor x \rfloor$}} denotes the integer part of {{\small $x$}}.  
Consecutive differences vanish for every even index ({{\small $a_{2m}=a_{2m+1}$}}), creating flat local segments that violate the pointwise lower bound of Definition~\ref{defi:BL}. 
The same phenomenon occurs for any constant plateau length, e.g., when {{\small $a_\ell = n - \lfloor \ell / c \rfloor$}} with {{\small $c\ge2$}}.  
In contrast, such local repetitions remain compatible with the average-Lipschitz framework,  
which tolerates small regions of flatness as long as row separation holds on average.
\end{enumerate}

These simple Robinson matrices, though violating the bi-Lipschitz condition, remain realistic  for seriation: their structure is far from degenerate and does not suggest any inherent barrier to  recovering the latent order under the statistical model~\eqref{eq:model}. 
The average-Lipschitz condition accommodates these matrices, extending the bi-Lipschitz class to a broader and more natural family of Robinson structures.


\section{The SABRE algorithm}
\label{section-sabre}

We now present our estimation procedure, \textit{Seriate by Aggregating Bisections and Re-Evaluating} ({{\small SABRE}}).
The algorithm takes as input the noisy permuted Robinson matrix 
{{\small $A$}} from~\eqref{eq:model} and outputs an estimator 
{{\small $\hat\pi$}} of the latent permutation {{\small $\pi$}}. 
Its construction consists of three stages:\vspace{-0.5em}
\begin{enumerate}\itemsep0.7em
  \item \textsc{Distance estimation.} Construct an empirical estimator {{\small $\widehat D$}} of the
    population row-distance matrix
    \begin{equation}\label{def-dist-[n]}
    D^*_{ij} := \sqrt{n}\,\|F_{\pi_i}-F_{\pi_j}\| , \qquad i,j\in[n],
    \end{equation}
    where {{\small $D^*_{ij}$}} measures the discrepancy between rows {{\small $i$}} and {{\small $j$}} of the signal matrix {{\small $F_\pi=\E A$}}.
    Under the Robinson structure of {{\small $F$}}, the distances {{\small $D^*_{ij}$}} tend to be small when {{\small $\pi_i$}} and {{\small $\pi_j$}} 
    are close in the latent order, and larger when these positions are farther apart.

  \item \textsc{First seriation.} From {{\small $\widehat D$}}, build a matrix {{\small $H$}}, intended to approximate the ideal comparison matrix {{\small $H^*$}} defined from the latent permutation {{\small $\pi$}} by
  \begin{equation}
  \label{defi-H*}
  H^*_{ij} =
  \begin{cases}
    -1, & \text{if } \pi_i < \pi_j, \\
    1, & \text{if } \pi_j < \pi_i,
  \end{cases}
  \qquad H^*_{ii}=0.
  \end{equation}
  While {{\small $H^*$}} is fully determined by {{\small $\pi$}} and has no undecided entries,
  its empirical counterpart {{\small $H$}} may contain zeros when comparisons are left undecided.

  \item \textsc{Refined seriation.} Re-evaluate  undecided comparisons
using {{\small $H$}} and {{\small $A$}} to produce a refined comparison
matrix {{\small $\widehat H$}}, and derive {{\small $\hat\pi$}} from the
row sums of {{\small $\widehat H$}}.
\end{enumerate}

Each of the three stages above is presented in detail in Sections~\ref{subsection:l2-dist}, \ref{section:idea-algo}, and~\ref{subsection:stage2}.

Note that since the latent permutation {{\small $\pi$}} is identifiable only up to a reversal, 
the comparison matrix {{\small $H^*$}} is identifiable only up to a global factor {{\small $\pm 1$}}. 
Accordingly, {{\small SABRE}} reconstructs {{\small $H^*$}} modulo this inherent indeterminacy.

\textit{Computational complexity.} 
Each stage can be implemented in time {{\small $\mathcal{O}(n^3)$}}, so the overall complexity of 
{{\small SABRE}} is polynomial in {{\small $n$}}.

\paragraph*{Intuition.}
Empirical distances {{\small $\widehat D$}} are informative only for pairs {{\small $(i,j)$}} for which the population distance {{\small $D^*_{ij}$}} exceeds the noise scale.
The first-seriation stage uses these informative pairs to infer the corresponding comparisons and store them in the partial matrix {{\small $H$}}.
In contrast, when {{\small $D^*_{ij}$}} is at or below the noise level, {{\small $\widehat D_{ij}$}} cannot reliably separate the two rows, and the associated entries of {{\small $H$}} remain undecided.
The refinement step addresses these remaining cases: using the ordering information already present in {{\small $H$}}, it identifies the entries of {{\small $A$}} that are informative for  resolving the undecided comparisons, thereby producing the refined matrix {{\small $\widehat H$}} from which {{\small $\hat\pi$}} is obtained.


\paragraph*{Inputs and outputs.}
For clarity, we summarize the inputs and outputs of {{\small SABRE}}.
\smallskip 
{{\small
\begin{algorithmic}[1]
\REQUIRE {{\small $A\in\R^{n\times n}$}}; upper bound {{\small $\sigma$}} from~\eqref{eq:model}; tuning parameters {{\small $(\delta_1,\delta_2,\delta_3,\delta_4)$}}.
\ENSURE {{\small $\hat\pi\in[n]^n$}}.
\STATE Stage~1 (distance estimation) uses {{\small $A$}} and outputs {{\small $\widehat D$}}.
\STATE Stage~2 (first seriation) uses {{\small $(\widehat D,\delta_1,\delta_2,\delta_3)$}} and outputs {{\small $H$}}.
\STATE Stage~3 (refined seriation) uses {{\small $(A,H,\delta_4,\sigma)$}} and outputs {{\small $\widehat H$}} and {{\small $\hat\pi$}}.
\end{algorithmic}
}}

Only the refined-seriation stage uses {{\small $\sigma$}}, which in our model~\eqref{eq:model} denotes an upper bound on the noise scale in {{\small $A$}}.



\subsection{Distance estimation}
\label{subsection:l2-dist} 

The first step of {{\small SABRE}} is to estimate the latent distance matrix {{\small $D^*$}} introduced in~\eqref{def-dist-[n]}. 
Expanding the squared distance using the inner-product representation of the Euclidean norm, we obtain for any {{\small $i,j\in[n]$}},
\[
    n^{-1}(D^*_{ij})^2 
    \ = \ \|F_{\pi_i}-F_{\pi_j}\|^2 \ = \ \langle F_{\si_i}, F_{\si_i} \rangle 
        + \langle F_{\si_j}, F_{\si_j} \rangle 
        - 2 \langle F_{\si_i}, F_{\si_j} \rangle ,
\]
where {{\small $\langle F_u, F_v \rangle = \sum_{k=1}^n F_{uk}F_{vk}$}} denotes the inner product between rows  {{\small $F_u$}} and {{\small $F_v$}} of {{\small $F$}}. 
Thus, {{\small $D^*_{ij}$}} decomposes into two quadratic terms and one cross term.

\paragraph*{Estimating cross terms.}
The cross term {{\small $\langle F_{\si_i}, F_{\si_j}\rangle$}} admits the unbiased estimator 
\[\langle A_i, A_j \rangle = \sum_{k = 1}^n A_{ik} A_{jk}.\]
In contrast, the quadratic terms {{\small $\langle F_{\si_i}, F_{\si_i}\rangle$}} cannot be directly estimated by {{\small $\langle A_i,A_i\rangle$}}, since
\[
    \mathbb{E}[A_{ik}^2] = F_{\si_i \si_k}^2 + \sigma^2 \mathbb{E}[E_{ik}^2],
\]
so {{\small $\langle A_i,A_i\rangle$}} is biased. Unless the noise variances {{\small $\E[E_{ik}^2]$}} are known, constructing an unbiased estimator of {{\small $D^*$}} is in general impossible; see \cite{issartel2021estimation}.
In our sub-Gaussian setting, these variances are indeed unknown and may differ across entries.

\paragraph*{Nearest-neighbor approximation.}
A natural idea is to approximate the quadratic terms by cross terms. 
Let {{\small $m_i = \argmin_{t \neq i} D^*_{it}$}} be the nearest neighbor of {{\small $i$}} with respect to {{\small $D^*$}}. 
Then 
\[
    \langle F_{\si_i}, F_{\si_i}\rangle \approx \langle F_{\si_i}, F_{\si_{m_i}} \rangle 
    \ \approx \ \langle A_i, A_{m_i}\rangle.
\]
Since the indices {{\small $m_i$}} are unknown, we construct empirical proxies {{\small $\widehat m_i$}} from the data, and estimate the squared distance via
\begin{equation}
    \label{def-dist-estim-main-paper}
   n^{-1} (\widehat D_{ij})^2  \ := \    \langle A_i, A_{\widehat m_i}\rangle + \langle A_j, A_{\widehat m_j}\rangle - 2\langle A_i, A_j\rangle .
\end{equation}
This estimator follows the population expansion of {{\small $n^{-1}(D^*_{ij})^2$}} introduced at the start of the subsection, replacing each term by its empirical counterpart.

\paragraph*{Constructing the proxy {{\small $\widehat m_i$}}.}
Under strong homogeneity assumptions on {{\small $F$}} (e.g., Toeplitz structure), one may consider
{{\small $\widehat m_i = \argmax_{s\neq i}\langle A_i,A_s\rangle$}},
reflecting the intuition that the largest inner product typically corresponds to the most similar row, and therefore to the nearest index in the latent order.
However, such proxies may become unreliable in heterogeneous settings, including the Lipschitz framework considered here.

To obtain a more stable approximation, we adopt a robust nearest-neighbor proxy following~\cite{issartel2021estimation}.
For each {{\small $i\in[n]$}},
\[
\widehat m_i  \ :=  \ \argmin_{j \in [n], j \neq i} \ \max_{k \in [n] \setminus \{i,j\}} \, |\langle A_k, A_i - A_j \rangle| \; .
\]
This construction defines the quantity {{\small $\widehat m_i$}} used in the distance estimator~\eqref{def-dist-estim-main-paper}.  
The full pseudocode for this distance-estimation stage is provided in Appendix~\ref{appendix:pseudocode}.

\paragraph*{Guarantees.} 
Under the average-Lipschitz assumption of Section~\ref{section-assumption},
the estimator {{\small $\widehat D$}} approximates the 
distance matrix {{\small $D^*$}} with the accuracy established in
Section~\ref{sub-section-distance-analysis}.


\subsection{First seriation.}
\label{section:idea-algo}

In {{\small SABRE}}, the first-seriation step takes as input the estimator {{\small $\widehat D$}} constructed in Stage~1.  
Since this step can be applied to any  estimate of the  distance matrix
{{\small $D^*$}}, we describe it for a generic estimate {{\small $D$}}.

The objective is to construct a matrix {{\small $H$}} that agrees with the  comparison matrix {{\small $H^*$}} defined in~\eqref{defi-H*},
wherever the distance estimate {{\small $D$}} provides sufficiently reliable information.  
The construction proceeds through three components: (i) bisections, (ii) orientation, (iii) aggregation.


\paragraph*{Step (i): bisections.}
For each {{\small $i\in[n]$}}, the goal is to recover the two sets of indices lying to the left and to the right of {{\small $i$}} in the latent order {{\small $\pi$}}.
To this end, we construct a graph {{\small $\mathcal G_i$}} on the vertex set {{\small $[n]\setminus{\{i\}}$}} whose large connected components are intended to approximate these two 
regions.

{{\scshape Construction of the graph.}}\hspace{0.3em} 
For a fixed index {{\small $i$}}, the graph {{\small $\mathcal G_i$}} is constructed as follows.
\vspace{-1.8em}
\begin{enumerate}\setlength{\itemsep}{0.6em}
  \item \textit{Vertex set.}  
  The vertices are all indices in {{\small $[n]\setminus\{i\}$}}.

  \item \textit{Edges.}  
  Two vertices {{\small $k$}} and {{\small $l$}} are joined by an edge whenever they are close to each other and both sufficiently far from {{\small $i$}}.  
  Formally, given thresholds {{\small $(\delta_1,\delta_2)$}}, we add an edge between {{\small $k$}} and {{\small $l$}} if
  \[
    D_{kl}\le \delta_1
    \qquad\text{and}\qquad
    D_{ik}\ge \delta_2,\ \ D_{il}\ge \delta_2 .
  \]
  Figure~\ref{fig:graph-Gi} illustrates this step in the idealized case  {{\small $D_{ij}=|\pi_i-\pi_j|$}}.

  \item \textit{Filtering of small components.}  Connected components that do not reach indices far from {{\small $i$}} are typically small, carry little useful information, and are therefore removed.
  Formally, a third threshold {{\small $\delta_3$}} is used to discard any component that does not
  contain a vertex {{\small $k$}} with {{\small $D_{ik}\ge\delta_3$}}, thereby ensuring that only
  components extending far enough from {{\small $i$}} are kept.
  Among the remaining components, we keep the two largest ones, denoted {{\small $(G_i,G'_i)$}}.
\end{enumerate}

\begin{figure}[H]
\centering
\begin{tikzpicture}[scale=1.0]

    \draw (-4,0) -- (7,0);

    \draw (1,-0.1) -- (1,0.1);
    \node[] at (0.9,0.31) {$\pi_i$};

    \draw (6,-0.1) -- (6,0.1);
    \node[] at (6,-0.32) {$\pi_k$};
    \draw (6.9,-0.1) -- (6.9,0.1);
    \node[] at (7,-0.32) {$\pi_l$};

    \draw[>=latex, dashed, <->, teal] (6,0.25) -- (6.9,0.25);
    \node[teal] at (6.45,0.5) {$\leq \delta_1$};

    \draw[>=latex, dashed, <->, teal] (1,0.9) -- (6,0.9);
    \node[teal] at (3.7,1.15) {$\geq \delta_2$};

    \draw (-3.9,-0.1) -- (-3.9,0.1);
    \node[] at (-3.9,-0.32) {$\pi_{\tilde k}$};

    \draw (0.26,-0.1) -- (0.26,0.1);
    \node[] at (0.22,-0.32) {$\pi_{k'}$};

    \draw (1.6,-0.1) -- (1.6,0.1);
    \node[] at (1.75,-0.32) {$\pi_{\tilde l}$};

    \draw[>=latex, dashed, <->, purple] (-3.9,-1.38) -- (1.6,-1.38);
    \node[purple] at (-1.1,-1.63) {$> \delta_1$};


    \draw[>=latex, dashed, <->, purple] (0.26,-0.75) -- (1,-0.75);
    \node[purple] at (0.66,-1) {$< \delta_2$};
    


\end{tikzpicture}
\caption{Illustration of the construction of the graph {{\small $\mathcal G_i$}} in the idealized case {{\small $D_{ij}=|\pi_i-\pi_j|$}}, where positions are shown along the latent order axis.  
Nodes {{\small $k,l$}} form a {{\small $\delta_1$}}-close pair that is {{\small $\delta_2$}}-far from $i$, 
so they lie on the same side of {{\small $i$}} and are connected in {{\small $\mathcal G_i$}}.
By contrast, indices lying on opposite sides of {{\small $i$}} do not satisfy the edge conditions 
and are therefore not connected in {{\small $\mathcal G_i$}}.  
For example, {{\small $\tilde k,\tilde l$}} violate the {{\small $\delta_1$}}~condition, 
while {{\small $k',\tilde l$}} fail the {{\small $\delta_2$}}~condition.}
\label{fig:graph-Gi}
\end{figure}

\vspace{-1.3em}
{{\scshape Rationale.}} \hspace{0.3em} 
When two indices {{\small $k,l$}} are both far from {{\small $i$}} and close to each other, they must lie on the same side of {{\small $i$}} in the latent order.
As a result, any edge in {{\small $\mathcal G_i$}} connects vertices belonging to the same side, and each connected component of {{\small $\mathcal G_i$}} is confined to a single side of {{\small $i$}}.

If two indices are consecutive or near-consecutive in the latent order,
their latent distance {{\small $D^*_{kl}$}} is expected to be small,
so these pairs typically satisfy {{\small $D_{kl}\le\delta_1$}} and become
connected in~{{\small $\mathcal G_i$}}.
Chains of such {{\small $\delta_1$}}-close neighbors then propagate along the 
ordering, so that indices on a common side tend to form a single connected 
component.

Because the matrix {{\small $D$}} only provides a coarse approximation of which indices lie far from or close to {{\small $i$}}, 
the graph {{\small $\mathcal G_i$}} may contain small spurious components.
The threshold {{\small $\delta_3$}} filters out such components, retaining only those that 
include vertices sufficiently far from {{\small $i$}}.

The two remaining components {{\small $(G_i,G'_i)$}} then provide meaningful approximations of the left and right sets of indices relative to {{\small $i$}}.



\paragraph*{Step (ii): orientation.}
Each unordered pair {{\small $(G_i,G'_i)$}} is now relabeled as an oriented
pair {{\small $(L_i,R_i)$}}, where
{{\small $(L_i,R_i)\in\{(G_i,G'_i),(G'_i,G_i)\}$}}.
No new sets are constructed: orientation simply assigns a left/right label
to the two components, with the intent that {{\small $L_i$}} (resp.\ {{\small $R_i$}})
correspond to indices lying to the left (resp.\ right) of {{\small $i$}} in the
latent order.
The goal is for all pairs {{\small $(L_i,R_i)$}}, {{\small $i=1,\ldots,n$}}, to share a
consistent global orientation.

To this end, we fix an arbitrary reference index
{{\small $c$}} and arbitrarily orient {{\small $(G_c,G'_c)$}} as {{\small $(L_c,R_c)$}}.
For any {{\small $i\neq c$}}, the orientation {{\small $(L_i,R_i)$}} is then chosen so
that {{\small $L_i$}} aligns with {{\small $L_c$}}.
Intuitively, the component in {{\small $\{G_i,G_i'\}$}} having the larger
overlap with {{\small $L_c$}} should be labeled {{\small $L_i$}}, and the
remaining component labeled {{\small $R_i$}}.
The precise decision rule is given in
Appendix~\ref{appendix:pseudocode}.


\paragraph*{Step (iii): aggregation.}
Once all pairs {{\small $(L_i,R_i)$}} have been consistently oriented, we infer the relative order of any pair {{\small $(i,j)$}} whenever their left/right sets provide enough information.
Following the definition~\eqref{defi-H*} of {{\small $H^*$}}, we set {{\small $H_{ij}=-1$}} if these sets
indicate that {{\small $i$}} lies to the left of {{\small $j$}}, and {{\small $H_{ij}=1$}} in the opposite case.
Pairs for which no conclusion can be drawn are left undetermined ({{\small $H_{ij}=0$}}). 
The full pseudocode for this first-seriation stage is provided in Appendix~\ref{appendix:pseudocode}.

\paragraph*{Choice of tuning parameters.}
The thresholds {{\small $(\delta_1,\delta_2,\delta_3)$}} are chosen so that
{{\small $\delta_1 \lesssim \delta_2 \lesssim \delta_3$}},
with each stage operating at increasingly larger distance scales.
If an upper bound {{\small $\omega_n$}} on the distance estimation error is
available -- as established in Section~\ref{sub-section-distance-analysis} --
then empirical distances are reliable whenever
{{\small $D_{ij} \gtrsim \omega_n$}}.
This suggests taking {{\small $\delta_1 \asymp \omega_n$}}.

\paragraph*{Guarantees.}
Under suitable regularity conditions on the input  matrix {{\small $D$}}, all nonzero entries of {{\small $H$}} coincide with those of the comparison matrix {{\small $H^*$}}; see Section~\ref{subsection-first-seriation-analysis} for a formal statement.



\subsection{Refined seriation.}
\label{subsection:stage2}

The first-seriation stage has produced a partial comparison matrix {{\small $H$}} containing all relations that can be reliably recovered from {{\small $\widehat D$}}.  
In this final stage, we leverage the partial ordering encoded in {{\small $H$}} to determine the comparisons that remain unresolved. 
The outcome is a refined comparison matrix {{\small $\widehat H$}}, from which the estimator {{\small $\hat\pi$}} is obtained.


\paragraph*{Preliminary comparison sets and statistics.}
For each unresolved pair {{\small $(i,j)$}} with {{\small $H_{ij}=0$}}, the partial comparison matrix {{\small $H$}} provides directional information 
by identifying indices that lie on the same side of both {{\small $i$}} and {{\small $j$}}. This leads to the preliminary sets
\[
L_{ij}=\{k : H_{ik}=H_{jk}=1\}, 
\qquad 
R_{ij}=\{k : H_{ik}=H_{jk}=-1\},
\]
which collect indices placed to the left or to the right of both points in the latent ordering. Aggregating the differences
\[
\sum_{k\in L_{ij}} (A_{ik}-A_{jk})
\quad\text{or}\quad
\sum_{k\in R_{ij}} (A_{ik}-A_{jk})
\]
would, in principle, provide a coherent indication of the relative position of {{\small $i$}} and {{\small $j$}}, since the rows of a Robinson matrix vary 
monotonically.  These aggregated differences thus form natural \textit{comparison statistics} which motivate the refined seriation stage.


\paragraph*{Limitations of the preliminary comparison sets.}
The sets {{\small $L_{ij}$}} and {{\small $R_{ij}$}} are not directly usable for analysis:  being constructed from the coarse comparison matrix {{\small $H$}}, and thus
from the data matrix {{\small $A$}}, they remain statistically coupled with the rows {{\small $A_i$}} and {{\small $A_j$}} used in the above comparison 
statistics.  As a consequence, these comparison statistics depend on {{\small $(A_i,A_j)$}} both through the sets {{\small $L_{ij},R_{ij}$}} and 
through the summed differences themselves, creating a statistical dependence that prevents the use of standard concentration arguments.  
To mitigate this issue and obtain comparison sets that are essentially independent of {{\small $(A_i,A_j)$}}, we rely on a light sample-splitting scheme, described below.


\paragraph*{Sample splitting and modified comparison sets.}
A natural approach is to recompute the distance estimate on a subset of the matrix
{{\small $A$}} that excludes the rows {{\small $(A_i,A_j)$}} used in the comparison statistics. 
This yields an estimate independent of {{\small $(A_i,A_j)$}}.
Using it, we construct modified sets {{\small $\tilde L_{ij}$}} and {{\small $\tilde R_{ij}$}} 
that retain the directional property of {{\small $(L_{ij},R_{ij})$}} while substantially reducing their problematic dependence on {{\small $(A_i,A_j)$}}. 
The formal definition of {{\small $(\tilde L_{ij},\tilde R_{ij})$}} is given in Appendix~\ref{appendix:pseudocode}.

Computationally, a naive approach would be to use the leave-one-out subsets {{\small $S^{ij}=[n]\setminus\{i,j\}$}} to recompute the distance estimate for each unresolved pair {{\small $(i,j)$}}.
This would increase the computational cost by a factor equal to the number of such pairs,
which can be as large as {{\small $O(n^2)$}}.
To avoid this, we partition {{\small $[n]$}} at random into three disjoint, balanced subsets {{\small $(S^1,S^2,S^3)$}}.  
For each unresolved pair {{\small $(i,j)$}}, we then select a set among {{\small $(S^1,S^2,S^3)$}} that contains neither {{\small $i$}} nor {{\small $j$}} -- denote this set by {{\small $S^{t_{ij}}$}} -- and recompute the distance estimator on the corresponding submatrix of {{\small $A$}} indexed by {{\small $S^{t_{ij}}\times S^{t_{ij}}$}}.

\paragraph*{Computational cost.}
This sample-splitting step does not alter the overall complexity:  constructing {{\small $\tilde L_{ij},\tilde R_{ij}$}} has the same {{\small $O(n^3)$}} cost as working directly with the preliminary sets {{\small $L_{ij},R_{ij}$}}, so the runtime of {{\small SABRE}} remains
{{\small $O(n^3)$}}.


\paragraph*{Comparison statistics and decision rule.}
Given the modified comparison sets {{\small $(\tilde L_{ij}, \tilde R_{ij})$}}, we form the aggregate differences as in the preliminary construction,
setting 
{{\small $l=\sum_{k\in \tilde L_{ij}} (A_{ik}-A_{jk})$}}
and 
{{\small $r=\sum_{k\in \tilde R_{ij}} (A_{ik}-A_{jk})$}}.
Either statistic may provide a reliable indication of whether {{\small $i$}} should be placed before or after {{\small $j$}}. 
We therefore apply the rule
\[
\widetilde H_{ij}
=
\begin{cases}
-\mathrm{sign}(l) &\text{if } |l|\ge 5\sigma\sqrt{n\log n},\\[0.2em]
\phantom{-}\mathrm{sign}(r) &\text{if } |r|\ge 5\sigma\sqrt{n\log n},\\[0.2em]
0 &\text{otherwise,}
\end{cases}
\]
which assigns a comparison whenever one of the statistics {{\small $l$}} or {{\small $r$}} exceeds a high-probability upper bound on the noise fluctuations ({{\small $5\sigma\sqrt{n\log n}$}} in our implementation).
This rule is applied to all unresolved pairs (i.e., for which {{\small $H_{ij}=0$}}), yielding a 
matrix {{\small $\widetilde H$}} that  updates only those entries for which {{\small $H$}} was previously unassigned. 
For pairs already resolved in {{\small $H$}}, we set {{\small $\widetilde H_{ij}=0$}}.
A complete implementation of this refined-seriation stage is given  in Appendix~\ref{appendix:pseudocode}.


\paragraph*{Output.}
The refined comparison matrix is defined as {{\small $\widehat H = H + \widetilde H$}},  
where {{\small $\widetilde H$}} collects all updated entries.  
The estimator of the latent permutation is then obtained from the row sums of  
{{\small $\widehat H$}}:
\[
\hat\pi_i = \frac{\widehat H_i \mathbf{1}_n + n + 1}{2},
\qquad i\in[n],
\]
where {{\small $\mathbf{1}_n$}} denotes the all-ones vector in  {{\small $\R^n$}}.


\paragraph*{Guarantees.}
Under the average-Lipschitz assumption of Section~\ref{section-assumption}, this refined-seriation stage correctly resolves all pairs that are sufficiently separated in the latent order.  
The corresponding guarantees for {{\small $\widehat H$}} are stated in Section~\ref{section:refined-seriation-analysis-mainpaper}, 
and their implications for {{\small $\hat\pi$}} are  summarized in Section~\ref{section:bilip-mat}.



\section{Minimax optimal rates for average-Lipschitz matrices}
\label{section:bilip-mat}

We first derive non-asymptotic upper bounds for {{\small SABRE}} under the average-Lipschitz model with respect to the {{\small $L_{\max}$}} loss, and then specialize them to the bi-Lipschitz subclass. 
We also establish a matching minimax lower bound and extend the analysis to alternative losses and to the relaxed framework of approximate permutations.

\subsection{Upper bounds for {{\small SABRE}}}
\label{subsec:upper-bounds}

We now state our main result, which characterizes the performance of the algorithm {{\small SABRE}} in estimating the latent permutation {{\small $\pi$}} under model~\eqref{eq:model}. 
Let {{\small $\hat\pi$}} denote the output of {{\small SABRE}} with tuning parameters
{{\small $\delta_1 = n^{3/4}\log n$}} and {{\small $\delta_{k+1} = \delta_k \log n$}} for {{\small $k\in[3]$}}.
Here the parameters {{\small $(\alpha, \beta, r, r', \sigma)$}} are treated as fixed constants independent of {{\small $n$}}.

\begin{thm}
\label{thm-matrix-setting}
For any  {{\small $(\alpha, \beta, r, r', \sigma)$,}} there exists a constant {{\small $C_{\alpha \beta r r' \sigma}$}}  depending only on these parameters such that the following holds for all {{\small $n\geq C_{\alpha \beta r r' \sigma}$}}.  
If {{\small $F \in \mathcal{AL}(\alpha,\beta,r,r')$}} as in Definition~\ref{defi:PBL}, then with probability at least {{\small $1-1/n^2$}}, the output {{\small $\hat\pi$}} of {{\small SABRE}} satisfies
\begin{equation*}
   L_{\max}(\hat \si,\si) \ \ \leq \ \  40 \, \frac{\va}{\alpha}\,\sqrt{\frac{\log n}{n}} \enspace .
\end{equation*} 
\end{thm}

In other words, the polynomial-time algorithm {{\small SABRE}} localizes each position {{\small $\pi_i$}} within an error at most {{\small $40(\sigma/\alpha)\sqrt{\log n/n}$}}, a rate that improves with signal strength {{\small $\alpha$}}, sample size {{\small $n$}}, and smaller values of {{\small $\sigma$}}, which serves as an upper bound on the noise scale in~\eqref{eq:model}.
The numerical constant {{\small $40$}} is not optimized.

Importantly, {{\small SABRE}} requires no prior knowledge of the regularity of {{\small $F$}}, yet its performance adapts to it: more regular matrices reduce estimation difficulty.
As shown below, this rate is minimax optimal not only over the average-Lipschitz class, but also for simpler one-parameter families of Robinson matrices.

The tuning parameters {{\small $(\delta_1,\delta_2,\delta_3,\delta_4)$}} are set in a universal way, depending only on the sample size {{\small $n$}}. 
More precisely, Theorem~\ref{thm-matrix-setting} is stated under fixed model parameters 
{{\small $(\alpha,\beta,r,r',\sigma)$}}, and the choices of 
{{\small $(\delta_1,\delta_2,\delta_3,\delta_4)$}} are explicit functions of {{\small $n$}}, 
without relying on unknown regularity parameters such as {{\small $\alpha$}} or {{\small $\beta$}}. 
The rationale and constraints underlying these choices are discussed in 
Section~\ref{section-algo-analysis}. 
A more general version of the theorem, allowing 
{{\small $(\alpha,\beta,r,r',\sigma)$}} to vary with {{\small $n$}}, is given in 
Appendix~\ref{appendix-general-theorem}, where the tuning parameters must satisfy 
inequalities involving these model parameters.

\begin{remark}[Dependence on model parameters]
The bound in Theorem~\ref{thm-matrix-setting} depends on {{\small $(\alpha,\beta,r,r',\sigma)$}} both directly, through the factor {{\small $\sigma/\alpha$}} in the rate, and indirectly, through the sample size condition {{\small $n \ge C_{\alpha\beta r r'\sigma}$}}.
Intuitively, this latter condition ensures that the first-seriation step of {{\small SABRE}} produces a meaningful ordering rather than a random alignment. 
The threshold {{\small $C_{\alpha\beta r r'\sigma}$}} grows when {{\small $\alpha$}} decreases or {{\small $\beta$}} increases, reflecting the increased difficulty of {{\small $L_{\max}$}}-seriation under weaker regularity.
\end{remark}



We now instantiate Theorem~\ref{thm-matrix-setting} on the subclass 
{{\small $\mathcal{BL}(\alpha,\beta)$}} of bi-Lipschitz matrices.

\begin{cor}\label{cor:bi-lipschitz-matrix} 
For any {{\small $(\alpha, \beta, \sigma)$}}, there exists a constant 
{{\small $C_{\alpha  \beta \sigma}$}} depending only on {{\small $(\alpha,\beta,\sigma)$}} 
such that the following holds for all {{\small $n\geq C_{\alpha \beta \sigma}$}}.  
If {{\small $F\in \mathcal{BL}(\alpha,  \beta)$}} as in Definition~\ref{defi:BL}, 
then with probability at least {{\small $1-1/n^2$}}, the output {{\small $\hat \si$}} of {{\small SABRE}} satisfies 
\begin{equation}\label{eq:seriation-rate:upper-bound_paper_version}
   L_{\max}(\hat \si,\si) \ \ \leq \ \ 160 \, \frac{\va}{\alpha}\,\sqrt{\frac{\log n}{n}} \enspace .
\end{equation} 
\end{cor}

To our knowledge, {{\small $\hat\pi$}} is the first polynomial-time estimator achieving 
the minimax seriation rate {{\small $\sqrt{\log(n)/n}$}} over bi-Lipschitz matrices -- 
a rate previously attained only by exhaustive search methods~\cite{giraud2021localization}.
Moreover, Corollary~\ref{cor:bi-lipschitz-matrix} provides an explicit, optimal dependence of the error 
on the key problem parameters {{\small $(\alpha,\sigma)$}}, a feature absent from previous work. 

Section~\ref{section-algo-analysis} provides a proof sketch of Theorem~\ref{thm-matrix-setting}.  
The proof of Corollary~\ref{cor:bi-lipschitz-matrix} relies on the inclusion 
{{\small $\mathcal{BL}(\alpha,\beta)\subset\mathcal{AL}(\alpha/4,\beta,r,r')$}}, 
which is established in Appendix~\ref{appendix:proof-main-theorem}.


\subsection{Optimality} 
\label{subsection-optimality}  

We now establish that the rate {{\small $(\va/\alpha)\sqrt{\log(n)/n}$}} obtained in Theorem~\ref{thm-matrix-setting} is minimax optimal.
In particular, this rate cannot be improved, even within very simple subclasses of the average-Lipschitz family.

Given the generality of the {{\small $\mathcal{AL}$}} class, one might expect faster estimation to be possible in lower-dimensional parametric settings.
This intuition, however, turns out to be false: even for the one-parameter  family of Robinson matrices
\begin{equation*}
    (F_\alpha)_{ij} \ = \  1-  \frac{\alpha |i-j|}{n}\,, \qquad  i,j \in [n],
\end{equation*} 
no estimator can achieve a rate faster than {{\small $(\va/\alpha)\sqrt{\log(n)/n}$}}.
Knowing the parameter {{\small $\alpha$}} in advance -- so that the signal matrix is exactly {{\small $F = F_\alpha$}} -- does not alter this conclusion. This is formalized in the following theorem, which establishes the minimax lower bound.

\begin{thm}
\label{1st:lowerBound} 
For any {{\small $\alpha,\sigma >0$}}, there exist a constant {{\small $C_{\alpha\sigma}$}} (depending only on {{\small $\alpha, \sigma$}}) and a numerical constant {{\small $c>0$}}, such that for all {{\small $n\ge C_{\alpha\sigma}$}},
\[
 \inf_{\hat{\si}} \ \sup_{\si } \  \P_{(F_\alpha,\si)} \!\left\{\, L_{\max}(\hat{\si},\si) \ \geq\ c \frac{\va}{\alpha}\sqrt{\frac{\log n}{n}} \,\right\} \ \ \geq \ \ \frac{1}{2}\ .
\]
Here the infimum is over all estimators {{\small $\hat \pi$}}, 
the supremum is over all permutations {{\small $\pi$}} of {{\small $[n]$}}, 
and {{\small $\P_{(F_\alpha,\si)}$}} denotes the distribution in model~\eqref{eq:model} with {{\small $F=F_\alpha$}}.
\end{thm}

Note that {{\small $F_\alpha$}} belongs to the bi-Lipschitz class {{\small $\mathcal{BL}(\alpha,\beta)$}} for any {{\small $\beta \ge \alpha$}}; hence, the same lower bound also holds for the broader average-Lipschitz class.

The proof, given in Appendix~\ref{appendix:lower-bound}, is carried out under independent Gaussian noise {{\small $E_{ij}\sim\mathcal N(0,\sigma^2)$}} for {{\small $i<j$}}.
The lower bound is established in this homoscedastic Gaussian model, whereas our upper bound in Theorem~\ref{thm-matrix-setting} holds in the more general sub-Gaussian setting allowing heterogeneous noise levels.
Since these bounds coincide up to constant factors, possible heterogeneity in the noise variances has no effect on the minimax rate.


\subsection{Extensions to other losses}

We now extend our analysis to other loss functions commonly used  in reordering problems.   
We first consider permutation-based discrepancies, namely Kendall's tau and normalized {{\small $\ell_1$}} losses,
then turn to matrix-level errors measured by the Frobenius loss.

\smallskip
\noindent\textit{Kendall's tau and {{\small $\ell_1$}} losses.}
To quantify the disagreement between an estimator {{\small $\hat\pi$}} and the true ordering {{\small $\pi$}}, 
a standard measure in ranking problems~\cite{chen2022optimal} is Kendall's tau distance:
\[
K(\hat \pi, \pi) \ = \ \frac{1}{n} \sum_{1 \leq i<j \leq n} \mathbf{1}_{\{\mathrm{sign}(\pi_i - \pi_j)\,\mathrm{sign}(\hat \pi_i - \hat \pi_j)  < 0\}} \,,
\]
so that {{\small $n K(\hat \pi, \pi)$}} counts the number of inversions between {{\small $\hat \pi$}} and {{\small $\pi$}}. Here, {{\small $\mathbf{1}_{\{\cdot\}}$}} denotes the indicator function. 
Consistently with the identifiability up to a reversal (Section~\ref{section-objective}), we define the Kendall's tau loss as
\begin{equation}
\label{kendall-loss}
    L_K(\hat \pi, \pi) \ = \ K(\hat \pi, \pi) \ \wedge \ K(\hat \pi, \pi^{\mathrm{rev}}) \,.
\end{equation}
A related measure is the normalized {{\small $\ell_1$}} loss,
\[
L_{1}(\hat \pi, \pi) \ = \ \frac{1}{n} \sum_{i=1}^n |\hat \pi_i - \pi_i| \ \ \wedge\ \ \frac{1}{n} \sum_{i=1}^n |\hat \pi_i - \pi_i^{\mathrm{rev}}| \,,
\]
which satisfies the classical relation~\cite{diaconis-persi}
\begin{equation}
\label{equivalence-losses}
  \tfrac{1}{2} L_{1}(\hat \pi, \pi) \ \ \leq \ \  L_{K}(\hat \pi, \pi) \ \ \leq \ \ L_{1}(\hat \pi, \pi)\,.
\end{equation}
Hence the two losses are equivalent up to a constant factor, allowing us to translate our {{\small $L_{\max}$}} guarantees into Kendall’s tau loss.

\begin{cor}
\label{coro-kendall-tau}
Under the assumptions of Theorem~\ref{thm-matrix-setting}, the estimate {{\small $\hat \si$}} from {{\small SABRE}} satisfies, with probability at least {{\small $1-1/n^2$}},
\begin{equation*}
   L_{K}(\hat \si,\si)  \ \ \leq \ \  40 \, \frac{\va }{\alpha}\,\sqrt{n \log n} \,.
\end{equation*} 
\end{cor}

Refining the proof of Theorem~\ref{1st:lowerBound} through standard techniques should yield a lower bound of order {{\small $(\sigma/\alpha)\sqrt{n}$}} for {{\small $L_K$}} (see Remark~\ref{lower-bound-proof-conje}).
Combined with Corollary~\ref{coro-kendall-tau}, this suggests that {{\small SABRE}} is minimax optimal for {{\small $L_K$}} up to a logarithmic factor.
Whether this logarithmic gap is intrinsic or merely an artifact of current proof techniques remains an open question.

\begin{remark}
\label{lower-bound-proof-conje}
A possible route to obtain the {{\small $\sqrt{n}$}} lower bound is to combine the equivalence~\eqref{equivalence-losses} between {{\small $L_K$}} and {{\small $L_1$}} losses with a refinement of the proof of Theorem~\ref{1st:lowerBound}.  
Specifically, by applying Varshamov–Gilbert's lemma~\cite{massart}, one can construct a well-spaced set of {{\small $e^{O(n)}$}} permutations under {{\small $L_1$}}, which yields a lower bound of order {{\small $(\va/\alpha)\sqrt{n}$}} for {{\small $L_K$}}.
\end{remark}

\noindent\textit{Frobenius loss.}
We now turn to the Frobenius loss, which measures matrix-level discrepancies between orderings:
\begin{equation}
\label{Frob-loss}
    L_{F}(\hat \pi, \pi) \ \ = \ \ \|  F_{\hat \pi} - F_{\pi}\| \ \wedge \ \| F_{\hat \pi} - F_{\pi^{\mathrm{rev}}} \| \,,
\end{equation}
where {{\small $\|\cdot\|$}} denotes the Frobenius norm and  {{\small $F$}} is the signal matrix in~\eqref{eq:model}. 
This loss was also considered in~\cite{cai2022matrix}. 

\begin{cor}
\label{coro-frob-loss}
Under the assumptions of Corollary~\ref{cor:bi-lipschitz-matrix}, the estimate {{\small $\hat \si$}} from {{\small SABRE}} satisfies, with probability at least {{\small $1-1/n^2$}},
\begin{equation*}
   L_{F}(\hat \si,\si)  \ \ \leq \ \  320 \, \frac{\beta \va }{\al}\,\sqrt{n \log n} \,.
\end{equation*} 
\end{cor}

In particular, if {{\small $F$}} is also Toeplitz, then {{\small $\alpha=\beta$}} and 
Corollary~\ref{coro-frob-loss} yields
\[
   L_{F}(\hat \si,\si)  \ \ \leq \ \  320 \, \sigma \sqrt{n \log n} \,.
\]
This rate matches the minimax bound for the Toeplitz class established in~\cite{cai2022matrix}, which was derived without computational constraints. To our knowledge, {{\small SABRE}}  is the first polynomial-time method achieving this accuracy.



\subsection{Extension to approximate permutations}
\label{extension-section}

We now relax the assumption that the latent order is an exact permutation.
Exact permutations can impose strong homogeneity constraints in Lipschitz-type models, as consecutive rows of a bi-Lipschitz matrix are almost identical.
This limits the ability to model heterogeneous structures, such as networks with both dense and sparse regions.
To address this, we extend our framework to \textit{approximate permutations}, where latent positions may repeat or leave gaps, provided they remain well spread.

\begin{definition}[Approximate permutation]
\label{spatial-sparsity}
For {{\small $\zeta \ge 0$}}, let {{\small $\ap(\zeta)$}} denote the set of vectors {{\small $\pi\in[n]^n$}} satisfying  
{{\small $\max_{k\in[n]} \min_{i\in[n]} |\pi_i - k| \le \zeta$}}. 
\end{definition}

Thus {{\small $\pi$}} is an approximate permutation with spacing {{\small $\zeta$}}:
when {{\small $\zeta=0$}} we recover the exact case, while larger {{\small $\zeta$}} allows for repeated or missing positions.
A quantitative illustration of this relaxation, showing how it allows for structural heterogeneity in {{\small $F$}}, is given in Appendix~\ref{appendix:homogeneity}.

\smallskip
\textit{Algorithmic modification.} {{\small SABRE}} extends naturally to this setting.
In the refined seriation step (Section~\ref{subsection:stage2}), the balanced tripartition {{\small $(S^1,S^2,S^3)$}} is replaced by leave-one-out subsets {{\small $S^{ij}=[n]\setminus\{i,j\}$}}.
This ensures that comparisons between any two items {{\small $(i,j)$}} use all remaining observations as references,  which mitigates the sparsity effects 
caused by unevenly spaced latent positions (up to the scale {{\small $\zeta$}}). 
The modification increases the computational cost from {{\small $O(n^3)$}} to {{\small $O(n^5)$}}, thereby trading computational efficiency for statistical robustness. Further details and the  pseudocode are given in Appendix~\ref{section-new_appendix-extension-rates}.

\smallskip
\textit{Result.}
Although more expensive computationally, the resulting estimator preserves the same statistical guarantees.
Let {{\small $\mathcal{AL}_\pi(\alpha,\beta,r,r')$}} denote the natural extension of {{\small $\mathcal{AL}(\alpha,\beta,r,r')$}} to approximate permutations; its formal definition is provided in Appendix~\ref{section-new_appendix-extension-rates}.

\begin{thm}
\label{thm-extension}
If {{\small $\pi\in\ap(\zeta)$}} with {{\small $\zeta=o(n)$}} and {{\small $F \in \mathcal{AL}_\pi(\alpha,\beta,r,r')$}}, then, with probability at least {{\small $1-1/n^2$}},
\begin{equation*}
L_{\max}(\hat\pi,\pi) \ \lesssim\ (\va/\alpha)\sqrt{\frac{\log n}{n}} \enspace.
\end{equation*}
\end{thm}

Hence, in the relaxed framework of approximate permutations, {{\small SABRE}} retains the optimal seriation rate established in Theorem~\ref{thm-matrix-setting}.
Unlike previous works~\cite{janssen2020reconstruction,giraud2021localization}, our result does not require the latent positions to lie on a regular grid or follow a uniform sample, thereby extending the applicability of {{\small SABRE}} to heterogeneous data. 
A formal version of this theorem is provided in Appendix~\ref{appendix-extension-approx-rates}.





\section{Discussion}
\label{section-discussion}

\subsection{Summary of contributions}

We derive sharp minimax rates for seriation under the {{\small $L_{\max}$}} loss over the class of average-Lipschitz matrices. 
We also introduce {{\small SABRE}}, a new polynomial-time algorithm that attains these rates with runtime {{\small $O(n^3)$}} and requires no knowledge of the regularity parameters of {{\small $F$}}. 
These results clarify how the statistical difficulty of seriation depends on the regularity of {{\small $F$}} and how structural assumptions shape the attainable rates.  

Within the average-Lipschitz framework, our guarantees also cover the bi-Lipschitz subclass.  
In this setting, {{\small SABRE}} is, to our knowledge, the first polynomial-time procedure attaining the minimax rate {{\small $\sqrt{\log(n)/n}$}}, previously reachable only via exhaustive search~\cite{giraud2021localization}.

The {{\small $L_{\max}$}} bounds yield performance guarantees for other natural losses.
For Kendall's tau, the rate is minimax up to logarithmic factors.
For the Frobenius loss, our bounds specialize to the Toeplitz subclass and attain the minimax rate, originally obtained through exhaustive enumeration~\cite{cai2022matrix}.

Our analysis also extends to more flexible latent structures through approximate permutations, which allow repeated or missing positions.  
In this relaxed framework, {{\small SABRE}} retains the optimal {{\small $L_{\max}$}} rate under mild spacing conditions on the latent positions, illustrating the robustness of the procedure.



\subsection{Discussion of assumptions}

The bi-Lipschitz condition of Definition~\ref{defi:BL} imposes strong pointwise constraints on the latent matrix {{\small $F$}}. 
As shown in Section~\ref{section:examples-BL-AL}, it excludes natural Robinson matrices that remain far from degenerate and do not indicate any fundamental barrier to recovering the latent order.
To address this limitation, we introduce the average-Lipschitz framework, which relaxes the bi-Lipschitz requirements through averaged controls, while retaining sufficient regularity for {{\small $L_{\max}$}}-seriation.

Toeplitz constraints~\cite{cai2022matrix,giraud2021localization,natik2021consistency} impose a homogeneous structure 
by requiring {{\small $F_{ij}=g(|i-j|)$}} for some non-increasing function {{\small $g$}}, so that similarities depend only on index gaps. 
In contrast, the average-Lipschitz condition imposes no such global template: the rows {{\small $F_1,\ldots,F_n$}} are not required to share a common profile.

The key structural requirement for our distance-based seriation is a local correspondence between the ordering gaps in {{\small $\pi$}} and the distance matrix {{\small $D^*$}} in~\eqref{def-dist-[n]}, as formalized in Definition~\ref{cond-LB}.
It is noteworthy that this local condition suffices for our analysis, despite not assuming the typical Euclidean-type distance properties such as additivity or transitivity.
For example, our framework allows {{\small $|\pi_i-\pi_k| = |\pi_i-\pi_j| + |\pi_j-\pi_k|$}} without implying {{\small $D^*_{ik} = D^*_{ij} + D^*_{jk}$}}, and {{\small $|\pi_i-\pi_k| \ge |\pi_j-\pi_k|$}} without entailing {{\small $D^*_{ik} \ge D^*_{jk}$}}.
While such relations may hold approximately in more structured settings (e.g., Toeplitz), where {{\small $D^*$}} inherits regular geometric patterns, they need not hold for general Robinson matrices.
Avoiding these global assumptions allows us to cover heterogeneous cases in which {{\small $D^*$}} lacks Euclidean-type structure.


\subsection{Structural regularity and loss sensitivity in seriation}
\label{section-toeplitz-discussion}

Previous minimax analyses for seriation have been driven by near-degenerate Robinson matrices, whose extreme form determines the lower bounds.
In the Toeplitz class, for instance, the minimax rate is set by corner cases such as matrices that are zero everywhere except along a few diagonals, 
leaving open whether more regular Toeplitz matrices admit faster seriation rates. 
More broadly, it remains unclear how seriation difficulty depends on the structural properties of a Robinson matrix.
Our Lipschitz analysis is a first step toward clarifying this issue and identifying classes of Robinson matrices whose minimax behavior depends on their level of regularity.

We next compare our Lipschitz analysis under {{\small $L_{\max}$}} loss with the Toeplitz framework of~\cite{cai2022matrix}, highlighting the impacts of structural regularity and of the loss function.

\paragraph*{Comparison with Cai et al. (2023).} The work~\cite{cai2022matrix} on Robinson Toeplitz models establishes minimax rates for exact reordering under Frobenius loss.  
In particular,~\cite{cai2022matrix} introduces the signal-to-noise ratio
\[
m(\mathcal{F},\mathcal{S}) = \min_{F\in\mathcal{F},\,\pi\neq\pi'\in\mathcal{S}} \|F_{\pi}-F_{\pi'}\|_F ,
\]
where {{\small $\mathcal{F}$}} denotes a class of Toeplitz matrices and {{\small $\mathcal{S}$}} a set of permutations,
and shows, via an exhaustive search over permutations, that exact recovery is possible if, and only if
\begin{equation}
\label{eq:tony-vai-snr}
    m(\mathcal{F},\mathcal{S}) \gtrsim \sigma\sqrt{n\log n}.
\end{equation}

Although this provides a valuable benchmark, our results improve upon this framework in two  ways.
First, when matrices are Toeplitz and bi-Lipschitz, Corollary~\ref{coro-frob-loss} shows that {{\small SABRE}} achieves exact recovery at the minimax SNR level~\eqref{eq:tony-vai-snr} in polynomial time.  This closes the statistical-computational gap for a broad class of structured matrices.

Second, the threshold~\eqref{eq:tony-vai-snr} characterizes statistical indistinguishability under Frobenius loss, 
which is driven by global perturbations. 
In contrast, {{\small $L_{\max}$}} remains sensitive to local variations that may be invisible in Frobenius norm.
As a result, accurate seriation may remain feasible even when the Frobenius-based SNR criterion is not met.
We illustrate this with a simple Toeplitz bi-Lipschitz matrix for which {{\small SABRE}} achieves exact recovery despite violating the SNR condition~\eqref{eq:tony-vai-snr}.

\smallskip 
\textit{Example.} 
Consider {{\small $F$}} with entries {{\small $F_{ij} = 1 - |i-j|/n$}},
and the parameter space   {{\small $(\mathcal{F},\mathcal{S}) = (\{F\},\{\mathrm{id},(1,k)\})$}}, where {{\small $(1,k)$}} swaps the first and {{\small $k$}}-th indices.  
Set {{\small $k = 41 \sigma \sqrt{n\log n}$}}. 
The SNR condition~\eqref{eq:tony-vai-snr} is violated by a factor of order {{\small $\sqrt{n}$}},
since {{\small $m(\mathcal{F},\mathcal{S}) = \|F - F_{(1,k)}\|_F \lesssim \sigma\sqrt{\log n}$}}.
Yet,  {{\small SABRE}} achieves exact recovery (Theorem~\ref{thm-matrix-setting}). 




\subsection{Related work on the maximum error}

The {{\small $L_{\max}$}} loss has been investigated in latent-space models such as the bi-Lipschitz framework of~\cite{giraud2021localization} and the graphon model in~\cite{janssen2020reconstruction,natik2021consistency}.  

\cite{giraud2021localization} provides the first optimal-rate guarantees for the {{\small $L_{\max}$}} loss in a bi-Lipschitz setting, 
but their procedure covering the full bi-Lipschitz class relies on an exhaustive search over permutations.  
In the same work, the authors also study a spectral method as a polynomial-time alternative, but its guarantees hold only under additional structure, namely approximately Toeplitz matrices with a sufficiently large spectral gap.

In graphon models, the goal is to recover an ordering of nodes from the noisy adjacency matrix of a network.  
\cite{janssen2020reconstruction} proposes a thresholding estimator based on the squared adjacency matrix, and obtains an error  
{{\small $L_{\max}(\hat\pi,\pi)\lesssim (\log n)^5/\sqrt{n}$}} with high probability under regularly spaced latent positions and technical assumptions on the graphon.  
\cite{natik2021consistency} analyzes a spectral algorithm under {{\small $\mathcal{C}^1$}}-smoothness and monotonicity conditions on the graphon.  
These results rely on latent-space graphon models, which differ from the Lipschitz setting considered here.


\subsection{Limitations and future directions}

While {{\small SABRE}} attains minimax-optimal rates under our assumptions, its implementation requires selecting tuning parameters and specifying a threshold {{\small $\sigma$}} that may be any upper bound on the noise level.
We provide an explicit choice of tuning parameters depending only on the sample size {{\small $n$}}.
However, our guarantees for this choice hold only when {{\small $n$}} exceeds a constant {{\small $C_{\alpha\beta r r'\sigma}$}} determined by the underlying model parameters.
Although such a mild dependence is standard in minimax analyses, it still constitutes a limitation for finite-sample applicability. 
A sharper analysis of parameter calibration and finite-sample robustness remains an open direction.

Further relaxing structural assumptions is of interest.  
Approximate permutations (Section~\ref{extension-section}) already cover piecewise monotone patterns and small local blocks, but hybrid settings combining ordering with clustering structures are  natural next steps.

Exploring faster algorithms is another relevant direction.  
While {{\small SABRE}} runs in cubic time, additional structure in {{\small $D^*$}} (e.g., Toeplitz or nearly Toeplitz forms) could reduce the computational cost.  
Alternative similarity measures beyond {{\small $D^*$}} may be worth investigating, since our analysis does \textit{not} rely on classical distance properties such as additivity, transitivity, or the triangle inequality.

Because {{\small SABRE}} targets the stringent {{\small $L_{\max}$}} criterion, some of its machinery may be unnecessary  when the loss is weaker, suggesting that simpler or faster procedures could suffice for Frobenius or Kendall losses.




\section{Proof sketch}
\label{section-algo-analysis}

We now present the arguments underlying Theorem~\ref{thm-matrix-setting}.

\subsection{Informal derivation of the rate}
\label{informal-derivation-of-rates}

We first outline how the three stages of {{\small SABRE}} combine to yield the rate stated in Theorem~\ref{thm-matrix-setting}. 
The corresponding guarantees for each stage are stated in the subsequent subsections.

\paragraph*{Stages~1-2. Distance estimation and coarse seriation.}
Under the average-Lipschitz condition, the population distance matrix {{\small $D^*$}} satisfies the local distance-equivalence property (Definition \ref{cond-LB}), linking distance values to ordering gaps {{\small $|\pi_i-\pi_j|$}}. 
Proposition \ref{small-prp-dist-estim-S=n} further shows that, with high probability, the empirical estimator {{\small $\widehat D$}} preserves this structure up to an additive deviation {{\small $\omega_n\!\asymp_{\beta,\sigma}\!n^{3/4}(\log n)^{1/4}$}}. 
Hence, {{\small $\widehat D$}} carries reliable local information about the latent order {{\small $\pi$}}, 
which will be exploited in the first seriation step.

Since {{\small $\widehat D$}} is informative only for pairs whose true distances satisfy {{\small $D^*_{ij} \gtrsim \omega_n$}}, the first seriation step  operates at this resolution, choosing the tuning parameters {{\small $(\delta_1,\delta_2,\delta_3)$}} on the same order.

Through the local distance-equivalence property, this distance scale {{\small $D^*_{ij}\!\gtrsim\!\omega_n$}} corresponds to ordering gaps {{\small $|\pi_i-\pi_j|\!\gtrsim\!\omega_n/\alpha$}} in the latent permutation.  
Hence, the first seriation step can at best recover comparisons above this ordering resolution, that is, at the finest scale permitted by the accuracy of {{\small $\widehat D$}}.
Proposition \ref{prop:end-stage1-new} formalizes this intuition, showing that, under this calibration of the tuning parameters, the resulting matrix {{\small $H$}} coincides with {{\small $H^*$}} for all such pairs (those with ordering gaps of order at least {{\small $\omega_n/\alpha$}}). This coarse estimate serves as input to the final refinement stage.

\paragraph*{Stage~3. Refined seriation.}
The final stage re-examines the undecided pairs using aggregate statistics of the form
{{\small $l=\sum_{k\in L}(A_{ik}-A_{jk})$}},
where the index set {{\small $L$}} -- inferred from the coarse ordering of Stage~2 -- is, with high probability, located on one side of both {{\small $i$}} and {{\small $j$}} in the latent order.
If one ignores the dependence between the statistic {{\small $l$}} and the set {{\small $L$}}, 
we obtain the following simple heuristic. 
Under the average-Lipschitz separation condition, the expected signal satisfies 
{{\small $\E[l]\ge \alpha|\pi_i-\pi_j|$}}, 
while the stochastic fluctuations of the noise are upper bounded by 
{{\small $O(\sigma\sqrt{n\log n})$}} with high probability. 
Hence, whenever {{\small $|\pi_i - \pi_j|\gtrsim(\sigma/\alpha)\sqrt{n\log n}$}}, 
the signal dominates the noise and the comparison between {{\small $i$}} and {{\small $j$}} is recovered correctly.
Aggregating these guarantees over all pairs yields the rate 
{{\small $L_{\max}(\hat\pi,\pi)\lesssim(\sigma/\alpha)\sqrt{\log(n)/n}$}}, as stated in Theorem~\ref{thm-matrix-setting}.

A key difficulty at this stage is the statistical dependence between {{\small $l$}} and the data used to form the set {{\small $L$}}.  
To make the above argument rigorous, this dependence is mitigated through a light sample-splitting step that restores approximate independence between the statistic and the underlying data, allowing standard concentration arguments to apply.
Proposition~\ref{prop:complete-seriation} formalizes this analysis and shows that, under this construction, all pairs separated by at least {{\small $O((\sigma/\alpha)\sqrt{n\log n})$}} are ordered correctly with high probability.

\paragraph*{Choice of tuning parameters.} 
The tuning parameters used in Theorem~\ref{thm-matrix-setting},
namely {{\small $\delta_1 = n^{3/4}\log n$}} and 
{{\small $\delta_{k+1} = \delta_k \log n$}}, 
satisfy the admissibility constraints from 
Sections~\ref{subsection-first-seriation-analysis} and 
\ref{section:refined-seriation-analysis-mainpaper}.  
A full verification is provided in Appendix~\ref{appendix-general-theorem}.


\subsection{Distance estimation}
\label{sub-section-distance-analysis}

We analyze the estimation of the population distance matrix {{\small $D^*$}}, defined in~\eqref{def-dist-[n]}. 
This analysis relies on a local distance-gap property stated below.

\begin{definition}[Local distance equivalence]
\label{cond-LB}
For  {{\small $\alpha,\beta,r>0$}} and {{\small $\omega\ge 0$}}, 
let {{\small $\lde(\alpha,\beta,\omega,r)$}} denote the set of matrices 
{{\small $D\in\R^{n\times n}$}} such that, for all {{\small $i,j$}} satisfying 
{{\small $|\pi_i-\pi_j|\wedge D_{ij} \le nr$}},
\[  \alpha|\pi_i-\pi_j|-\omega
  \ \ \le \ D_{ij} \ \le \ \
  \beta|\pi_i-\pi_j|+\omega \; .\]
\end{definition}

Under the average-Lipschitz condition, the population distance matrix {{\small $D^*$}} can be shown to satisfy this local distance equivalence.
In addition, the empirical estimator {{\small $\widehat D$}}, introduced in Section~\ref{subsection:l2-dist}, 
admits a uniform concentration bound around {{\small $D^*$}}, 
with deviations of order {{\small $\omega_n\asymp_{\beta,\sigma}n^{3/4}(\log n)^{1/4}$}}.   
Combining these two properties yields the following result.

\begin{prop}
\label{small-prp-dist-estim-S=n}
For any {{\small $(\alpha,\beta,r,r',\sigma)$}}, there exists a constant {{\small $C_{\beta \sigma}$}} depending only on {{\small $(\beta ,\sigma)$}} such that the following holds for all {{\small $n \geq 8 \vee \tfrac{1}{r}$}}.
If {{\small $F\in\mathcal{AL}(\alpha,\beta,r,r')$}}, then with probability at least {{\small $1-1/n^4$}}, the estimator {{\small $\widehat D$}} satisfies the uniform bound
\begin{equation}
    \label{distance-estimate-simple-error-constant}
    \max_{i,j} |\widehat D_{ij} - D^*_{ij}| \ \leq \ \omega_n,
    \qquad \quad  
    \omega_n := C_{\beta \sigma}\, n^{3/4} (\log n)^{1/4},
\end{equation}
and moreover {{\small $\widehat D \ \in \  \lde(\alpha,\beta,\omega_n,r \wedge r')$}}.
\end{prop}

This proposition shows that, for all locally close pairs, 
the empirical distances {{\small $\widehat D_{ij}$}} scale approximately linearly with the ordering gaps {{\small $|\pi_i-\pi_j|$}},
up to multiplicative factors {{\small $\alpha$}} and {{\small $\beta$}}, and an additive deviation {{\small $\omega_n$}}.
The proof is deferred to Appendix~\ref{appendix-distance-estimation}.




\subsection{First seriation step}
\label{subsection-first-seriation-analysis}

We analyze the first seriation step of Section~\ref{section:idea-algo}.  
The next result gives a deterministic guarantee, valid for any input matrix {{\small $D\in\lde(\alpha,\beta,\omega,r)$}} (Definition~\ref{cond-LB}).

\begin{prop}
\label{prop:end-stage1-new} 
For any  {{\small $(\alpha,\beta, r)$}}, 
there exist constants {{\small $C_{\beta}$, $C_{\alpha\beta}$, $C_{\alpha r}$,  $C'_{\alpha\beta}$, $C'_{\alpha r}>0$}} 
depending only on {{\small $(\alpha,\beta,r)$}} such that the following holds for all {{\small $n$}} and {{\small $\omega$}}. 
If the input {{\small $D$}} satisfies {{\small $D\in\lde(\alpha,\beta,\omega,r)$}}, 
and the tuning parameters {{\small $(\delta_1,\delta_2,\delta_3)$}} fulfill
\begin{equation}
\label{delta-constraint-1rst-seriation-simùplified}
C_{\beta}\; \om \le \delta_1, \qquad
C_{\alpha\beta}\,\delta_1 \le \delta_2 \le C_{\alpha r}\,n, \qquad
C'_{\alpha \beta}\,\delta_2 \le \delta_3 \le C'_{\alpha r}\,n ,
\end{equation}
then the output {{\small $H$}} of the first seriation step satisfies, for some {{\small $s \in \{\pm\}$}},
\[
H_{ij} \ = \ s H^*_{ij} \,, \qquad 
\forall i,j \in [n] \text{ such that } H_{ij}\neq 0 \ \text{ or } \ |\pi_i-\pi_j|\geq (\delta_2+\omega)/\alpha .
\]
\end{prop}

The sign {{\small $s\in\{\pm\}$}} accounts for the natural symmetry of the seriation problem,
since the latent ordering is identifiable only up to a global reversal.

This proposition ensures that the matrix {{\small $H$}} produced by the first seriation step is correct on its support and for all pairs with ordering gaps at least {{\small $(\delta_2+\omega)/\alpha$}}, 
a threshold that increases with the additive slack {{\small $\omega$}} and decreases with the signal level  {{\small $\alpha$}} from the condition {{\small $D\in\lde(\alpha,\beta,\omega,r)$}}.
It therefore characterizes the coarse recovery properties of the first seriation stage and provides the foundation for the refined analysis developed in the next section.
The proof of Proposition~\ref{prop:end-stage1-new} is given in Appendix~\ref{appendix-aggregation-bisections}.

\begin{remark}
Informally, the constraints~\eqref{delta-constraint-1rst-seriation-simùplified} balance two competing requirements: 
the thresholds {{\small $(\de_1, \de_2, \de_3)$}} must be large enough compared to the deviation parameter {{\small $\omega$}}  
to ensure reliable comparisons, yet not so large that the first seriation step loses resolution 
and fails to discriminate nearby items. A detailed version of~\eqref{delta-constraint-1rst-seriation-simùplified} is given in Appendix~\ref{appendix-general-theorem}.
\end{remark}


\subsection{Refined seriation step}
\label{section:refined-seriation-analysis-mainpaper}

The refined seriation step of Section~\ref{subsection:stage2} uses the coarse ordering from Stages~1-2 to determine the remaining comparisons.  
The proposition below quantifies the accuracy of this refinement, when the inputs {{\small $(D,H)$}} are sufficiently accurate approximations of {{\small $(D^*,H^*)$}}.
Recall that {{\small $\omega_n = C_{\beta \sigma} n^{3/4}(\log n)^{1/4}$}} from~\eqref{distance-estimate-simple-error-constant}.
For convenience, we set {{\small $r=r'$}}, with the case {{\small $r\neq r'$}} obtained by replacing {{\small r}} with {{\small $r\wedge r'$}} where needed.

\begin{prop}
\label{prop:complete-seriation}
For any {{\small $(\alpha,\beta,r,\sigma)$}}, there exist constants 
{{\small $\tilde C_{\beta}$, $\tilde C_{\alpha\beta}$, $\tilde C_{\alpha r}$, $C_{\beta \sigma}>0$}} 
depending only on {{\small $(\alpha,\beta,r,\sigma)$}} such that the following holds for all {{\small $n$}}.
Assume that {{\small $F\in\mathcal{AL}(\alpha,\beta,r,r)$}},
and that the input matrices {{\small $(D, H)$}} satisfy {{\small $D \in \lde(\alpha,\beta,\omega_n,r)$}} and 
\begin{equation}\label{prop:input-H}
H_{ij} = s H^*_{ij}, 
\qquad 
\forall i,j \ \text{ such that } \ H_{ij}\neq 0 \ \text{ or } \ |\pi_i-\pi_j| \geq \rho ,
\end{equation}
for some {{\small $s\in\{\pm\}$}} and {{\small $\rho \in [0, rn]$}}.
If  the tuning parameter {{\small $\delta_4$}} fulfills
\begin{equation}
\label{constraint-2nd-seriation-simplified-new}
 \tilde C_{ \beta }\; \rho +   \tilde  C_{\alpha \beta}\; (\om_n \vee \sqrt{n\log n})\;  \ \leq \ \delta_4 \ \leq \ \tilde C_{\alpha r} n ,
\end{equation}
then the output of the final seriation step satisfies
\begin{equation*}
\P\Big\{ \widehat H_{ij} = s H^*_{ij},\ 
\forall i,j:\ |\pi_i-\pi_j|\ge\tfrac{40\sigma}{\alpha}\sqrt{n\log n} \Big\}
\ge 1-\tfrac{9}{n^3}.
\end{equation*}
\end{prop}

Condition~\eqref{prop:input-H} ensures that the input {{\small $H$}}  
already coincides with {{\small $H^*$}} both on its support and for sufficiently separated pairs.  
Building on this, the final seriation step refines these coarse guarantees, 
yielding correct pairwise comparisons for all indices whose ordering distance exceeds 
{{\small $O((\sigma/\alpha)\sqrt{n\log n})$}}.
The proof of Proposition~\ref{prop:complete-seriation} is given in Appendix~\ref{appendix:proof-complete-seriation}.

\begin{remark}
In essence, condition~\eqref{constraint-2nd-seriation-simplified-new} 
requires {{\small $\delta_4$}} to be large enough to dominate the cumulative errors from the first seriation stage, 
yet not so large that the refinement step loses resolution.  
A detailed version of~\eqref{constraint-2nd-seriation-simplified-new}, with explicit constants, 
is given in Appendix~\ref{appendix-general-theorem}.
\end{remark}




\bigskip


\begin{funding}
 The work of C. Giraud has been partially supported  by grant ANR-19-CHIA-0021-01 (BiSCottE, ANR) and ANR-21-CE23-0035 (ASCAI, ANR).
The work of  N. Verzelen has been partially supported by ANR-21-CE23-0035 (ASCAI, ANR). 
\end{funding}

\bigskip 


\bibliographystyle{imsart-number} 
\bibliography{biblio_seriation.bib}       


\newpage

\begin{appendix}

\section{Tuning parameters and general version of Theorem~\ref{thm-matrix-setting}}
\label{appendix-general-theorem}

Theorem~\ref{thm-matrix-setting} was first stated under the simplifying assumption that the model parameters {{\small $(\alpha,\beta,r,r',\sigma)$}} are fixed constants independent of {{\small $n$}}.  
This made it possible to provide explicit tuning values {{\small $(\delta_1,\delta_2,\delta_3,\delta_4)$}} that do not depend on unknown quantities such as {{\small $(\alpha,\beta,r,r')$}}.  
In fact, the result extends to the general situation where these parameters may depend on {{\small $n$}}.

We now collect the technical requirements on the tuning parameters 
{{\small $(\delta_1,\delta_2,\delta_3,\delta_4)$}} 
used in the analyses of the first and refined seriation stages of Sections~\ref{subsection-first-seriation-analysis}
and~\ref{section:refined-seriation-analysis-mainpaper}.  
For simplicity of exposition, we present these calibration constraints in the case {{\small $r=r'$}}, 
and keep the symbol {{\small $r$}} throughout 
{{\small(\ref{delta-constraint-1rst-seriation}--\ref{constraint-2nd-seriation})}}.  
All statements remain valid for distinct radii {{\small $r$}} and {{\small $r'$}}, 
in the worst-case sense that any occurrence of {{\small $r$}} may be replaced by 
{{\small $r_0:=r\wedge r'$}} wherever needed.

\medskip

\textit{Tuning parameters for the first seriation step.} The tuning parameters {{\small $(\delta_1,\delta_2,\delta_3)$}} are required to satisfy the following condition, which is a detailed version of~\eqref{delta-constraint-1rst-seriation-simùplified} with explicit constants:  
\begin{align}\label{delta-constraint-1rst-seriation}
 \omega + \beta \ \leq \ \delta_1 , \qquad 
 \omega + \frac{\beta}{\alpha}(\delta_1+\omega) \ < \ \delta_2 , \qquad 
 2 \ \leq \ \frac{\delta_2+\omega}{\alpha} \ \leq \ \frac{n}{8} , \\
 1 \vee (\delta_2+\omega) \ \leq \ (1 \wedge \alpha) r n , \qquad 
 \omega + \frac{\beta}{\alpha}(\delta_2+\omega) \ < \ \delta_3 \ \leq \ \big(r \wedge \tfrac{\alpha}{8}\big)n - \omega .
 \nonumber
\end{align}
Here {{\small $\omega$}} denotes the deviation parameter appearing in the local 
distance-equivalence condition {{\small $D\in\lde(\alpha,\beta,\omega,r)$}} 
(Definition~\ref{cond-LB}).  
In the statistical analysis of {{\small SABRE}}, it will be instantiated as 
{{\small $\omega=\omega_n$}} from Proposition~\ref{small-prp-dist-estim-S=n}.


\smallskip 
\textit{Tuning parameters for the second seriation step.} In the refined stage, the algorithm takes as input an intermediate matrix {{\small $H$}} 
obtained from the first seriation step, whose accuracy is characterized by 
{{\small $\rho = (\delta_2+\omega)/\alpha$}}, as established in Proposition~\ref{prop:end-stage1-new}. 
The new threshold {{\small $\delta_4$}} introduced in this step is then required, 
together with {{\small $\rho$}}, to satisfy the following condition, which is a detailed version of~\eqref{constraint-2nd-seriation-simplified-new} with explicit constants:
\begin{align}
\label{constraint-2nd-seriation}
\delta_4 + 4\beta \sqrt{n\log n} + 4\omega_n &< \alpha n ( r \wedge \tfrac{1}{32}) , \nonumber \\
   2\omega_n  + \beta\Big( 2\rho + \alpha^{-1}(4\beta \sqrt{n\log n} + 2\omega_n )\Big) &\leq  \delta_4 . 
\end{align}

We now state the general version of Theorem~\ref{thm-matrix-setting}, where the model parameters may depend on {{\small $n$}}.
\begin{thm}
\label{thm-matrix-setting-general}
For any {{\small $n\ge C$}} and {{\small $(\alpha,\beta,r,r',\sigma)$}}, 
where {{\small $C>0$}} is a  numerical constant,
if {{\small $F\in\mathcal{AL}(\alpha,\beta,r,r')$}}, then,  under the  conditions~(\ref{delta-constraint-1rst-seriation}--\ref{constraint-2nd-seriation}),   the conclusion of Theorem~\ref{thm-matrix-setting} still holds.
\end{thm}

Thus, {{\small SABRE}} achieves the same convergence rate even when the model parameters vary with {{\small $n$}}, provided the tuning parameters obey the constraints~(\ref{delta-constraint-1rst-seriation}--\ref{constraint-2nd-seriation}).  These conditions are rather intricate; to make them more transparent, we summarized in \eqref{delta-constraint-1rst-seriation-simùplified} and \eqref{constraint-2nd-seriation-simplified-new}  how they reduce in the constant-parameter case. Below we explain how they lead to the choice of tuning values made in Theorem~\ref{thm-matrix-setting}.


\subsection{Choice of tuning parameters.} 
\label{subsection-choice-tuning-param-from-simplified-condition}
When the model parameters {{\small $(\alpha,\beta,r,r',\sigma)$}} are assumed to be fixed constants independent of {{\small $n$}}, the  constraints (\ref{delta-constraint-1rst-seriation}--\ref{constraint-2nd-seriation}) admit  explicit  solutions. 
Setting {{\small $\omega=\omega_n$}} using the distance-estimation rate {{\small $\omega_n \asymp_{\beta,\sigma} n^{3/4}(\log n)^{1/4}$}} (see~\eqref{distance-estimate-simple-error-constant}), the constraints~(\ref{delta-constraint-1rst-seriation}--\ref{constraint-2nd-seriation}) are satisfied for any choice of {{\small $(\delta_1,\delta_2,\delta_3,\delta_4)$}} such that
\[
C_{\beta\sigma}\,n^{3/4}(\log n)^{1/4} \le \delta_1, \qquad
C_{\alpha\beta}\,\delta_1 \le \delta_2 \le C_{\alpha r}\,n, \qquad
C'_{\alpha,\beta}\,\delta_2 \le \delta_3 \le C'_{\alpha r}\,n,
\]
and, for {{\small $\rho=(\delta_2+\omega_n)/\alpha$}},
\[
\tilde C_{\alpha\beta}\,\delta_2 \ \le \ \delta_4 \ \le \ \tilde C_{\alpha r} n,
\]
for some constants {{\small $C_{\beta\sigma},
C_{\alpha\beta}, C_{\alpha r},
C'_{\alpha,\beta}, C'_{\alpha r}, \tilde C_{\alpha\beta}, \tilde C_{\alpha r} > 0$}}.
These conditions are fulfilled for all sufficiently large {{\small $n$}} provided that the following relations hold:
\begin{equation}\label{asymptotic-inputs}
\frac{\delta_1}{n^{3/4}(\log n)^{1/4}} \to \infty, 
\qquad 
\frac{\delta_{k+1}}{\delta_k} \to \infty \ \ (k=1,2,3), 
\qquad 
\frac{\delta_4}{n} \to 0.
\end{equation}
Hence, a convenient and admissible  calibration is {{\small $\delta_1 = n^{3/4}\log n$}} and {{\small $\delta_{k+1}=\delta_k \log n$}} for {{\small $k\in[3]$}}.
This choice is used in the statement of Theorem~\ref{thm-matrix-setting}.




\section{Pseudocode}
\label{appendix:pseudocode}

This section provides the full pseudocode of {{\small SABRE}}, following the three stages introduced in Section~\ref{section-sabre}. We begin with a schematic overview summarizing the inputs and outputs of each stage, to guide navigation through the detailed procedures. The pseudocode for each stage is given in the subsections that follow.

\smallskip 

\paragraph*{Schematic overview of {{\small SABRE}}}
{{\small
\begin{algorithmic}[1]
\REQUIRE Matrix {{\small $A\in\R^{n\times n}$}}; upper bound {{\small $\sigma$}} from~\eqref{eq:model}; tuning parameters {{\small $(\delta_1,\delta_2,\delta_3,\delta_4)$}}.
\ENSURE {{\small $\hat\pi\in[n]^n$}}.
\STATE Construct a distance estimate {{\small $\widehat D$}} from {{\small $A$}}.
\STATE Build a partial comparison matrix {{\small $H$}} from {{\small $\widehat D$}} using {{\small $(\delta_1,\delta_2,\delta_3)$}}.
\STATE Refine the undecided comparisons using {{\small $H$}}, {{\small $A$}}, and  the parameters {{\small $(\sigma, \delta_4)$}} to obtain {{\small $\widehat H$}}.
\STATE Derive {{\small $\hat\pi$}} from the row sums of {{\small $\widehat H$}}.
\end{algorithmic}
}}

\smallskip 

We now provide the full pseudocode.

\begin{algorithm}[H]
\caption{ {{ \small \texttt{SABRE}$(A,\delta_1,\delta_2,\delta_3,\delta_4,\sigma)$}} }
\label{algo:sabre}
\begin{algorithmic}[1]

  \REQUIRE $(A, \sigma, \delta_1,\delta_2,\delta_3,\delta_4)$
  \ENSURE  $\hat\pi\in[n]^n$

  \STATE $\widehat D \leftarrow$ \texttt{Estimate} \texttt{Distance}$(A,[n])$

  \STATE $H \leftarrow$ \texttt{Aggregate} \texttt{Bisections}$(\widehat D,\delta_1,\delta_2,\delta_3)$

  \STATE $\widehat H \leftarrow$ \texttt{Re-evaluate} \texttt{Comparisons}$(H,\widehat D,A,\sigma,\delta_4)$

  \STATE $\hat\pi_i \leftarrow \big(\widehat H_i \mathbf{1}_n + n + 1\big)/2$ for all $i$

\end{algorithmic}
\end{algorithm}



\subsection{Implementation  of the distance estimator}
\label{distance-estimator-pseudocode}

This appendix provides the full pseudocode for the distance estimator used in Stage~1 of {{\small SABRE}}.  
The procedure is written for a general subset {{\small $S\subset[n]$}}, as this form is used in the data-splitting scheme of Sections~3.3.

\begin{algorithm}[H]
\caption{ {{\small \texttt{Estimate} \texttt{Distance}}}}
\label{algo:distance}
\begin{algorithmic}[1]
    \REQUIRE Symmetric matrix  $A\in\R^{n\times n}$ and subset  $S\subset[n]$
    \ENSURE Symmetric matrix $\widehat D=\widehat D(S)$

    \STATE Build the restricted matrix  $A^S$:
    \[
        A^S_{ij} = 
        \begin{cases}
            A_{ij}, & (i,j)\in S\times S,\\[1mm]
            0,      & \text{otherwise}.
        \end{cases}
    \]

    \FOR{ $i\in S$}
        \STATE 
        $\widehat m_i \;=\; \argmin_{\,j\in S,\ j\neq i}\;
        \max_{\,k\in S\setminus\{i,j\}}\,
        |\langle A^S_k,\, A^S_i - A^S_j\rangle|$
    \ENDFOR

    \FOR{$i,j\in S$ \textbf{with}  $i<j$}
        \STATE 
        $(\widehat D_{ij})^2
        \;=\;
        n\Big(\,
            \langle A^S_i, A^S_{\widehat m_i}\rangle
            +
            \langle A^S_j, A^S_{\widehat m_j}\rangle
            -
            2\langle A^S_i, A^S_j\rangle
        \Big)$
        \STATE $\widehat D_{ji}=\widehat D_{ij}$
    \ENDFOR

    \STATE Set  $\widehat D_{ij}=0$ whenever $i \notin S$ or $j \notin S$
\end{algorithmic}
\end{algorithm}

When {{\small $S=[n]$}}, the procedure returns the estimator {{\small $\widehat D$}} used in Stage~1.


\subsection{Implementation  of the first-seriation stage}
\label{appendix-algo:orientation}

We provide below the pseudocode for the first-seriation stage of {{\small SABRE}},  which constructs the partial comparison matrix {{\small $H$}} from a generic distance estimate {{\small $D$}}. 
The routine {{\small \texttt{Aggregate} \texttt{Bisections}}} implements precisely this procedure, as described in
Section~\ref{section:idea-algo}, and invokes the subroutine {{\small \texttt{Orientation}}}  to assign consistent left/right labels across indices. We write {{\small $0_{n\times n}$}} for the {{\small $n\times n$}} zero matrix.

\begin{algorithm}[H]
    \caption{{{\small \texttt{Aggregate Bisections}}}} \label{algo:SAB} 
    \begin{algorithmic}[1]
        \REQUIRE $(D, \delta_1,\delta_2,\delta_3)$, where $D$ is an $n\times n$ symmetric matrix
        \ENSURE $(L_i,R_i)_{i\in[n]}$ and $H \in \{-1,0,1\}^{n\times n}$
        \FOR {$i \in [n]$}
            \STATE Build a graph $\mathcal G_i$ with node set $[n]\setminus\{i\}$ by linking $k,l$ whenever $D_{kl}\leq \delta_1$ and $D_{ik}\vee D_{il}\geq \delta_2$.
            \STATE Collect all connected components of $\mathcal G_i$ that contain a vertex  $k$ with $D_{ik}\geq \delta_3$.
            \STATE Let $G_i,G'_i$ be the two largest such components.
        \ENDFOR
        \STATE $(L_i,R_i)_{i\in[n]} \leftarrow$ \texttt{Orientation}$((G_i,G'_i)_{i\in[n]})$ \label{line:orientation}
        \STATE Initialize $H$ as the $n\times n$ null matrix: $H=0_{n\times n}$
        \FOR{$i,j \in [n]$ with  $i\neq j$ and $(L_i,R_i)\neq (\emptyset,\emptyset)$}
            \STATE If $i\in L_j$ or $j\in R_i$, then $H_{ij}=-1$
            \STATE Else if $i\in R_j$ or $j\in L_i$, then $H_{ij}=1$
        \ENDFOR
    \end{algorithmic}
\end{algorithm}

Given a set {{\small $G$}}, {{\small $|G|$}} denotes its cardinality.

{{\small\underline{\textbf{Subroutine} \texttt{Orientation}}}}
\renewcommand{\algorithmiccomment}[1]{#1}
{{\footnotesize
    \begin{algorithmic}[1]
        \REQUIRE $(G_i,G_i')_{i \in [n]}$
        \ENSURE $(L_i,R_i)_{i\in [n]}$
        \STATE \label{pick-c-orientation-line} Pick any index $c\in \argmax_{i\in[n]} \left(|G_i| \wedge |G_i'|\right)$, and set  $L_c = G_c$ and $R_c = G_c'$. 
        \FOR{$i\in[n]$ \ such that  $i \neq c $\ and $G_i' = \emptyset$}
            \IF{ $i\in L_c$} 
                \STATE{$L_i=\emptyset$, \quad   $R_i=G_i$}  \label{line-i-one-non-empty}
            \ELSE
                \STATE{$R_i=\emptyset$, \quad   $L_i =G_i$} \label{line-i-one-non-empty-bis}
            \ENDIF 
        \ENDFOR
         \FOR{$i\in[n]$ \ such that  $i \neq c $\ and $G_i' \neq \emptyset$}
            \IF{$G_i\cap L_c=\emptyset$  \ \textbf{or}\ $G_i' \cap R_c=\emptyset$}
                \STATE   $R_i=G_i$, \quad    $L_i=G_i'$  
            \ELSE
                \STATE  $L_i=G_i$, \quad    $R_i=G_i'$ \label{line-i-two-non-empty-bis}
            \ENDIF 
        \ENDFOR
    \end{algorithmic}
    }}

In line~\ref{pick-c-orientation-line}, we pick an index {{\small $c$}} expected to lie near the middle of the ordering (i.e., {{\small $\pi_c\simeq n/2$}}). We arbitrarily fix the orientation of {{\small $(L_c,R_c)$}}, then propagate it to all {{\small $i\neq c$}} so that {{\small $(L_i,R_i)$}} shares the same orientation. As shown later, {{\small $G_i\neq\emptyset$}} for all {{\small $i$}}, so only two cases occur: {{\small $G_i'=\emptyset$}} (line~$2$) or {{\small $G_i'\neq\emptyset$}} (line~$9$).


\subsection{Implementation of the refined-seriation stage}
\label{appendix:pseudocode:refined-seriation}

This appendix provides the pseudocode and notation  for Stage~3 of {{\small SABRE}}.  
The implementation consists of two components, as described in Section~\ref{subsection:stage2}:
(i) a comparison subroutine evaluating the statistics associated with a pair 
{{\small $(i,j)$}}, and 
(ii) the construction of the modified comparison sets {{\small $\tilde L_{ij},\tilde R_{ij}$}} obtained through a light sample-splitting scheme.

\paragraph*{Comparison subroutine.}
Given two rows {{\small $A_i,A_j$}} and sets {{\small $L,R$}}, the
routine computes
\[
l=\sum_{k\in L}(A_{ik}-A_{jk}),
\qquad
r=\sum_{k\in R}(A_{ik}-A_{jk}),
\]
and applies the rule below.

{{\small\underline{\textbf{Subroutine} \texttt{Evaluate Comparison}}}}
{{\footnotesize
\begin{algorithmic}[1]
    \REQUIRE $(A_i, A_j, L, R, \sigma)$ 
    \ENSURE $H_{ij} \in \{-1,0,1\}$
    \STATE $H_{ij}=0$
    \STATE $l=\sum_{k\in L}(A_{ik}-A_{jk})$, \qquad $r=\sum_{k\in R}(A_{ik}-A_{jk})$
    \IF{$|l|\ge 5\sigma\sqrt{n\log n}$}
        \STATE $H_{ij}=-\mathrm{sign}(l)$
    \ELSIF{$|r|\ge 5\sigma\sqrt{n\log n}$}
        \STATE $H_{ij}= \mathrm{sign}(r)$
    \ENDIF
\end{algorithmic}
}}

\smallskip 

\paragraph*{Construction of the modified comparison sets.}
To reduce the dependence between {{\small $(A_i,A_j)$}} and the preliminary sets
{{\small $L_{ij},R_{ij}$}}, we recompute distance estimates on submatrices of {{\small $A$}} that exclude rows {{\small $i,j$}}.
We draw a random balanced partition
\[
[n]=S^1\,\cup\,S^2\,\cup\,S^3,
\qquad
\lfloor n/3\rfloor \le |S^t|\le \lceil n/3\rceil,
\]
and split the data matrix {{\small $A$}} into blocks {{\small $(A^t)_{t\in[3]}$}}, where each block {{\small $A^t$}} denotes the submatrix of {{\small $A$}} supported on {{\small $S^t\times S^t$}}. 
From each block we compute a distance estimate {{\small $\widehat D^t$}} using the Algorithm~\ref{algo:distance},
\[
\widehat D^t=\texttt{Estimate} \ \texttt{Distance}(A,S^t),\qquad t=1,2,3.
\]
For each unresolved pair {{\small $(i,j)$}} with {{\small $H_{ij}=0$}}, we select {{\small $t_{ij}$}} such that {{\small $i,j\notin S^{t_{ij}}$}}, and define
\[
p_{ij} \in \argmin_{p\in S^{t_{ij}}} D_{ip},
\]
where {{\small $D$}} denotes the distance matrix supplied to this stage (in {{\small SABRE}}, {{\small $D$}} is the Stage~1 estimator {{\small $\widehat D$}}).
The modified comparison sets are then
\[
\tilde L_{ij}
=
\{k\in L_{ij}\cap S^{t_{ij}} : \widehat D^{t_{ij}}_{p_{ij}k}\ge\delta_4\},
\qquad
\tilde R_{ij}
=
\{k\in R_{ij}\cap S^{t_{ij}} : \widehat D^{t_{ij}}_{p_{ij}k}\ge\delta_4\}.
\]

\paragraph*{Remark.}
The construction of the sets {{\small $\tilde L_{ij},\tilde R_{ij}$}} is designed to retain the directional information coming from {{\small $L_{ij},R_{ij}$}}
while using only indices in {{\small $S^{t_{ij}}$}} and the independent estimate {{\small $\widehat D^{t_{ij}}$}}.  The residual dependencies are shown to be benign in Appendix~\ref{appendix:proof-complete-seriation}.


\paragraph*{Refined-seriation routine.}
The final stage of {{\small SABRE}} applies the comparison subroutine to all
unresolved pairs using the modified sets defined above. The full pseudocode is given below.
We recall that {{\small $0_{n\times n}$}} denotes the {{\small $n\times n$}} zero matrix.

\begin{algorithm}[H]
    \caption{{{\small \texttt{Re-evaluate Comparisons}}}} \label{algo:complete-seriation}
    \begin{algorithmic}[1]
        \REQUIRE $(H,D,A,\sigma,\delta_4)$  
        \ENSURE $\widehat H \in \{-1,0,1\}^{n\times n}$
        \STATE Initialize $\widetilde H=0_{n\times n}$
        \STATE \label{partition-3} Randomly partition $[n]$ into $(S^1,S^2,S^3)$ with balanced sizes\label{line-partition-algo}
        \FOR{$t \in \{1,2,3\}$}
        \STATE \label{Dt} Compute $\widehat D^t=$ \texttt{Estimate} \texttt{Distance}$(A,S^t)$ \hfill \COMMENT{Algorithm~\ref{algo:distance}}
        \ENDFOR
        \FOR{$i<j$ with $H_{ij}=0$}
        \STATE \label{left-right-set-estimates} $L_{ij}=\{k : H_{ik}=H_{jk}=1\}$,\quad
               $R_{ij}=\{k : H_{ik}=H_{jk}=-1\}$
        \STATE Choose $t_{ij}$ with $i,j\notin S^{t_{ij}}$
        \STATE \label{k-ij} $p_{ij} \in \argmin_{p\in S^{t_{ij}}} D_{ip}$
        \STATE \label{proxy-set} $\tilde L_{ij}=\{k\in L_{ij}\cap S^{t_{ij}}: 
              \widehat D^{t_{ij}}_{p_{ij}k}\ge\delta_4\}$
        \STATE $\tilde R_{ij}=\{k\in R_{ij}\cap S^{t_{ij}}: 
              \widehat D^{t_{ij}}_{p_{ij}k}\ge\delta_4\}$
        \STATE $\widetilde H_{ij} =$
               \texttt{Evaluate} \texttt{Comparison}$(A_i,A_j,
               \tilde L_{ij},\tilde R_{ij},\sigma)$
        \STATE $\widetilde H_{ji} = -\widetilde H_{ij}$
    \ENDFOR
       \STATE $\widehat H = H + \widetilde H$
    \end{algorithmic}
\end{algorithm}


\section{Proofs of the main results}
\label{appendix:proof-main-theorem}

This appendix contains the proofs of the upper-bound results stated in Section~\ref{section:bilip-mat}.  
We first prove Theorem~\ref{thm-matrix-setting}, then establish Corollary~\ref{cor:bi-lipschitz-matrix} by showing that bi-Lipschitz matrices are a subclass of the average-Lipschitz model. We finally derive the guarantees for Kendall's tau and Frobenius losses stated in Corollaries~\ref{coro-kendall-tau} and~\ref{coro-frob-loss}.

\subsection{Proof of Theorem~\ref{thm-matrix-setting}}
\label{append-proof-thm} 

The argument combines the three key propositions established in the main text --
the distance estimation bound (Proposition~\ref{small-prp-dist-estim-S=n}),  
the first-stage guarantee (Proposition~\ref{prop:end-stage1-new}),  
and the refined-stage guarantee (Proposition~\ref{prop:complete-seriation}).  
Together, these results yield an error bound on the  estimate 
{{\small $\widehat H$}} produced by {{\small SABRE}}.  
To transfer this accuracy to the final permutation estimator, we use the following lemma.

\begin{lem}
\label{lem-last-step-H-perm}
    Let {{\small $\nu > 0$}} and  {{\small $H \in \{-1,0,1\}^{n \times n}$}} such that
    \begin{equation}\label{eq-lem-H-pi}
        \exists s \in \{\pm\}: \ \ \     \quad H_{ij} = s H^*_{ij}\ , \ \ \ \forall i,j \  \text{ such that }  \, |\pi_i - \pi_j| \geq \nu n .
    \end{equation}
    Define {{\small $\si^H=(\si^H_1,\ldots,\si^H_n)$}} with {{\small $\pi_i^H = (H_i \mathbf{1}_n + n+1) /2$}} for all {{\small $i$}}. Then,  {{\small $L_{\max}(\si^H,\si)\leq \nu$}}. 
\end{lem}
Given any  {{\small $H \in \{-1,0,1\}^{n \times n}$}}, and  denoting   {{\small $H_i$}} its {{\small $i$}}-th row,    Lemma~\ref{lem-last-step-H-perm}  shows that the vector {{\small $\si^H$}} is as accurate as {{\small $H$}}.
The proof   is in Appendix~\ref{appendix:technical-lemmes:olf-sumplement-mat}.

\smallskip

Let  {{\small $\cE$}} be the event where  {{\small $\widehat H$}} satisfies~\eqref{eq-lem-H-pi} with {{\small $\nu = (40\sigma/\alpha)\sqrt{\log(n)/n}$}}.
Applying Lemma~\ref{lem-last-step-H-perm} for  {{\small $\hat \pi = \pi^{\wH}$}}, we obtain the desired error bound  {{\small $L_{\max}(\hat \si,\si) \leq (40 \sigma /\alpha) \sqrt{ \log(n)/n }$}}  with probability at least {{\small $\P\ac{\cE}$}}. It remains to prove that {{\small $\P\ac{\cE^c} \leq  1/n^2$}}.

\smallskip
Let  {{\small $\cE'$}} be the event  where  the empirical distance matrix satisfies {{\small $\wD \in \lde(\al,\beta,\om_n,r_0)$}} with {{\small $r_0 = r \wedge r'$}}. 
By Proposition~\ref{small-prp-dist-estim-S=n},  {{\small $\P\ac{\cE'^c}\leq 1/n^4$}}. 
Conditionally on {{\small $\cE'$}}, Proposition~\ref{prop:end-stage1-new} implies that the first estimate {{\small $H$}}  is correct on its support, and for all pairs {{\small $i,j$}} such that {{\small $|\pi_i - \pi_j|\geq \rho$}} for {{\small $\rho = (\de_2+\om_n)/\al$}}. 
Here and throughout, {{\small $\P_{|\mathcal{G}}\ac{\cdot}$}} denotes conditional probability given an event {{\small $\mathcal{G}$}}.  
Proposition~\ref{prop:complete-seriation} then yields {{\small $\P_{|\cE'}\ac{\cE^c} \le 9/n^3$}}.

Using the law of total probability, we obtain
\[\P\ac{\cE^c} \leq \P_{|\cE'}\ac{\cE^c} + \P\ac{\cE'^c}  \leq 1/n^2 , \qquad  n\geq 10.\]
This completes the proof of Theorem~\ref{thm-matrix-setting}. \hfill $\square$

\medskip
The same reasoning applies to the generalized setting of Theorem~\ref{thm-matrix-setting-general}: 
one only needs to replace the assumption of fixed parameters by the calibration conditions~(\ref{delta-constraint-1rst-seriation}--\ref{constraint-2nd-seriation}) 
and takes {{\small $\omega=\omega_n$}} as in  Proposition~\ref{small-prp-dist-estim-S=n}.


\subsection{Proof of Corollary~\ref{cor:bi-lipschitz-matrix}}
\label{append-proof-cor-bi-lip}

To prove Corollary~\ref{cor:bi-lipschitz-matrix}, it suffices to show that
{{\small $\mathcal{BL}(\alpha,\beta)\subseteq\mathcal{AL}(\alpha/4,\beta,r,r')$}} for some constants {{\small $(r, r')$}}
and then apply Theorem~\ref{thm-matrix-setting}.
This inclusion follows from the next lemma.

\begin{lem}[BL $\Rightarrow$ AL]
\label{lem:BL-AL}
Fix {{\small $\alpha, \beta > 0$}}. There exists a constant {{\small $c_{\alpha\beta}>0$}}, depending only on {{\small $(\alpha,\beta)$}}, such that for all {{\small $n\ge 11 \vee c_{\alpha\beta}$}},
any {{\small $F\in\mathcal{BL}(\alpha,\beta)$}} belongs to {{\small $\mathcal{AL}(\alpha/4,\beta,r,r')$}} for every  {{\small $r\in(0,1/4]$}} and {{\small $r'$}} satisfying
\[
r' < \min\!\left\{
\frac{\alpha}{\sqrt{2}}\,r,\
\frac{\alpha^{3/2}}{8\sqrt{\alpha+\beta}}
\right\}.
\]
\end{lem}

This completes the proof of the inclusion 
{{\small $\mathcal{BL}(\alpha,\beta)\subseteq\mathcal{AL}(\alpha/4,\beta,r,r')$}},
and Corollary~\ref{cor:bi-lipschitz-matrix} then follows directly from
Theorem~\ref{thm-matrix-setting} with {{\small $\alpha$}} replaced by {{\small $\alpha/4$}}.

\begin{proof}[Proof of Lemma~\ref{lem:BL-AL}.]
\textit{Average {{\small $\ell_2$}} upper bound.}
By definition of the bi-Lipschitz class, for all {{\small $i,j\in[n]$}},
\[
\|F_i-F_j\|\le \beta\,\frac{|i-j|}{\sqrt{n}} ,
\]
which exactly matches the average {{\small $\ell_2$}} upper bound in Definition~\ref{defi:PBL}.

\smallskip 
\textit{Average {{\small $\ell_1$}} lower bound.}
Let {{\small $r\in(0,1/4]$}}  and {{\small $i<j$}} with {{\small $j-i\le r n$}}.  
For any subsets {{\small $L\subset\{k<i\}$}} and {{\small $R\subset\{k>j\}$}},
the bi-Lipschitz inequalities yield
\[
\sum_{k\in L}(F_{ik}-F_{jk})
\ \vee\
\sum_{k\in R}(F_{jk}-F_{ik})
\ \ge\
(|L|\vee|R|)\,\alpha\,\frac{|i-j|}{n}.
\]
Here {{\small $|\cdot|$}} denotes set cardinality.
Set {{\small $c_0=1/32$}} and take
{{\small $L_0=\{k<i-c_0n\}$}} and {{\small $R_0=\{k>j+c_0n\}$}}.
Since {{\small $j-i \le r n$}}, we necessarily have 
either {{\small $i \ge (1-r)n/2$}} or {{\small $j \le (1+r)n/2$}}.
Consequently,
\[|L_0|\vee|R_0|
\ \ge\
\Big(\frac{1-r}{2}-c_0\Big)n-1
\ \ge\
\frac{n}{4}
\qquad (n\ge 11).\]

Hence
\[
\sum_{k<i-c_0n}(F_{ik}-F_{jk})
\ \vee\
\sum_{k>j+c_0n}(F_{jk}-F_{ik})
\ \ge\
\frac{\alpha}{4}\,|i-j|,
\]
which establishes the average {{\small $\ell_1$}} lower bound with {{\small $\alpha'=\alpha/4$}}.

\smallskip
\textit{Non-collapse at large distance.} Let {{\small $r\in(0,1/4]$}} and {{\small $i<j$}} such that {{\small $|i-j| > rn$}}.
Assume {{\small $n\ge c_{\alpha\beta}$}} so that Lemma~\ref{lem:BL-separation} applies.
If {{\small $|i-j|\le n/2$}}, then
\[
\|F_i-F_j\|\ \ge\ \frac{\alpha}{\sqrt{2}}\frac{|i-j|}{\sqrt{n}}
\ \ge\
\frac{\alpha}{\sqrt{2}}\,r\,\sqrt{n}.
\]
If instead {{\small $|i-j|>n/2$}}, the same lemma gives
\[
\|F_i-F_j\|\ \ge\
\frac{\alpha^{3/2}}{4\sqrt{\alpha+\beta}}\frac{|i-j|}{\sqrt{n}}
\ \ge\
\frac{\alpha^{3/2}}{8\sqrt{\alpha+\beta}}\sqrt{n}.
\]
Therefore, the non-collapse condition in Definition~\ref{defi:PBL} holds with any {{\small $r'$}} satisfying
\[
r' \  < \ \min\!\left\{
\frac{\alpha}{\sqrt{2}}\,r,\
\frac{\alpha^{3/2}}{8\sqrt{\alpha+\beta}}
\right\}.
\]
Lemma~\ref{lem:BL-AL} follows.
\end{proof}


The proof of Lemma~\ref{lem:BL-AL} relies on the following 
 {{\small $\ell_2$}} separation property for bi-Lipschitz matrices.

\begin{lem}[Average {{\small $\ell_2$}} lower bound under bi-Lipschitz matrices]
\label{lem:BL-separation}
Let  {{\small $F\in\mathcal{BL}(\alpha,\beta)$}}.
Then for all {{\small $i,j\in[n]$}}, 
\[
\|F_i-F_j\|\ \ge  \ C_{\alpha \beta  n}  \frac{|i-j|}{\sqrt{n}} , \quad \textup{ where } \ \ 
C_{\alpha \beta n} = 
\begin{cases}
 \frac{\alpha}{\sqrt{2}}, & \text{if } |i-j| \le n/2,\\[3pt]
\sqrt{\left(\frac{\alpha}{2(\alpha+\beta)} - \frac{2}{n}\right)_+} \  \frac{\alpha}{2}, & \text{otherwise}.
\end{cases}
\]
In particular,  there exists some constant $c_{\alpha \beta}$  depending only on $(\alpha, \beta)$, such that for all $n \geq c_{\alpha \beta}$, and  $i,j\in[n]$,
\[
\|F_i-F_j\|\ \ge  \ C_{\alpha \beta}  \frac{|i-j|}{\sqrt{n}} , \quad \textup{ where } \ \ 
C_{\alpha \beta} = 
\begin{cases}
 \frac{\alpha}{\sqrt{2}}, & \text{if } |i-j| \le n/2,\\[3pt]
\frac{\alpha^{3/2}}{4\sqrt{\alpha+\beta}}, & \text{otherwise}.
\end{cases}
\]
\end{lem}

\begin{proof}[Proof of Lemma~\ref{lem:BL-separation}.]
Fix {{\small $i<j$}} and write {{\small $L := |i-j|$}}.

$\circ$ If {{\small $L \le n/2$}}, then the number of indices {{\small $k\in[n]$}} such that {{\small $k \le i$}} or {{\small $k \ge j$}} is lower bounded by {{\small $n/2$}}. In addition, for {{\small $k\le i$}} or {{\small $k\ge j$}},
  \[
  |F_{ik}-F_{jk}|\ \ge\   \frac{\alpha L}{n} ;
  \]
hence {{\small $\|F_i-F_j\|\ \ge\ \sqrt{\frac{n}{2}} \frac{\alpha L}{n}  \ = \ \frac{\alpha L}{\sqrt{2 n}}$}}.

$\circ$ Now assume that {{\small $L > n/2$}}. By telescoping over adjacent rows,
\begin{equation}
\label{telescopic-sum}
F_{ik}-F_{jk}=\sum_{t=i}^{j-1}\big(F_{t,k}-F_{t+1,k}\big).
\end{equation}
From the definition of {{\small $\mathcal{BL}(\alpha,\beta)$}} applied to each adjacent pair {{\small $(t,t+1)$}}:
\[
F_{t,k}-F_{t+1,k}\in
\begin{cases}
[\alpha/n,\ \beta/n], & \text{if } k\le t,\\[3pt]
[-\beta/n,\ -\alpha/n], & \text{if } k\ge t+1.
\end{cases}
\]
Consider the interior case {{\small $i<k<j$}}.
In the sum~\eqref{telescopic-sum}:\\
-- for {{\small $t=i,\dots,k-1$}} (exactly {{\small $k-i$}} terms) we have {{\small $F_{t,k}-F_{t+1,k}\in[-\beta/n,-\alpha/n]$}};\\
-- for {{\small $t=k,\dots,j-1$}} (exactly {{\small $j-k$}} terms) we have {{\small $F_{t,k}-F_{t+1,k}\in[\alpha/n,\beta/n]$}}.

Adding these terms gives the bounds
\begin{equation}
\label{encadrement-valeur-interieur}
\frac{\alpha(j-k)-\beta(k-i)}{n}
\ \le\
F_{ik}-F_{jk}
\ \le\
\frac{\beta(j-k)-\alpha(k-i)}{n}.
\end{equation}

We want to find indices  {{\small $k$}} for which the above difference takes large values.
Choose {{\small $\theta\in(0,1/2]$}} and define
\[
\mathcal{K}_i :=\{\,k:\ i<k\le i+\lfloor \theta L\rfloor\,\},\qquad
\mathcal{K}_j:=\{\,k:\ j-\lfloor \theta L\rfloor\le k<j\,\}.
\]
These two disjoint sets each contain {{\small $\lfloor \theta L\rfloor$}} consecutive columns near {{\small $i$}} and {{\small $j$}}.

For {{\small $k\in \mathcal{K}_i$}}, we have {{\small $k-i\le \theta L$}} and {{\small $j-k\ge (1-\theta)L$}}.
Applying the lower bound in~\eqref{encadrement-valeur-interieur} gives
\[
F_{ik}-F_{jk}\ \ge\ \frac{\alpha(1-\theta)L-\beta\theta L}{n}
=\frac{(\alpha(1-\theta)-\beta\theta)L}{n}.
\]
For {{\small $k\in \mathcal{K}_j$}}, we have {{\small $k-i\ge (1-\theta)L$}} and {{\small $j-k\le \theta L$}}.
Applying the upper bound in~\eqref{encadrement-valeur-interieur} gives
\[
F_{ik}-F_{jk}\ \le\ -\,\frac{(\alpha(1-\theta)-\beta\theta)L}{n}.
\]
Thus on {{\small $\mathcal{K}_i\cup \mathcal{K}_j$}} we have
\begin{equation}
\label{final-difference-up-bound}
|F_{ik}-F_{jk}|\ \ge\ \frac{(\alpha(1-\theta)-\beta\theta)_+\,L}{n}.
\end{equation}

We now choose {{\small $\theta$}} so that {{\small $\alpha(1-\theta)-\beta\theta = \alpha/2$}}.
Solving this gives {{\small $\displaystyle \theta:=\frac{\alpha}{2(\alpha+\beta)}$}}.

Finally, we sum over {{\small $\mathcal{K}_i$}} and {{\small $\mathcal{K}_j$}}.
Using~\eqref{final-difference-up-bound} and {{\small $|\mathcal{K}_i\cup \mathcal{K}_j|\ge 2\lfloor \theta L\rfloor\geq (2\theta L -2)_+$}},
\[
\|F_i-F_j\|^2
=\sum_{k=1}^n (F_{ik}-F_{jk})^2
\ \ge\ |\mathcal{K}_i\cup \mathcal{K}_j|\cdot\Big(\frac{\alpha L}{2n}\Big)^2
\ \ge\ \left(\theta - \frac{1}{L}\right)_+\frac{\alpha^2}{2}\,\frac{L^3}{n^2}.
\]
Since {{\small $L\ge n/2$}}, 
\[
\|F_i-F_j\|
\ \ge\ \sqrt{\left(\theta - \frac{2}{n}\right)_+} \ \frac{\alpha}{2}\,\frac{L}{\sqrt{n}}
.
\]
This completes the proof of Lemma~\ref{lem:BL-separation}.
\end{proof}



\subsection{Proofs of Corollary~\ref{coro-kendall-tau} and~\ref{coro-frob-loss}} 
\label{appendix-coro-losses-proof}
$\circ$ Using the equivalence~\eqref{equivalence-losses} between the loss functions {{\small $L_1$}} and {{\small $L_K$}}, we have
\[L_{K}(\hat \pi, \pi) \ \ \leq \ \  L_1(\hat \pi, \pi) \ \  \leq  \ \    n  L_{\max}(\hat \pi, \pi) \ \ \leq \ \  40\frac{\sigma}{\alpha}\sqrt{n \log n}\enspace,\]
where the last inequality holds with probability at least {{\small $1-1/n^2$}}, by  Theorem~\ref{thm-matrix-setting}. Corollary~\ref{coro-kendall-tau} follows. 

 $\circ$  By definition of the Frobenius norm, which is {{\small $\|M\|= \sqrt{\sum_{1\leq i,j \leq n} M_{ij}^2}$}} for any matrix {{\small $M\in \mathbb{R}^{n \times n}$}}, we have 
\begin{align}
\label{L-F}
    L_{F}(\hat \pi, \pi) \ \  &\leq \ \ \big{(} n  \max_{1\leq i, j\leq n} |F_{\pi_i \pi_j} -F_{\hat \pi_i \hat  \pi_j} | \big{)} \ \ \wedge  \ \ \big{(}n \max_{1\leq i, j\leq n} |F_{\pi_i^{\mathrm{rev}} \pi_j^{\mathrm{rev}}} -F_{\hat \pi_i \hat  \pi_j} | \big{)} \enspace.
\end{align}
Using the triangle inequality and {{\small $F\in \mathcal{BL}(\alpha,  \beta)$}} we also have
\begin{align*}
    |F_{\pi_i \pi_j} -F_{\hat \pi_i \hat  \pi_j} | \ \ \leq  |F_{\pi_i \pi_j} -F_{\hat \pi_i  \pi_j} | \ +  \ |F_{\hat \pi_i \pi_j} -F_{\hat \pi_i  \hat \pi_j} |
    &\leq \ \ 2 \frac{\beta}{n} \max_{i\in[n]} | \pi_i - \hat \pi_i |\enspace.
\end{align*}
Repeating the same argument for {{\small $\pi^{\mathrm{rev}}$}} (instead of {{\small $\pi$}}) we obtain a similar upper bound.  Plugging the two  into~\eqref{L-F} we obtain
\begin{align*}
    L_{F}(\hat \pi, \pi) \ \  &\leq \ \ 2 \beta \big{(} \max_{i\in[n]} | \pi_i - \hat \pi_i |  \ \wedge \ \max_{i\in[n]} | \pi_i^{\mathrm{rev}} - \hat \pi_i | \big{)} \\
    &= \ \ 2 \beta n L_{\max}( \hat \pi, \pi) \\
    &\leq \ \  320 \beta \frac{ \sigma}{\alpha} \sqrt{n \log n} \enspace.
\end{align*}
where the last inequality holds with probability at least {{\small $1-1/n^2$}}, by  Corollary~\ref{cor:bi-lipschitz-matrix}. 
This completes the proof of  Corollary~\ref{coro-frob-loss}. \hfill $\square$


\section{Distance estimation}
\label{appendix-distance-estimation}

This appendix provides the proof of Proposition~\ref{small-prp-dist-estim-S=n}.  
We first establish that the average-Lipschitz condition implies local distance equivalence for the population matrix {{\small $D^*$}}.
We then prove that the  estimator {{\small $\widehat D$}} concentrates uniformly around {{\small $D^*$}}, 
and extend this result to the case of subsampled matrices {{\small $A^S$}} used in later stages of {{\small SABRE}}.


\subsection{From average-Lipschitz to local distance equivalence}

To connect the local distance equivalence property to our structural assumption, 
we show that the average-Lipschitz condition implies local distance equivalence.

\begin{lem}[From average-Lipschitz to local distance equivalence]
\label{lem:AL-implies-LDE}
If {{\small $F\in\mathcal{AL}(\alpha,\beta,r,r')$}}, then its row-distance matrix  {{\small $D^*$}}, defined in~\eqref{def-dist-[n]}, belongs to {{\small $\lde(\alpha,\beta,0,r_0)$}} 
with {{\small $r_0= r \wedge r'$}}.
\end{lem}

Hence, under the average-Lipschitz condition, 
the distance matrix {{\small $D^*$}} satisfies the local distance equivalence of Definition~\ref{cond-LB}.

\begin{proof}[Proof of Lemma~\ref{lem:AL-implies-LDE}]

We verify the two inequalities of Definition~\ref{cond-LB}
for all pairs {{\small $(i,j)$}} such that  
{{\small $|\pi_i-\pi_j|\wedge D^*_{ij}\le n r_0$}}.

\smallskip
\textit{Upper bound.}
If {{\small $|\pi_i-\pi_j|\le n r$}},  
the average {{\small $\ell_2$}} upper bound in the definition of  
{{\small $\mathcal{AL}(\alpha,\beta,r,r')$}} gives
\[
D^*_{ij}
=\sqrt{n}\,\|F_{\pi_i}-F_{\pi_j}\|
\ \le\ \beta\,|\pi_i-\pi_j|,
\]
which is exactly the desired upper inequality with {{\small $\omega=0$}}.

\smallskip
\textit{Lower bound.}
Assume that {{\small $\pi_i < \pi_j$}} and 
{{\small $|\pi_i-\pi_j|\le n r$}}.  
The average {{\small $\ell_1$}} lower bound in
{{\small $\mathcal{AL}$}} ensures that,
there exists a side
{{\small $S\in\big\{\{k<\pi_i-c_0n\},\,\{k>\pi_j+c_0n\}\big\}$}} 
such that
\[
\sum_{k\in S} | F_{\pi_i k}-F_{\pi_j k} | \ \ge\ \alpha\,|\pi_i-\pi_j| .
\]
Hence {{\small $\|F_{\pi_i}-F_{\pi_j}\|_1\ge\alpha\,|\pi_i-\pi_j|$}}, and by Cauchy--Schwarz,
\[
D^*_{ij}
=\sqrt{n}\,\|F_{\pi_i}-F_{\pi_j}\|
\ \ge\ \|F_{\pi_i}-F_{\pi_j}\|_1
\ \ge\ \alpha\,|\pi_i-\pi_j| .
\]
Thus the lower bound also holds with {{\small $\omega=0$}}.

\smallskip
\textit{Role of {{\small $r_0= r \wedge r'$}}.}
The non-collapse condition in  
{{\small $\mathcal{AL}$}} guarantees that  
{{\small $|\pi_i-\pi_j|>n r\Rightarrow D^*_{ij}>n r'$}}.
Therefore, whenever  
{{\small $|\pi_i-\pi_j|\wedge D^*_{ij}\le n r_0$}},  
we necessarily have {{\small $|\pi_i-\pi_j|\le n r$}},  
and the previous inequalities apply.

\smallskip
Combining the two bounds yields
\[
\alpha\,|\pi_i-\pi_j|
\ \le\ D^*_{ij}\ \le\ 
\beta\,|\pi_i-\pi_j|
\qquad\text{whenever}\qquad
|\pi_i-\pi_j|\wedge D^*_{ij}\le n r_0.
\]
Hence {{\small $D^*\in\lde(\alpha,\beta,0,r_0)$}}.
\end{proof}


\subsection{Uniform concentration of the empirical distance estimator}
\label{appendix:distance-estimation-fullsample}

This section provides the detailed statement and proof of the uniform concentration bound 
announced in Proposition~\ref{small-prp-dist-estim-S=n}, which characterizes the deviation of the empirical distance estimator {{\small $\widehat D$}} from its population counterpart {{\small $D^*$}}.

The assumption {{\small $F\in\mathcal{AL}(\alpha,\beta,r,r')$}} is  stronger than needed here;
we only use its local {{\small $\ell_2$}} upper bound.

\paragraph*{Statement of the distance error bound.}
There exists a numerical constant {{\small $C$}} such that, for any {{\small $\beta,r>0$}} and {{\small $n\ge 8\vee(1/r)$}},
if {{\small $F\in[0,1]^{n\times n}$}} satisfies the local {{\small $\ell_2$}} upper bound of the average-Lipschitz condition,
then the estimator {{\small $\widehat D$}} satisfies the following uniform concentration bound:
\begin{equation}\label{distance-estimate-simple-error}
 \max_{i,j \in [n]} |\widehat D_{ij} - D^*_{ij}| 
 \ \leq\ 
 C\Big( \sqrt{ \beta n} \,+\, \sqrt{ (\sigma+1)\sigma}\, n^{3/4} (\log n)^{1/4} \Big),
\end{equation}
with probability at least {{\small $1-1/n^4$}}.  
This directly  yields the bound stated in Proposition~\ref{small-prp-dist-estim-S=n} with {{\small $\omega_n := C_{\beta\sigma}\, n^{3/4} (\log n)^{1/4}$}}
for the  explicit constant  {{\small $C_{\beta\sigma} = C (\sqrt{\beta} +  \sqrt{ (\sigma+1)\sigma})$}}.

The bound~\eqref{distance-estimate-simple-error} reflects two sources of error:
a deterministic bias, {{\small $\sqrt{\beta  n }$}}, due to the nearest-neighbor approximation used in the construction of {{\small $\widehat D$}} in Section~\ref{subsection:l2-dist},
and random fluctuations, {{\small $\sqrt{ (\sigma+1)\sigma}\, n^{3/4} (\log n)^{1/4}$}}, induced by the sub-Gaussian noise matrix {{\small $E$}}.  
As a consequence, when {{\small $(\beta,\sigma)$}} are fixed constants, the normalized error
{{\small $n^{-1}\max_{i,j}|\widehat D_{ij}-D^*_{ij}|$}} vanishes as {{\small $n\to\infty$}}.

\paragraph*{Proof of the distance error bound~\eqref{distance-estimate-simple-error}.} 
The argument relies on the following auxiliary lemma, which controls the squared-distance errors
of the  estimator {{\small $\widehat D$}}.

\begin{lem}
\label{lem:dist-bound}
Let {{\small $n \ge 4$}}.  
Denote by {{\small $m_i \in [n]$}} a nearest neighbor of {{\small $i$}} with respect to {{\small $D^*$}},
that is, {{\small $m_i \in \argmin_{t\ne i} D^*_{it}$}}.
Then the estimator {{\small $\widehat D$}} computed from {{\small $A$}} by Algorithm~\ref{algo:distance} satisfies
\[
 \frac{1}{n}\,\max_{i,j\in[n]} \big|(D^*_{ij})^2 - (\widehat D_{ij})^2\big|
 \ <\ C'\,\omega'_n,
\]
for some numerical constant {{\small $C'$}}, with probability at least {{\small $1-1/n^4$}},
where
\begin{equation}\label{dist-error-bound-small-lemma}
 \omega'_n
 \ :=\
 \tfrac{\|F\|_\infty}{\sqrt{n}}\,
 \max_{i\in[n]} D^*_{i m_i}
 \ +\
 \Big(\sigma + \tfrac{\|F\|_\infty}{\sqrt{n}}\Big)
 \sigma \sqrt{n \log n}.
\end{equation}
Here {{\small $\|F\|_\infty := \max_{i\in[n]} \|F_i\|$}} denotes the largest Euclidean norm of the rows of {{\small $F$}}.
\end{lem}

Lemma~\ref{lem:dist-bound} extends the result of~\cite{issartel2021estimation} to real-valued matrices {{\small $A$}} corrupted by sub-Gaussian noise.
A sketch of its proof is given at the end of this section. 

\smallskip
We now derive~\eqref{distance-estimate-simple-error} from Lemma~\ref{lem:dist-bound}. Recall that {{\small $F\in[0,1]^{n\times n}$}}.
Then each row satisfies {{\small $\|F_i\|\le\sqrt{n}$}}, so {{\small $\|F\|_\infty\le\sqrt{n}$}}.\\
Recall that the rows of {{\small $F$}} satisfy the  {{\small $\ell_2$}} upper bound
{{\small $\|F_{\pi_i}-F_{\pi_j}\|\le (\beta/\sqrt n)|\pi_i-\pi_j|$}}
for all {{\small $|\pi_i-\pi_j|\le rn$}}.
By definition of the nearest neighbor {{\small $m_i$}}, we have {{\small $D^*_{i m_i} \le D^*_{i k}$}} for all {{\small $k\neq i$}}.
In particular, taking an adjacent index {{\small $i_c$}} with {{\small $|\pi_i-\pi_{i_c}|=1\le rn$}},
the  {{\small $\ell_2$}} upper bound yields {{\small $D^*_{i i_c}\le \beta |\pi_i-\pi_{i_c}| = \beta$}},
and hence {{\small $D^*_{i m_i}\le \beta$}}.\\
Plugging these bounds into~\eqref{dist-error-bound-small-lemma} gives
\[ \omega'_n \ \le\  \beta  \ +\ (\sigma+1)\sigma \sqrt{n\log n} .\]
Using the numeric inequality {{\small $|a-b|\le \sqrt{|a^2-b^2|}$}} for any {{\small $a,b \ge 0$}},
we obtain from Lemma~\ref{lem:dist-bound}
\[
 \max_{i,j}|D^*_{ij}-\widehat D_{ij}|
 \ <\
 \sqrt{C' n\,\omega'_n}.
\]
Applying {{\small $\sqrt{a+b}\le\sqrt{a}+\sqrt{b}$}}
yields the bound~\eqref{distance-estimate-simple-error}. 
\hfill$\square$

\medskip 
\paragraph*{Proof sketch of Lemma~\ref{lem:dist-bound}.}
We conclude this section with the main steps of Lemma~\ref{lem:dist-bound},
which follow the argument of~\cite{issartel2021estimation} and are included here for completeness.

In the nearest-neighbor approximation of Section~\ref{subsection:l2-dist},
we replaced the quadratic term {{\small $\langle F_{\pi_i}, F_{\pi_i}\rangle$}}
by the cross term {{\small $\langle F_{\pi_i}, F_{\pi_{m_i}}\rangle$}},
which introduces an approximation error bounded by
\begin{align}\label{neighbError}
 \big|\langle F_{\pi_i}, F_{\pi_i}-F_{\pi_{m_i}}\rangle\big|
 \ \le\
 \|F_{\pi_i}\|\,\|F_{\pi_i}-F_{\pi_{m_i}}\|
 \ \le\
 \|F\|_\infty \max_i D^*_{i m_i}/\sqrt{n},
\end{align}
using Cauchy–Schwarz.
The empirical cross term {{\small $\langle A_i, A_j\rangle$}} (for {{\small $i\neq j$}}) concentrates around
{{\small $\langle F_{\pi_i},F_{\pi_j}\rangle$}}, and with probability at least {{\small $1-1/n^4$}},
\begin{equation}\label{eq-scalar-prod-concentration}
 \forall\,i<j:\quad
 |\langle A_i,A_j\rangle-\langle F_{\pi_i},F_{\pi_j}\rangle|
 \ \le\
 C\Big(\sigma+\tfrac{\|F\|_\infty}{\sqrt{n}}\Big)
 \sigma\sqrt{n\log n}.
\end{equation}
The proof of~\eqref{eq-scalar-prod-concentration} is provided in Appendix~\ref{appendix:proof-corssed-term-conentration-dist-elementary}.
Combining~\eqref{neighbError} and~\eqref{eq-scalar-prod-concentration} provides the main ingredients underlying the bound~\eqref{dist-error-bound-small-lemma}.
This concludes the proof sketch of Lemma~\ref{lem:dist-bound}.
\hfill$\square$



\subsection{Extension of the distance estimation to subsets}
\label{appendix-extension-dist-estim}

We now extend the uniform concentration bound~\eqref{distance-estimate-simple-error} 
to the case where Algorithm~\ref{algo:distance} is applied to a submatrix 
{{\small $A^S$}} built from a subset {{\small $S \subset [n]$}}. 
We denote the resulting estimator by {{\small $\widehat D(S)$}}.

\paragraph*{Population counterpart.}
This estimator is not meant to approximate {{\small $D^*$}} directly, 
but the corresponding population quantity {{\small $D^*(S)$}}, defined by
\begin{equation}\label{def-dist-neighb}
D^*_{ij}(S)
\ :=\
n\sqrt{\frac{1}{|S|}\sum_{k\in S}(F_{\pi_i\pi_k}-F_{\pi_j\pi_k})^2},
\qquad
i,j\in S,
\end{equation}
(where {{\small $|S|$}} is the cardinality of {{\small $S$}}),
and {{\small $D^*_{ij}(S)=0$}} whenever {{\small $i\notin S$}} or {{\small $j\notin S$}}.
This quantity naturally extends Definition~\eqref{def-dist-[n]} of {{\small $D^*$}} to subsets of indices.
For {{\small $S=[n]$}}, we recover {{\small $D^*(S)=D^*$}}.

\paragraph*{Local distance equivalence on subsets.}
We denote by {{\small $\lde(\alpha,\beta,\omega,r,S)$}} 
the class of matrices {{\small $D\in\R^{n\times n}$}} 
whose entries {{\small $D_{ij}$}} for {{\small $i,j\in S$}} 
satisfy the local distance equivalence condition of Definition~\ref{cond-LB}.
This is the natural extension of {{\small $\lde(\alpha,\beta,\omega,r)$}} to subsets {{\small $S\subset[n]$}}.

\paragraph*{Extension of the distance bound.}
We denote by {{\small $\eta$}}  any upper bound on the  nearest-neighbor distances within {{\small $S$}}, 
\begin{equation}\label{nearest-neighbor-dist-upper-bound}
\max_{i\in S}\ \min_{k\in S,\ k\ne i}\ |\pi_i-\pi_k|
\ \le\ 
\eta(S):=\eta\; .
\end{equation}
There exists a numerical constant {{\small $C$}} such that, for any {{\small $S\subset[n]$}}, any {{\small $\beta,r, \eta >0$}} and {{\small $n\ge 8\vee (\eta/r)$}},
if {{\small $F\in[0,1]^{n\times n}$}} satisfies the local {{\small $\ell_2$}} upper bound of the average-Lipschitz condition on {{\small $S$}},
 then the estimator {{\small $\widehat D(S)$}} satisfies the following uniform concentration bound
\begin{equation}\label{distance-estimate-condition}
\max_{i,j\in S}|\widehat D_{ij}(S)-D^*_{ij}(S)|
\ \le\
C\!\left(\sqrt{\beta\eta n}
+\sqrt{(\sigma+1)\sigma}\,n^{3/4}(\log n)^{1/4}\right)
\ := \ \omega_n(\eta),
\end{equation}
with probability at least {{\small $1-1/n^4$}}.
This error bound has the same structure as in the full-sample case 
({{\small $S=[n]$}} in Appendix~\ref{appendix:distance-estimation-fullsample}), with the bias term now scaled by {{\small $\sqrt{\eta}$}}.

\paragraph*{Special cases.}
When {{\small $S=[n]$}} (as in the first stage of {{\small SABRE}}), 
we have {{\small $\eta=1$}} and~\eqref{distance-estimate-condition} reduces to~\eqref{distance-estimate-simple-error}.
When {{\small $S\subset[n]$}} is a random subset of size at least {{\small $|S|\ge n/3$}} 
(as in the third stage of {{\small SABRE}}), 
we can show that {{\small $\eta(S)\lesssim \sqrt{n\log n}$}} with high probability, 
leading again to the same  rate
\begin{equation}\label{om-n-def}
\omega_n
=
C_{\beta\sigma}\,n^{3/4}(\log n)^{1/4}.
\end{equation}

\paragraph*{Summary.}
The following lemma summarizes this generalization.

\begin{lem}
\label{lem-hypothesis:dist-equiv-estimator} 
For any {{\small $S\subset[n]$}}, any {{\small $\beta,r, \eta >0$}} and {{\small $n\ge 8\vee (\eta/r)$}},
if {{\small $F\in[0,1]^{n\times n}$}} satisfies the local {{\small $\ell_2$}} upper bound of the average-Lipschitz condition on {{\small $S$}},
 then the estimator {{\small $\widehat D(S)$}} satisfies the concentration bound~\eqref{distance-estimate-condition} with probability at least {{\small $1-1/n^4$}}. As a consequence, 
for any {{\small $\alpha, \beta, r,\omega>0$}} such that {{\small $D^*(S)\in\lde(\alpha,\beta,\omega,r,S)$}}, 
we have {{\small $\widehat D(S)\in\lde(\alpha,\beta,\omega+\omega_n(\eta),r,S)$}} with the same probability.
\end{lem}
Its proof follows the same argument as that of the full-sample bound~\eqref{distance-estimate-simple-error},
with an additional factor {{\small $\sqrt{\eta}$}} arising from the larger nearest-neighbor distances within the subset {{\small $S$}}.
This completes the extension of the distance estimation bound to arbitrary subsets {{\small $S\subset[n]$}}.




\section{First seriation}
\label{appendix-aggregation-bisections}

This appendix provides the proof of Proposition~\ref{prop:end-stage1-new}, which establishes the deterministic guarantee for the first seriation step.

\paragraph*{Notation.}
In this appendix, we write {{\small $\{\textup{P}(k)\}$}} for the set {{\small $\{k\in[n] : \textup{P}(k)\ \text{holds}\}$}}.  
Given a set {{\small $G$}}, we denote by {{\small $\pi_G$}} the collection {{\small $\{\pi_g\}_{g\in G}$}}, and {{\small $|G|$}} its cardinality.  


\subsection{Proof of Proposition~\ref{prop:end-stage1-new}}\label{append-proof-1rst-seriation}
In this appendix, we set {{\small $\rho := (\de_2+\om)/\al$}}.  Lemma~\ref{lem-bisection} shows that {{\small $(G_i,G_i')$}} is a rough bisection  of {{\small $[n]\setminus\{i\}$}} with an error less than {{\small $\rho$}} with respect to the ordering $\pi$.

\begin{lem}[Bisections]\label{lem-bisection}
Let  {{\small $n, \alpha, \beta, r, \om$}} and {{\small  $\de_1,\de_2, \de_3$}} as in~\eqref{delta-constraint-1rst-seriation}. If  {{\small $D \in \lde(\alpha,\beta,\om,r)$,}} then for   all {{\small $i\in [n]$}}  we have {{\small $G_i \neq \emptyset$}}, and
 \begin{enumerate}
  \item  $\pi_{G_i}$ is on one side  of $\pi_i$, and includes all {{\small $k\in[n]$}} that lie at least {{\small $\rho$}} away from {{\small $\pi_i$}} on this side :
  \[  \{k  \leq \pi_i -  \rho\} \subset  \pi_{G_i} \subset \{k < \pi_i\} \qquad \textrm{ or } \qquad  \{k \geq \pi_i + \rho\} \subset \pi_{G_i} \subset \{k > \pi_i\} \enspace.\] 
 \item if {{\small $G_i' \neq \emptyset$}}, then {{\small $\pi_{G_i}$}} and {{\small $\pi_{G'_i}$}} are on opposite sides of $\pi_i$, and  each satisfies property~{{\small \textup{1}.}} 
 \item if {{\small $G_i'  = \emptyset$}}, then either {{\small $n\in \pi_{G_i}$}} with {{\small  $|\pi_i - 1|  < n/8$}}, or {{\small $1\in \pi_{G_i}$}} with {{\small  $|\pi_i - n|  <  n/8$}}. 
\end{enumerate}
\end{lem} 
Despite the regularity ensured by property~1, the set   {{\small $G_i$}} has possibly an exponentially high number of realizations. This high entropy is one key difficulty that motivated our data-splitting scheme in the refined seriation step (Section~\ref{subsection:stage2}).

To better illustrate Algorithm~\ref{algo:SAB}, Appendix~\ref{proof-lem-bisection} provides a proof of Lemma~\ref{lem-bisection}  in the ideal and simple scenario where {{\small $D \in \lde(1,1,0,1)$}}. The general scenario where {{\small $D \in \lde(\alpha,\beta,\om,r)$}} is in Appendix~\ref{append-proof-lem-bisection-general-case}. 

Given the above  bisections {{\small $(G_i, G'_i)_{i\in[n]}$}}, Lemma~\ref{lem-orientation} shows that the {{\small $n$}} pairs {{\small $(L_i, R_i)_{i\in[n]}$}} share the same 'left-right' orientation. (The  proof  is at the end of this  appendix.)

\begin{lem}[Orientation]
\label{lem-orientation}
Assume  (without loss of generality) that {{\small $\argmin_{i \in L_c \cup  R_c}    \si_{i}  \in  L_{c}$}}.
If the assumptions of Lemma~\ref{lem-bisection} are satisfied, then  for all {{\small $i\in[n]$}} we have {{\small $L_i \cup R_i \neq \emptyset$}}, and 
\begin{enumerate}
    \item if {{\small $L_i\neq \emptyset$}}, then {{\small $ \{k  \leq \pi_i -  \rho\} \subset  \pi_{L_i} \subset \{k <\pi_i\}$}}; 
    \item if {{\small $R_i\neq \emptyset$}}, then {{\small $ \{k \geq  \pi_i +  \rho\} \subset \pi_{R_i} \subset \{k > \pi_i\}$}}; 
    \item if {{\small $L_i  = \emptyset$}}, then  {{\small $|\pi_i - 1|  < n/8$}}; similarly, if {{\small $R_i  = \emptyset$}}, then {{\small $|\pi_i - n| < n/8$}}.  
\end{enumerate}
\end{lem}
For brevity in Lemma~\ref{lem-orientation}, we assumed  the canonical  `left-right' orientation where\\ {{\small $\argmin_{i \in L_c \cup  R_c}    \si_{i}  \in  L_{c}$}}. Consequently, we will prove the relation  {{\small $H_{ij}= H^*_{ij}$}} (for  {{\small $s=+$}}) in  Proposition~\ref{prop:end-stage1-new}.  In the rest of the paper,  we  sometimes drop $s$ and simply state our results for  the case {{\small $s=+$}}. Note however  that $\pi$ is identifiable  only up to a reversal, so it is possible that the reverse orientation holds in the assumption of Lemma~\ref{lem-orientation} (i.e. {{\small $\argmin_{i \in L_c \cup  R_c}    \si_{i}  \in  R_c$}}). In this case, it is possible to show a similar result as in Lemma~\ref{lem-orientation}, where the  {{\small $L_i$}}'s and {{\small $R_i$}}'s are switched, and  the conclusion of Proposition~\ref{prop:end-stage1-new}  becomes {{\small $H_{ij}= -H^*_{ij}$}} (instead of  {{\small $H_{ij}= H^*_{ij}$}}).

\smallskip

\textit{Proof of Proposition~\ref{prop:end-stage1-new}.} As explained above, we  consider   the canonical  `left-right' orientation, and we only prove Proposition~\ref{prop:end-stage1-new}  for this orientation {{\small $(s=+)$}}.
For all {{\small $i\in[n]$}}, any non-empty {{\small $L_i$}} is on the left side  of {{\small $i$}}, and any non-empty  {{\small $R_i$}} is on the right side of {{\small $i$}} (by Lemma~\ref{lem-orientation}). Therefore, every time Algorithm~\ref{algo:SAB} sets {{\small $H_{ij}$}} to {{\small $\pm 1$}}, it is  the correct value {{\small $H_{ij}= H^*_{ij}$}}. Thus, {{\small $H$}} satisfies  {{\small $H_{ij}= H^*_{ij}$}} whenever {{\small $H_{ij}\neq 0$}}. 

It remains to show that {{\small $H_{ij}\neq 0$}} for all {{\small $i,j$}} such that {{\small $|\pi_j - \pi_i| \geq \rho$}}. Observe that for any $i$, either the two sets {{\small $L_i, R_i$}} are non-empty, or exactly one is empty (by Lemma~\ref{lem-orientation}).  Let us consider  the scenario where at least one of {{\small $i$}} and {{\small $j$}} has two non empty sets, say {{\small $L_i \neq \emptyset$}} and {{\small $R_i \neq \emptyset$}}. Then, {{\small $j$}} belongs to {{\small $L_i$}} or {{\small $R_i$}} when {{\small $|\pi_j - \pi_i| \geq  \rho$}},  since {{\small $L_i \cup R_i$}} contains the elements that are {{\small $\rho$}} away from {{\small $i$}} (by Lemma~\ref{lem-orientation}). Thus, by construction of Algorithm~\ref{algo:SAB}, {{\small $H_{ij} \neq 0$}} in this first scenario.   

Let us consider the second scenario where each of $i$ and $j$ has (exactly) one empty set. If the two empty sets are on the same side, e.g.  {{\small $L_i = L_j = \emptyset$}}, then the  sets on the opposite side are necessarily non-empty: {{\small $R_i \neq \emptyset$, $R_j \neq \emptyset$}}. This yields either {{\small $i \in R_j$}} or {{\small $j\in R_i$}} when {{\small $|\pi_j - \pi_i| \geq  \rho$}} (again by Lemma~\ref{lem-orientation}). If the two empty sets are on opposite sides, say {{\small $L_i  = \emptyset$}} and {{\small $R_j  = \emptyset$}}, then  {{\small $|\pi_i - 1|  \leq n/8$}} and {{\small $|\pi_j - n| \leq n/8$}} (by Lemma~\ref{lem-orientation}). Hence, $\pi_j$ is to the right of $\pi_i$, and thus  {{\small $j \in R_i$}} when {{\small $|\pi_j - \pi_i| \geq \rho$}} (by Lemma~\ref{lem-orientation}). We conclude that  {{\small $H_{ij} \neq 0$}} in this second scenario.  The proof of Proposition~\ref{prop:end-stage1-new} is complete.\hfill $\square$


\subsection{Proof of Lemma~\ref{lem-bisection}  in the ideal scenario}
\label{proof-lem-bisection} 
We consider the ideal scenario where {{\small $D$}} equals the ordering distance, i.e., {{\small $D_{kl}=|\pi_k-\pi_l|$}} for all $k,l$. We thus  have {{\small $D\in \lde(\alpha, \beta, \om, r)$}} for the ideal values {{\small $\alpha = \beta = 1$}} and  {{\small $\om= 0$, $r=1$}}. In this case, note that {{\small $\rho := (\de_2+\om)/\al =\de_2$}}. Given {{\small $i\in[n]$}}, the next lemmas describe properties of the connected components of  {{\small $\Gcal_i$}}.
\begin{lem}\label{mini-lem-same-side} 
If {{\small $\de_1 < \de_2$}}, all the nodes of a connected component are on the same side of $i$.
\end{lem}
\textit{Proof.}  If two nodes $k, \ell$ of {{\small $\Gcal_i$}}  are connected  (by an edge), we have  {{\small $|\pi_k-\pi_{\ell}|\leq \de_1$}} and {{\small $|\pi_i-\pi_{\ell}|\vee |\pi_i-\pi_k| \geq \de_2$}}. So $k, \ell$ are on the same side of $i$ (since {{\small $\de_1 < \de_2$}}). Now, if two nodes $k', \ell'$ are in  a same connected component of {{\small $\Gcal_i$}}, there exists a path (of connected nodes) from $k'$ to $\ell'$.  Hence $k'$ and $\ell'$ lie on the same side, since any two consecutive nodes along the path do. \hfill $\square$

\begin{lem}\label{min-lem-include-distant-point}
If {{\small $\de_1 \geq 1$}}, all $k$ such that {{\small $\pi_k \leq \pi_i - \de_2$}} (resp. {{\small $\pi_k \geq \pi_i + \de_2$}}) are in the same connected component.
\end{lem}
\textit{Proof.}
Let {{\small $k, l$}} such that {{\small $\pi_k < \pi_l \leq \pi_i- \de_2$}}. There exist indices {{\small $\ell_0,\ell_1,\ldots, \ell_{|k-l|}$}} such that {{\small $\pi_{\ell_0} = \pi_k$}} and  {{\small $\pi_{\ell_{|k-l|}}= \pi_l$}} and {{\small $\pi_{\ell_{s+1}} - \pi_{\ell_s} =1$}} for all {{\small $s$}}. So  {{\small $|\pi_{\ell_{s+1}} - \pi_{\ell_s}| \leq \de_1$}} (since {{\small $1 \leq \de_1$}}), and  {{\small $\ell_{s}, \ell_{s+1}$}} are connected (by an edge) in {{\small $\Gcal_i$}}. By induction on $s$, we conclude that $k$ and $l$ belong to the same component of {{\small $\Gcal_i$}}.  \hfill $\square$

\begin{lem}\label{min-lem-exit-1or2-compo}
If {{\small $\de_1 \geq 1$}} and {{\small $\de_2 \leq \de_3 \leq n/8$}}, the number of  connected components  that include (at least one)  $k$ such that  {{\small $|\pi_k-\pi_i|\geq \de_3$}}, equals  {{\small $1$}} or {{\small $2$}}.
\end{lem}
\textit{Proof.}  Recall the general fact: any two connected components are either identical or disjoint.  If  a  connected component contains  $k$ such that  {{\small $|\pi_i- \pi_k|\geq \de_3$}}, then for {{\small $\de_3 \geq \de_2$}}, it contains {{\small $\{\pi_k \leq \pi_i - \de_2 \}$}} or {{\small $\{\pi_k \geq \pi_i + \de_2\}$}} by Lemma~\ref{min-lem-include-distant-point}. Hence, {{\small $\Gcal_i$}} has at most  two such components, and at least one since {{\small $\exists k$}} such that  {{\small $|\pi_k-\pi_i|\geq  \de_3$}} (indeed,  {{\small $|1 -\pi_i| \vee |n- \pi_i|  \geq n/8$}}). \hfill $\square$  \vspace{1pt} \smallskip  

\textit{Proof of Lemma~\ref{lem-bisection}} in the ideal scenario where   {{\small $D\in \lde(1, 1, 0, 1)$}}, and {{\small $1 \leq \de_1 < \de_2 \leq \de_3 \leq n/8$}}. It follows from Lemma~\ref{min-lem-exit-1or2-compo} that   {{\small $G_i \neq \emptyset$}} since {{\small $G_i$}} is (by definition) the largest  component containing (at least one)  $g$ such that  {{\small $|\pi_g-\pi_i|\geq \de_3$}}. Then, Lemma~\ref{mini-lem-same-side} gives   {{\small $\pi_{G_i} \subset \{k <\pi_i\}$}} or  {{\small $\pi_{G_i} \subset \{k >\pi_i\}$}}. Since  {{\small $|\pi_g-\pi_i|\geq  \de_3 \geq \de_2$}},   Lemma~\ref{min-lem-include-distant-point} tells us that {{\small $\{k  \leq \pi_i -  \de_2\} \subset  \pi_{G_i}$}} or {{\small $\{k \geq \pi_i +  \de_2\} \subset \pi_{G_i}$}}. Since {{\small $\rho =\de_2$}}, this establishes  property~{{\small 1.}} of Lemma~\ref{lem-bisection}. 

If  {{\small $G_i' \neq \emptyset$}}, we repeat the same argument, and obtain  {{\small $\{k  \leq \pi_i -  \rho\} \subset \pi_{G_i'} \subset \{k <\pi_i\}$}} or  {{\small $\{k \geq \pi_i +  \rho\} \subset  \pi_{G_i'} \subset \{k >\pi_i\}$}}. Then, the two components {{\small $G_i$}} and {{\small $G_i'$}}  are necessarily on opposite sides of $i$ with respect to the ordering $\pi$ (otherwise they would be equal). This gives property~{{\small 2}} of Lemma~\ref{lem-bisection}.

If {{\small $G_i'= \emptyset$}}, then    {{\small $D_{ik}< \de_3$}} for all {{\small $k$}} on opposite side of {{\small $G_i$}} (i.e., {{\small $\pi_k$}} and {{\small $\pi_{G_i}$}} are on  opposite sides of {{\small $\pi_i$}}). Among those {{\small $k$}},   take {{\small $k_0 \in\pi^{-1}\{1,n\}$}}. Then     {{\small $|\pi_i-\pi_{k_0} |  <  \de_3 \leq n/8$}}. Meanwhile, for {{\small $k_0' \in\pi^{-1}\{1,n\}$, $k'_0 \neq k_0$,}} we have {{\small $k_0'\in G_i$}} (by the {{\small 1}} of Lemma~\ref{lem-bisection}). Property~{{\small 3}}  follows. \hfill $\square$


\subsection{Proof of Lemma~\ref{lem-orientation}}
\label{proof:lemma-align} 
If the subroutine {{\small \texttt{Orientation}}} defines the sets {{\small $L_i,R_i$}} by setting either {{\small $(L_i,R_i)=(G_i,G'_i)$}} or {{\small $(L_i,R_i)=(G'_i,G_i)$}}, then Lemma~\ref{lem-bisection} ensures the following. If {{\small $L_i \neq \emptyset$}}, then either {{\small $\{k  \leq \pi_i -  \rho\} \subset  \pi_{L_i} \subset \{k < \pi_i\}$}} or {{\small $\{k \geq \pi_i +  \rho\} \subset \pi_{L_i} \subset \{k > \pi_i\}$}}. If {{\small $L_i = \emptyset$}}, then either {{\small  $\pi_i - 1   <  n/8$}} and {{\small  $n \in \pi_{R_i}$}},  or, {{\small  $n - \pi_i   <  n/8$}}  and {{\small  $1 \in \pi_{R_i}$}}. The same holds for {{\small $R_i$}}. Therefore, it remains to prove that  {{\small  $L_i$ }}  is to the left  of {{\small $i$}}, and  {{\small  $R_i$ }} is to the right of {{\small $i$}}.

Let {{\small $i$}} such that  {{\small $G'_i = \emptyset$}}. In this case, it suffices to show that   {{\small $L_c$}} or {{\small $R_c$}} contains $i$, and this set is   on the opposite  side of  {{\small $G_i$}}. 
Assume for a moment that {{\small $(\pi_c - 1) \wedge (n - \pi_c) \geq n/4$}}. Because either {{\small  $\pi_i -1 < n/8$}}  with {{\small $n \in \pi_{G_i}$}}, or {{\small  $n- \pi_i < n / 8$}}  with {{\small $1 \in \pi_{G_i}$}}, we obtain {{\small $|\pi_c - \pi_i| > n/8 \geq  \rho$}}. Thus  {{\small $i \in L_c \cup R_c$}} (by  Lemma~\ref{lem-bisection}).   If {{\small $i\in L_c$}}, then  {{\small $\pi_i < \pi_c \leq 3n/4$}}, and necessarily {{\small $\pi_i -1 < n/8$}} and  {{\small $n \in \pi_{G_i}$}}. Hence,  {{\small $L_c$}} and {{\small $G_i$}} are on opposite sides. The same argument applies if {{\small $i\in R_c$}}.

We now prove that {{\small $(\pi_c - 1) \wedge (n - \pi_c) \geq n/4$}}.
Let $j$ such that  {{\small $\pi_j = \lfloor n/2 \rfloor$}}. Then for {{\small $\rho \leq n/8$}},  each of {{\small $\pi_{G_j}$}} and {{\small $\pi_{G'_j}$}} contains  {{\small $\{k \leq \pi_j - n/8\}$}} or {{\small $\{k \geq \pi_j + n/8\}$}} (by Lemma~\ref{lem-bisection}). Consequently, {{\small $|G_j| \wedge |G'_j| \geq  n/2  - n/8 -2 \geq n/4$}} for {{\small $n \geq 16$}}. Hence {{\small $|G_c| \wedge |G'_c| \geq n/4$}}. Since {{\small $\pi_{G_c} \subset \{k <\pi_c\}$}} or  {{\small $\pi_{G_c} \subset \{k >\pi_c\}$}}, and the same holds for  {{\small $G'_c$}} (by Lemma~\ref{lem-bisection}),  we get {{\small $(\pi_c - 1) \wedge (n - \pi_c) \geq n/4$}}.

Let $i$ such that  {{\small $G'_i \neq \emptyset$}}. We now show that one of the sets     {{\small $L_c$}} or {{\small $R_c$}} has an empty intersection with one of {{\small $G_i$}} or {{\small $G_i'$}}, and that these two sets lie on opposite sides. First,  {{\small $G_i$}} and {{\small $G_i'$}} are on opposite sides of $i$ (by  Lemma~\ref{lem-bisection}). If e.g. {{\small $\pi_i < \pi_c$}}, then {{\small $R_c$}} has empty intersection with  the set {{\small $G_i$}} or {{\small $G_i'$}} that is to the left of $i$. Thus, we can readily check that there is always  an empty intersection between  {{\small $G_i$}} or {{\small $G_i'$}} and {{\small $L_c$}} or {{\small $R_c$}}. Observe that only sets lying on opposite sides can have an empty intersection  (since  {{\small $1 \in \pi_{L_c}$}} and {{\small $n \in \pi_{R_c}$}}, and, either {{\small $1 \in \pi_{G_i}$}} and {{\small $n \in \pi_{G_i'}$}}, or the reverse, {{\small $n\in \pi_{G_i}$}} and   {{\small $1 \in \pi_{G_i'}$}}). 



\section{Refined seriation}
\label{appendix:proof-complete-seriation}

This appendix provides the proof of Proposition~\ref{prop:complete-seriation}, which analyzes the performance of the refined seriation step.

We first establish accuracy guarantees for the subroutine {{\small \texttt{Evaluate} \texttt{Comparison}}}, whose pseudocode appears in Appendix~\ref{appendix:pseudocode:refined-seriation}. We then derive Proposition~\ref{prop:complete-seriation} in Appendix~\ref{proof-prop-complete-seriation}.

\paragraph*{Notation.}
In this appendix, we write {{\small $\{\textup{P}(k)\}$}} for the set {{\small $\{k\in[n] : \textup{P}(k)\ \text{holds}\}$}}.  
Given a set {{\small $G$}}, we denote by {{\small $\pi_G$}} the collection {{\small $\{\pi_g\}_{g\in G}$}},  and {{\small $|G|$}} its cardinality.

\subsection{Subroutine performance}
\label{appendix-perf-routine-evaluate}
Lemma~\ref{lem:comparison-evaluation} below shows that {{\small \texttt{Evaluate} \texttt{Comparison}}} performs well, 
if the inputs {{\small $L, R \subset [n]$}} and {{\small $i,j\in [n]$}} satisfy the  following three conditions: \\ (i)  the inclusions 
\begin{equation}\label{conditions-L-R}
    \pi_L  \subset \{k < \pi_i \wedge \pi_j\}\qquad \textrm{and} \qquad  \pi_R \subset \{k >  \pi_i \vee \pi_j\}
\end{equation}
where {{\small $L$}} (resp. {{\small $R$}}) is, with respect to the ordering {{\small $\pi$}},  on the left (resp. right) side  of {{\small $i,j$}} (here, we used the notation {{\small $\{k < \pi_i \wedge \pi_j\}:=\{k: k < \pi_i \wedge \pi_j\}$}});\\
(ii) the event
\begin{equation}\label{event:gaussian-noise-individual}
   \mathcal{E}(L, R) = \left\{\,  \max_{ B\in\{L, R \} } \ \,\frac{1}{\sqrt{2 |B|}} \ \,  \Big| \sum_{\ell \in B} (E_{j\ell}-E_{i \ell})\Big| \ \,\leq \ \,  \sqrt{12 \log n}\enspace \right\} 
\end{equation}
where we recall that {{\small $E$}}  stands for the noise in~\eqref{eq:model}; 
  (iii)  the lower bound 
\begin{equation}\label{condition:strong-robinson-L-R}
\sum_{k \in L} H^*_{ij} (F_{\pi_j \pi_k}-F_{\pi_i \pi_k})   \ \vee \  \sum_{k \in R} H^*_{ij} (F_{\pi_i \pi_k}-F_{\pi_j \pi_k})\ \ \geq  \ \  \der  |\pi_i-\pi_j|  
\end{equation}
which can be seen as  a similar but {{\small $\pi$}}-permuted version of the average {{\small $\ell_1$}} lower bound in  {{\small $F \in \mathcal{AL}(\alpha, \beta,r,r')$}} (Definition~\ref{defi:PBL}) with a parameter {{\small $\der$}} instead of {{\small $\al$}}. Recall that  {{\small $H^*_{ij} = 1 - 2 \mathbf{1}_{\si_i < \si_j}$}}.  

\begin{lem}
\label{lem:comparison-evaluation}
Let {{\small $i, j \in[n]$}} and {{\small $L, R\subset [n]$}} such that    (\ref{conditions-L-R}-\ref{event:gaussian-noise-individual}-\ref{condition:strong-robinson-L-R}) are fulfilled. Then  the output  {{\small $H_{ij}$}}  of
{{\small \texttt{Evaluate} \texttt{Comparison}}}  satisfies
\begin{equation*}
      H_{ij}= H^*_{ij}  \qquad  \textrm{ if}  \quad |\pi_i-\pi_j| \geq 10 \frac{\va}{\der}\sqrt{n \log n } \enspace.
\end{equation*} 
\end{lem}

Thus, for {{\small $i,j$}}   at   distance at least {{\small $10 (\va/\der)\sqrt{n \log n}$}}, the estimate {{\small $H_{ij}$}} equals  {{\small $H^*_{ij}$}}.  Lemma~\ref{lem:comparison-evaluation} is only for the orientation~\eqref{conditions-L-R}, however   due to the non-identifiability of the  `left-right' orientation, it is possible that {{\small $L,R$}} satisfies the reverse orientation, where  {{\small $\pi_R  \subset \{k < \pi_i \wedge \pi_j\}$ and $\pi_L \subset \{k >  \pi_i \vee \pi_j\}$}}, instead of (\ref{conditions-L-R}).  In this case, it is straightforward  to adapt the proof of Lemma~\ref{lem:comparison-evaluation}, and obtain a similar conclusion  with  {{\small $H_{ij}= -H^*_{ij}$}}  instead. 

\smallskip 

\textit{Proof of Lemma~\ref{lem:comparison-evaluation}.}
\label{proof-lem-single-evaluation} 
Let {{\small $i,j \in [n]$}} such that   {{\small $|\pi_i-\pi_j| \geq 10(\va / \der) \sqrt{n \log n}$}}. Assume   that {{\small $\pi_i < \pi_j$}}, which yields {{\small $H^*_{ij}=-1$}}. Recall that  {{\small $l= \sum_{k\in L} A_{ik}- A_{jk}$}}  \ and \  {{\small $r= \sum_{k\in R} A_{ik}- A_{jk}$}}. 
We first show that
\begin{equation}\label{l-test:lower-bound}
    l >  -5 \va \sqrt{n \log n} \enspace.
\end{equation}
Then, observe that the output of {{\small \texttt{Evaluate} \texttt{Comparison}}}  is correct, i.e. {{\small $H_{ij}= -1$}},  if and only if 
\begin{equation}\label{l-r-tests}
     l \geq 5 \sigma \sqrt{n \log n}  \qquad 
    \textup{ or } \qquad  \  r \leq - 5 \sigma \sqrt{n \log n} \enspace.
\end{equation}
So all we need to prove is  \eqref{l-test:lower-bound} and \eqref{l-r-tests}. 

$\circ$ \text{Proof of \eqref{l-test:lower-bound}.} 
Conditionally on the event {{\small $\mathcal{E}(L, R)$}} in~\eqref{event:gaussian-noise-individual}  we have
\begin{equation*}
    \Big| \sum_{\ell \in L} (E_{j\ell}-E_{i \ell})\Big| \ \,\leq \ \,  \sqrt{2 |L|} \sqrt{12 \log n} \ \ < \ \ 5 \sqrt{n \log n } \enspace,
\end{equation*}
since {{\small $|L|  \leq n$}} and {{\small $\sqrt{2} \sqrt{12} < 5 $}}. 
Therefore, on {{\small $\mathcal{E}(L, R)$}} we have
\begin{align}\label{last-eq:proof-new}
    l \ = \ \sum_{k \in L} A_{i k}-A_{j k} \ \ &= \ \  \sum_{k \in L} (F_{\si_i \si_k} - F_{\si_j \si_k}) \ \ + \ \ \sigma  \sum_{k \in L} (E_{ik} -E_{jk}) \nonumber \\
     \ \ &> \ \  \sum_{k \in L} (F_{\si_i \si_k} - F_{\si_j \si_k}) \ \ - \ \  5\va \sqrt{n\log n} \enspace.
\end{align}
Meanwhile, the Robinson shape of {{\small $F$}} tells us that {{\small $F_{\si_i \si_k} - F_{\si_j \si_k} \geq 0$}} for  {{\small $\pi_i < \pi_j$}} and  {{\small $\pi_k \in \pi_L \subset \{k < \pi_i \}$}}, where the inclusion comes from~\eqref{conditions-L-R}. Plugging this into~\eqref{last-eq:proof-new} we obtain {{\small $l >  -5 \va \sqrt{n \log n}$}}, which gives \eqref{l-test:lower-bound}. 

$\circ$ \text{Proof of \eqref{l-r-tests}.}  The assumption \eqref{condition:strong-robinson-L-R} with {{\small $H^*_{ij}=-1$}} gives 
\begin{equation*}
     \sum_{k \in L} (F_{\pi_i\pi_k}-F_{\pi_j \pi_k})   \vee  \sum_{k \in R} (F_{\pi_j \pi_k}-F_{\pi_i \pi_k})\ \ \geq  \ \   \der  |\pi_i-\pi_j| \enspace.
\end{equation*}
At least one of the two sums  satisfies this lower bound. If this is {{\small  $\sum_L$}}, then we plug this lower bound into~\eqref{last-eq:proof-new} and obtain {{\small $l  >   \der |\si_i  - \si_j | -  5\va \sqrt{n\log n}$}}. Therefore, {{\small  $l \geq 5 \sigma \sqrt{n \log n}$}} since  {{\small $|\pi_i-\pi_j| \geq 10(\va /  \der) \sqrt{n \log n}$}}.  If this is  {{\small $\sum_R$}}, we can use the same reasoning, and obtain {{\small $r \leq  -  \der |\si_i  - \si_j | +   5\va \sqrt{n\log n}$, }} which yields {{\small $r \leq -5 \sigma \sqrt{n \log n}$}}. This  completes the proof of~\eqref{l-r-tests}. 

The above proof is for the case   {{\small  $\pi_i < \pi_j$}}. The symmetric case   $\pi_i > \pi_j$ follows by the same argument.  The proof of Lemma~\ref{lem:comparison-evaluation} is complete. \hfill $\square$





\subsection{Proof of Proposition~\ref{prop:complete-seriation}}
\label{proof-prop-complete-seriation}

To prove Proposition~\ref{prop:complete-seriation}, we apply Lemma~\ref{lem:comparison-evaluation} to the comparison sets 
{{\small $(\tilde L_{ij},\tilde R_{ij})$}}, whose construction is described in Appendix~\ref{appendix:pseudocode:refined-seriation}.  
The main task is thus to verify that these sets satisfy the three conditions (\ref{conditions-L-R}), \eqref{event:gaussian-noise-individual} and (\ref{condition:strong-robinson-L-R}) of Lemma~\ref{lem:comparison-evaluation}.

Throughout the proof we assume {{\small $s=+$}} in Assumption~\eqref{prop:input-H}, so that the input matrix {{\small $H$}} is aligned with {{\small $H^*$}}. The treatment of the reversed orientation is  analogous and omitted for brevity. In the remainder of this appendix, we fix a pair {{\small $i<j$}} with {{\small $H_{ij}=0$}}.


\subsubsection*{Proof of condition (\ref{conditions-L-R}).}
We denote by {{\small $L^*_{ij}$}} the set of indices lying to the left of {{\small $i,j$}} in the ordering {{\small $\pi$}}:
\[
L^*_{ij}
:= \{\,k : H^*_{ik}=H^*_{jk}=1\,\}
= \{\,k : \pi_k < \pi_i \wedge \pi_j\,\}.
\]
Since {{\small $L_{ij}=\{k: H_{ik}=H_{jk}=1\}$}} and {{\small $H$}} is correct on its support by~(\ref{prop:input-H}), we have
\[
L_{ij} \subset L^*_{ij}.
\]
By construction of the modified comparison sets, we also have {{\small $\tilde L_{ij}\subset L_{ij}$}}.  
Combining the two inclusions gives
\[
\tilde L_{ij} \subset L^*_{ij},
\]
which is precisely the required condition for the left side in~(\ref{conditions-L-R}).

Repeating the same argument for {{\small $\tilde R_{ij}$}}, we obtain
\[
\tilde R_{ij} \subset R^*_{ij},
\qquad \quad 
R^*_{ij}:=\{\,k : \pi_k > \pi_i \vee \pi_j\,\}.
\]
Thus condition~(\ref{conditions-L-R}) holds.


\subsubsection*{Distance estimates and ordering distance.}
Before proving conditions~(\ref{event:gaussian-noise-individual}–\ref{condition:strong-robinson-L-R}), we collect several relations between our distance estimates and the ordering distance.
Recall that {{\small $p_{ij}\in\argmin_{p\in S^{t_{ij}}}D_{ip}$}}, where {{\small $t_{ij}\in[3]$}} is chosen so that {{\small $i,j\notin S^{t_{ij}}$}}; its construction in the refinement step is given in Appendix~\ref{appendix:pseudocode:refined-seriation}.  
We also recall that {{\small $c_0=1/32$}}, and  {{\small $\omega_n = C_{\beta \sigma} n^{3/4}(\log n)^{1/4}$}} as defined in~\eqref{distance-estimate-simple-error-constant}.

\begin{lem}
\label{implication:delta=>rho}
    Let {{\small $\cE_0$}} be the event where  the following three properties hold uniformly for all {{\small $i<j$}} such that  {{\small $H_{ij}=0$}}. 
    \begin{enumerate}
        \item {{\small $|\pi_{p_{ij}} - \pi_{i}| \ \ \leq \ \ \al^{-1}\pa{ 4\beta \sqrt{n \log n} + 2 \om_n}$}}
        \item   {{\small $\widehat D^{t_{ij}}_{p_{ij}k}   \geq \delta_4$}} $\Longrightarrow$  {{\small $|\pi_i - \pi_k| \wedge |\pi_j - \pi_k| \  \geq  \  \rho$}} 
        \item      {{\small $\{k: \widehat D^{t_{ij}}_{p_{ij} k} \geq \de_4\} \ \ \supset  \ \ \{k \in S^{t_{ij}}: (\pi_k < \pi_i \wedge \pi_j -c_0 n)  \textrm{ or }  (\pi_k > \pi_i \vee \pi_j +c_0 n) \}$}} 
    \end{enumerate}
    If   {{\small $D^* \in \lde(\alpha,\beta,0,r)$}} and   {{\small $D \in \lde(\alpha,\beta,\om_n,r)$}} and {{\small $H, \de_4$}} as in (\ref{prop:input-H}-\ref{constraint-2nd-seriation}), then  {{\small $\P\ac{\cE_0}\geq 1- \tfrac{2}{n^3}$}}. 
\end{lem}

The assumptions of Lemma~\ref{implication:delta=>rho} are indeed satisfied under the conditions of Proposition~\ref{prop:complete-seriation}.
Since {{\small $F\in\mathcal{AL}(\alpha,\beta,r,r)$}}, Lemma~\ref{lem:AL-implies-LDE} ensures that {{\small $D^*$}} belongs to {{\small $\lde(\alpha,\beta,0,r)$}}.
The requirements {{\small $D\in\lde(\alpha,\beta,\omega_n,r)$}} and~\eqref{prop:input-H} on the inputs {{\small $(D,H)$}} are already imposed in Proposition~\ref{prop:complete-seriation}.
Finally, the constraint~\eqref{constraint-2nd-seriation} on {{\small $\delta_4$}} is a  detailed version of the simplified condition~\eqref{constraint-2nd-seriation-simplified-new} stated in Proposition~\ref{prop:complete-seriation}.
Hence the conditions of Lemma~\ref{implication:delta=>rho} hold, and {{\small $\P\ac{\cE_0^c}\leq 2/n^3$}}.

The first property in Lemma~\ref{implication:delta=>rho} shows that {{\small $p_{ij}$}} lies within {{\small $O(\om_n)$}} of {{\small $i$}} in ordering distance, so {{\small $p_{ij}$}} is a good proxy for {{\small $i$}}.  
The second property implies that any {{\small $k$}} satisfying {{\small $\widehat D^{t_{ij}}_{p_{ij}k}\ge\delta_4$}} -- and therefore any {{\small $k\in\tilde L_{ij}\cup\tilde R_{ij}$}} -- is at ordering distance at least {{\small $\rho$}} from both {{\small $i$}} and {{\small $j$}}.  
The third property ensures that the threshold {{\small $\delta_4$}} does not exclude all indices far from {{\small $\{i,j\}$}}, so that at least one of the sets {{\small $\tilde L_{ij}$}} or {{\small $\tilde R_{ij}$}} contains many indices with {{\small $\pi_k<\pi_i\wedge\pi_j-c_0n$}} or {{\small $\pi_k>\pi_i\vee\pi_j+c_0n$}}.  
In other words, properties (ii) and (iii) mean that {{\small $\delta_4$}} is simultaneously large enough to discard all indices within ordering distance strictly smaller than {{\small $\rho$}} from {{\small $\{i,j\}$}}, and small enough to retain a substantial number of indices far to the left or to the right of {{\small $\{i,j\}$}}.
The proof of Lemma~\ref{implication:delta=>rho} is deferred to  Appendix~\ref{appendix-proof-lem-relations-dist-estimates-&-ordering-dist}.


\subsubsection*{Proof of condition~(\ref{event:gaussian-noise-individual}) \& data dependence.}
The definition of  {{\small $\tilde L_{ij}$}}  in  Appendix~\ref{appendix:pseudocode:refined-seriation}  involves the set {{\small $L_{ij}$}}, which depends on  {{\small $( A_i, A_j)$}}, and might therefore suggest a nontrivial dependence of {{\small $\tilde L_{ij}$}} on the data. However, this intuition is incorrect. 
Define
\begin{equation}
\label{L'ij-definition}
    L_{ij}' \ \ := \ \ \{  k\in L^*_{ij}  \cap S^{t_{ij}} \ \textrm{such that}  \ \widehat D^{t_{ij}}_{p_{ij}k} \geq \delta_4 \}  ,
\end{equation}
which coincides with the definition of {{\small $\tilde L_{ij}$}} except that {{\small $L_{ij}$}} is replaced by its population analogue {{\small $L_{ij}^*=\{k:\pi_k<\pi_i\wedge\pi_j\}$}}.
Lemma~\ref{tilde-Rij:reduce-dependence} shows that {{\small $\tilde L_{ij}=L_{ij}'$}} on the event {{\small $\cE_0$}}, so {{\small $\tilde L_{ij}$}} does not inherit any additional dependence on {{\small $(A_i,A_j)$}} through {{\small $L_{ij}$}}.
An analogous statement holds for {{\small $\tilde R_{ij}$}}, where we set
\begin{equation}
\label{def-R'-ij}
     R_{ij}' = \{  k\in R^*_{ij}  \cap S^{t_{ij}} \ \textrm{such that}    \ \widehat D^{t_{ij}}_{p_{ij}k} \geq \delta_4  \}\enspace.
\end{equation}

\begin{lem}
\label{tilde-Rij:reduce-dependence}
Assume {{\small $H$}} and {{\small $\delta_4$}} satisfy (\ref{prop:input-H}) and (\ref{constraint-2nd-seriation}) respectively.  
Then, on the event {{\small $\cE_0$}} of Lemma~\ref{implication:delta=>rho}, we have {{\small $\tilde L_{ij}=L_{ij}'$}} and {{\small $\tilde R_{ij}=R_{ij}'$}} for all {{\small $i<j$}} with {{\small $H_{ij}=0$}}.
\end{lem}

The proof of Lemma~\ref{tilde-Rij:reduce-dependence} is given in Appendix~\ref{appendix-remaining-proofs}. 
Even after this reduction, {{\small $\tilde L_{ij}$}} still depends on {{\small $(A_i,A_j)$}} through the random index {{\small $p_{ij}$}}. 
This dependence is mild, since {{\small $p_{ij}$}} can take at most {{\small $n$}} possible values.
Combining this observation with Lemma~\ref{tilde-Rij:reduce-dependence}, we are able to obtain a uniform control of the noise event~\eqref{event:gaussian-noise-individual} over all pairs {{\small $\{i,j\}$}} with {{\small $H_{ij}=0$}}. This is formally stated in the next lemma.

\begin{lem}
\label{noise-uniform-control}
Let {{\small $\mathcal{E}_0'$}} denote the event {{\small $\bigcap_{\{i<j : H_{ij}=0\}} \mathcal{E}(\tilde L_{ij},\tilde R_{ij})$}}, where {{\small $\mathcal{E}(\cdot,\cdot)$}} is defined in~(\ref{event:gaussian-noise-individual}). Under the assumptions of Lemma~\ref{tilde-Rij:reduce-dependence},
\[
  \P\ac{ \mathcal{E}_0' } \ \ge\ 1 - \frac{4}{n^3} \;.
\]
\end{lem}

The proof of Lemma~\ref{noise-uniform-control} is given in Appendix~\ref{proof-lem-noise-control}. This establishes condition~(\ref{event:gaussian-noise-individual}) for 
{{\small $(\tilde L_{ij},\tilde R_{ij})$}}.


\subsubsection*{Proof of condition~(\ref{condition:strong-robinson-L-R}).} 

We now verify the {{\small $\ell_1$}} separation condition~(\ref{condition:strong-robinson-L-R}).  Two effects must be controlled.  
The subsampling used in the refinement step removes a fixed proportion of indices, and the construction of {{\small $(\tilde L_{ij},\tilde R_{ij})$}} further discards all indices too close to {{\small $\{i,j\}$}} in ordering distance.  
We must therefore ensure that the remaining indices still preserve the {{\small $\ell_1$}} lower bound provided by the {{\small $\mathcal{AL}(\alpha,\beta,r,r')$}} assumption.

The following lemma shows that this is the case: with high probability, the sets {{\small $(\tilde L_{ij},\tilde R_{ij})$}} retain an average {{\small $\ell_1$}} separation of order {{\small $\alpha|\pi_i-\pi_j|$}}, up to a factor {{\small $1/4$}}, uniformly over all pairs with {{\small $H_{ij}=0$}}.

\begin{lem}
\label{lem:strong-robinson-after-subsampling}
Let {{\small $\mathcal{E}_0''$}} denote the event
\[
   \mathcal{E}_0'' 
   := 
   \bigcap_{\{i<j : H_{ij}=0\}} 
   \{\tilde L_{ij},\tilde R_{ij} \text{ satisfy condition~(\ref{condition:strong-robinson-L-R}) with } \der=\alpha/4\}.
\]
Under the assumptions of Lemma~\ref{tilde-Rij:reduce-dependence} and  {{\small $F\in\mathcal{AL}(\alpha,\beta,r,r)$}}, 
\[
   \P\ac{\mathcal{E}_0''} \ \ge\ 1 - \frac{3}{n^3} .
\]
\end{lem}

The proof of Lemma~\ref{lem:strong-robinson-after-subsampling} is in Appendix~\ref{appendix-remaining-proofs}.


\subsubsection*{Application of Lemma~\ref{lem:comparison-evaluation} and proof of Proposition~\ref{prop:complete-seriation}.} 
For all {{\small $i < j$}} such that {{\small $H_{ij}=0$}}, Algorithm~\ref{algo:complete-seriation} calls {{\small \texttt{Evaluate} \texttt{Comparison}}} on the pair {{\small $(\tilde L_{ij},\tilde R_{ij})$}} and returns {{\small $\widetilde H_{ij}$}}.  
Combining Lemmas~\ref{implication:delta=>rho}, \ref{noise-uniform-control} and~\ref{lem:strong-robinson-after-subsampling} and applying a union bound, we obtain
\[
\P\ac{\cE_0 \cap \cE_0' \cap \cE_0''} \ \ge\ 1 - \frac{9}{n^3}.
\]
Thus, with probability at least {{\small $1-9 /n^3$}}, the assumptions of Lemma~\ref{lem:comparison-evaluation} hold uniformly for all {{\small $i<j$}} such that {{\small $H_{ij}=0$}}, with {{\small $\der=\alpha/4$}}.  Lemma~\ref{lem:comparison-evaluation} then yields
\[
   \widetilde H_{ij}=H^*_{ij}, 
   \qquad 
   \forall\, i<j \text{ such that } H_{ij}=0
   \text{ and } |\pi_i-\pi_j|\ge 40\frac{\va}{\al}\sqrt{n\log n}.
\]
By symmetry of {{\small $\widetilde H$}} and {{\small $H^*$}}, the same holds for {{\small $i>j$}}. Since the supports of {{\small $\widetilde H$}} and {{\small $H$}} are disjoint and 
{{\small $H$}} is correct on its support, we obtain
\[
   H_{ij} + \widetilde H_{ij} = H^*_{ij},
   \qquad 
   \forall\, i,j \text{ such that } 
   |\pi_i-\pi_j|\ge 40\frac{\va}{\al}\sqrt{n\log n},
\]
with probability at least {{\small $1-9/n^3$}}. Therefore, {{\small $\widehat H := H+\tilde H$}} has the accuracy stated in Proposition~\ref{prop:complete-seriation}, which concludes the proof.



\subsection{Proofs of Lemmas~\ref{tilde-Rij:reduce-dependence} and~\ref{lem:strong-robinson-after-subsampling}}
\label{appendix-remaining-proofs}

This appendix contains the proofs of Lemmas~\ref{tilde-Rij:reduce-dependence} and~\ref{lem:strong-robinson-after-subsampling}, which were used in Appendix~\ref{proof-prop-complete-seriation}.

\begin{proof}[Proof of Lemma~\ref{tilde-Rij:reduce-dependence} ]
From the proof of~\eqref{conditions-L-R}, we already have {{\small $L_{ij}\subset L_{ij}^*$}}, and therefore {{\small $\tilde L_{ij}\subset L_{ij}'$}} by the definitions of 
{{\small $\tilde L_{ij}$}} and {{\small $L_{ij}'$}}. We now show the reverse inclusion. 
On the event {{\small $\cE_0$}} of Lemma~\ref{implication:delta=>rho}, property~(ii) gives
\[
    \widehat D^{t_{ij}}_{p_{ij}k}\ge\delta_4 
    \ \Longrightarrow\ 
    |\pi_i-\pi_k|\wedge|\pi_j-\pi_k|\ge\rho,
\]
which implies
\[
    L_{ij}^*\cap\{\widehat D^{t_{ij}}_{p_{ij}k}\ge\delta_4\}
    \subset
    \{\pi_k < (\pi_i\wedge\pi_j)-\rho\}.
\]
Since {{\small $L_{ij}=\{k:H_{ik}=H_{jk}=1\}$}} and 
{{\small $H_{st}=H^*_{st}$}} whenever {{\small $|\pi_s-\pi_t|\ge\rho$}} 
by~\eqref{prop:input-H}, we have
\[
    \{\pi_k < (\pi_i\wedge\pi_j)-\rho\} \subset L_{ij}.
\]
Combining the two inclusions yields
\[
    L_{ij}^*\cap\{\widehat D^{t_{ij}}_{p_{ij}k}\ge\delta_4\}\subset L_{ij},
\]
that is, {{\small $L_{ij}'\subset\tilde L_{ij}$}}.
Thus {{\small $\tilde L_{ij}=L_{ij}'$}}. The argument for 
{{\small $\tilde R_{ij}$}} and {{\small $R_{ij}'$}} is identical. This completes the proof of Lemma~\ref{tilde-Rij:reduce-dependence}.
\end{proof}

\smallskip

\begin{proof}[Proof of Lemma~\ref{lem:strong-robinson-after-subsampling}]
Assume for simplicity that {{\small $n/3$}} is an integer.  We first control how the {{\small $\ell_1$}} separation from {{\small $\mathcal{AL}(\alpha,\beta,r,r)$}}  behaves under subsampling. The following lemma shows that a random subset of size {{\small $n/3$}} preserves the {{\small $\ell_1$}} lower bound up to a factor {{\small $1/4$}}.

\begin{lem}
\label{assump-2-stable-by-sampling}
Let {{\small $F\in\mathcal{AL}(\alpha,\beta,r,r')$}} and let {{\small $S\subset[n]$}} be drawn uniformly without replacement with {{\small $|S|=n/3$}}. 
If  {{\small $n \ge C_{\alpha \beta}$}} for some constant {{\small $C_{\alpha \beta}$}},
then with probability at least {{\small $1-2/n^4$}}, the following holds uniformly for all {{\small $i<j$}} such that {{\small $j-i\le nr$}}:
 \begin{equation*}
       \sum_{ k < i - c_0 n}  (F_{ik}-F_{jk})\ \mathbf{1}_{k\in S} \ \ \vee  \  \sum_{ k > j + c_0 n}  (F_{jk}-F_{ik})\ \mathbf{1}_{k\in S} \quad  \geq \quad \frac{\alpha}{4}\,   |i-j| \enspace.
    \end{equation*}
\end{lem}
The proof of Lemma~\ref{assump-2-stable-by-sampling} is  given in Appendix~\ref{section:proof-lemma-E5}.

Consider now the balanced random partition {{\small $(S^1,S^2,S^3)$}} of {{\small $[n]$}}.
Since {{\small $(\pi_{S^1},\pi_{S^2},\pi_{S^3})$}} is also a uniform balanced partition, each {{\small $\pi_{S^t}$}} has the same marginal distribution as a uniform sample of size
{{\small $n/3$}} without replacement. 

Moreover, by assumption {{\small $H_{st}=H^*_{st}\neq0$}} whenever {{\small $|\pi_s-\pi_t|\ge\rho$}}.  
Thus, for any {{\small $i,j$}} with {{\small $H_{ij}=0$}}, we necessarily have {{\small $|\pi_i-\pi_j|<\rho\le rn$}}.  
Hence all pairs {{\small $(i,j)$}} with {{\small $H_{ij}=0$}} fall within the range where Lemma~\ref{assump-2-stable-by-sampling} is applicable.

Applying Lemma~\ref{assump-2-stable-by-sampling} independently to each {{\small $S^t$}}, {{\small $t\in[3]$}}, and taking a union bound yields that, with probability at least {{\small $1-6/n^4$}},
\begin{equation*}
   \sum_{k\in L_{ij}^*(c_0)} H^*_{ij}(F_{\pi_j \pi_k} - F_{\pi_i \pi_k})   \mathbf{1}_{k \in S^{t_{ij}}}  \vee  \sum_{k\in R_{ij}^*(c_0)} H^*_{ij}(F_{\pi_i \pi_k} - F_{\pi_j \pi_k})  \mathbf{1}_{k \in S^{t_{ij}}}   \geq  \ \frac{\alpha}{4}\, |\pi_i - \pi_j|
\end{equation*}
for all {{\small $i<j$}} with {{\small $H_{ij}=0$}}, where  {{\small $L_{ij}^*(c_0)=\{k:\pi_k<\pi_i\wedge\pi_j-c_0 n\}$}} and  {{\small $R_{ij}^*(c_0)=\{k:\pi_k>\pi_i\vee\pi_j+c_0 n\}$}}. In particular, when {{\small $n\geq 6$}}, the inequalities above hold simultaneously with probability at least {{\small $1-1/n^3$}}.

To transfer this bound to the modified sets {{\small $(\tilde L_{ij},\tilde R_{ij})$}}, it suffices to show the inclusions
\begin{equation}\label{incl-F1}
    S^{t_{ij}}\cap L_{ij}^*(c_0) \subset  \tilde L_{ij},
   \qquad
    S^{t_{ij}}\cap R_{ij}^*(c_0) \subset  \tilde R_{ij},
\end{equation}
since all summands are nonnegative.  
Indeed, for any {{\small $k\in L_{ij}^*$}} we have {{\small $H^*_{ij}(F_{\pi_j\pi_k}-F_{\pi_i\pi_k})\ge0$}} because {{\small $F$}} is Robinson and {{\small $\tilde L_{ij}\subset L_{ij}^*$}}, and the analogous property holds on the right side.

Fix {{\small $i<j$}} with {{\small $H_{ij}=0$}}.
On the event {{\small $\mathcal{E}_0$}} of Lemma~\ref{implication:delta=>rho}, Lemma~\ref{tilde-Rij:reduce-dependence} gives {{\small $\tilde R_{ij}=R_{ij}'$}}, with {{\small $R_{ij}'$}} defined in~\eqref{def-R'-ij}.  
Thus the right-hand inclusion in \eqref{incl-F1} is equivalent to
\[
    S^{t_{ij}}\cap R_{ij}^*(c_0)
   \ \subset\ 
   \{k\in R_{ij}^*\cap S^{t_{ij}} : \widehat D^{t_{ij}}_{p_{ij}k}\ge\delta_4\}.
\]
Since {{\small $R_{ij}^*(c_0)\subset R_{ij}^*$}}, this reduces to
\[
    S^{t_{ij}}\cap R_{ij}^*(c_0)
   \ \subset\
     \{k : \widehat D^{t_{ij}}_{p_{ij}k}\ge\delta_4\},
\]
and this inclusion follows directly from property~(iii) of Lemma~\ref{implication:delta=>rho}.  
The proof for the left-hand inclusion in \eqref{incl-F1} is identical.

Thus both inclusions in \eqref{incl-F1} hold on {{\small $\mathcal{E}_0$}}.
Since Lemma~\ref{implication:delta=>rho} yields {{\small $\P(\mathcal{E}_0)\ge 1-\tfrac{2}{n^3}$}}, a union bound gives
\[
   \P\ac{\mathcal{E}_0''}\ \ge\ 1-\frac{3}{n^3},
\]
for the event {{\small $\mathcal{E}_0''$}} defined in Lemma~\ref{lem:strong-robinson-after-subsampling}.
This concludes the proof of Lemma~\ref{lem:strong-robinson-after-subsampling}.
\end{proof}





\section{Proofs of technical lemmas}
\label{appendix:technical-lemmes:olf-sumplement-mat}

\subsection{Proof of Lemma~\ref{lem-last-step-H-perm}}

Here {{\small $\mathbf{1}_{\{\cdot\}}$}} denotes indicator functions, and 
{{\small $\mathbf{1}$}} stands for the all-ones vector {{\small $\mathbf{1}_n \in \mathbb{R}^n$}}.

Let {{\small $i\in[n]$}}. We have {{\small $H^*_{ii}=0$}} and {{\small $H^*_{ik} = 1 -2 \mathbf{1}_{\si_i < \si_k}$}} for  {{\small $k\neq i$}}, so 
\[ \sum_{k=1}^n  H^*_{ik}  \ = \ \sum_{k: \pi_k < \pi_i}  H^*_{ik} \ +  \sum_{k: \pi_k > \pi_i}  H^*_{ik}  \  = \ (\pi_i-1) - (n-\pi_i) = 2\pi_i -(n+1) \]
where we used that {{\small  $\pi$}} is a permutation of {{\small $[n]$}}. Hence,  {{\small  $\pi_i =  (H_i^* \mathbf{1} + n + 1)/2 = \pi^{H^*}_i$}} (where the last equality is true by definition of {{\small  $\pi^{H^*}$}}), and thus {{\small  $\pi =\pi^{H^*}$}}. Let {{\small  $S_i^- = \{k: \pi_i  < \pi_k < \pi_i + n \nu\}$}} and {{\small  $S_i^+ = \{k: \pi_i - n \nu < \pi_k < \pi_i \}$}}. If {{\small  $H$}} satisfies~\eqref{eq-lem-H-pi} for  {{\small  $s=+$}}, then, using the definition {{\small  $\pi_i^H =  (H_i \mathbf{1} + n + 1)/2 $}} we obtain
\[2(\pi_i - \pi_i^H)  \ = \ 2(\pi_i^{H^*} - \pi_i^H)  \ = \ H_i^* \mathbf{1}  -H_i \mathbf{1} \  = \ \sum_{k\in S_i^- }  (-1 - H_{ik}) + \sum_{k\in S_i^+}  (1 - H_{ik}) \]
since  {{\small  $H_{ik} = H^*_{ik}$}} for  {{\small  $k \notin S_i^- \cup S_i^+$}}, and {{\small  $H^*_{ik}=-1$}} for  {{\small  $k \in S_i^-$}}, and {{\small  $H^*_{ik}=1$}} for  {{\small  $k \in S_i^+$}}. The two sums have opposite signs (since {{\small  $H_{ik}\in \{-1, 0, 1\}$}}), hence their magnitudes cannot reinforce each other, and we have
\[2|\pi_i - \pi_i^H| \leq \   \big{|}\sum_{k\in S_i^- }  (-1 - H_{ik}) \big{|} \ \vee \ \big{|}\sum_{k\in S_i^+}  (1 - H_{ik})\big{|} \ \leq \  2(|S_i^-|\vee |S_i^+|) \enspace. \]
Since each $S_i^\pm$ contains at most $\lfloor n\nu \rfloor$ indices,
we have $|S_i^-|\vee |S_i^+|\le \lfloor n\nu \rfloor \le n\nu$.
Therefore
\[
2|\pi_i - \pi_i^H| \le 2n\nu.
\]
This bound holds for any {{\small  $i$}}, hence {{\small $\max_i |\pi_i - \pi_i^H|\leq n \nu$}}. 

The case $s=-$ is handled by symmetry. Specifically,  if  {{\small  $H$}} satisfies~\eqref{eq-lem-H-pi} for {{\small  $s=-$}}, consider the reverse permutation {{\small $\si^{\mathrm{rev}}$}} (defined by  {{\small $\si^{\mathrm{rev}}_i= n + 1 - \si_i$}} for all {{\small $i$}}). As we already saw,  {{\small  $\pi_i =  (H_i^* \mathbf{1} + n + 1)/2$}}; therefore  {{\small  $\pi^{\mathrm{rev}}_i =  (n + 1 - H_i^* \mathbf{1} )/2 =\pi^{-H^*}_i$}} (where the last equality holds  by definition of {{\small  $\pi^{-H^*}$}}). Repeating the same argument as above, we obtain
\[2|\pi^{\mathrm{rev}}_i - \pi_i^H|  \ =  \ |-H_i^* \mathbf{1}  -H_i \mathbf{1}| \ \leq \ 2 n \nu \enspace, \]
and thus  {{\small $\max_i |\si^{\mathrm{rev}}_i - \pi_i^H|\leq n \nu$}}. 

We conclude that {{\small $L_{\max}(\pi, \pi^H) \leq \nu$}} for any {{\small $s\in\{\pm\}$}}. The   proof of Lemma~\ref{lem-last-step-H-perm} is complete.


\subsection{Proof of Lemma~\ref{lem-bisection}}
\label{append-proof-lem-bisection-general-case}

Lemma~\ref{technical-lem-between-dist} gives  useful relations between  {{\small $D_{ij}$}} and  {{\small $|\pi_i-\pi_j|$}}.
\begin{lem}\label{technical-lem-between-dist}
If   {{\small $D\in \lde(\alpha,\beta,\om,r)$}} and {{\small $\de_1, \de_2$}} as in~\eqref{delta-constraint-1rst-seriation}, then for all $i,k, \ell$,  we have 
\begin{align}
D_{k \ell} \leq \delta_{1} & \Longrightarrow |\pi_{k}-\pi_{\ell}|\leq\kappa   \label{cond:C1}\\
|\pi_{i}-\pi_{\ell}|\leq\kappa &\Longrightarrow D_{i\ell}< \delta_{2} \label{cond:C2}\\
D_{i\ell} < \delta_{2}& \Longrightarrow |\pi_{i}-\pi_{\ell}|< \rho\label{cond:C3}\\
  \pi_{\ell} \leq n-1 \ \textup{ and } \ |\pi_{\ell}-\pi_{\ell_c}|= 1 &  \Longrightarrow \  \ D_{\ell \ell_{c}} \leq \delta_{1} \label{cond:C4}
\end{align}
where {{\small  $\ell_c$}}  is defined by {{\small  $\pi_{\ell_c} = \pi_{\ell} +1$}},  and
{{\small $\kappa :=(\de_1+\om)/\alpha$}}, and {{\small $\rho :=(\de_2+\om)/\alpha$}}.
\end{lem}

\textit{Proof of Lemma \ref{technical-lem-between-dist}.} If  {{\small $D_{k \ell}\leq \de_1$}} for {{\small $\de_1\leq rn$}},  then {{\small $ |\pi_{k}-\pi_{\ell}| \leq  \frac{D_{k \ell} +\om}{\alpha}\leq \frac{\de_1+\om}{\alpha} := \kappa$}}.
Similarly, if {{\small $|\pi_{k}-\pi_{\ell}|\leq\kappa$}} for {{\small $\kappa\leq rn$}},  then  {{\small  $D_{i \ell} \leq  \beta \kappa + \om < \delta_{2}$}}. If {{\small $D_{i \ell}<\de_2$}} for {{\small $\de_2\leq rn$}},  then
{{\small $|\pi_{i}-\pi_{\ell}|< \left(\de_2+\om\right)/\alpha :=\rho.$}} Finally,  {{\small $|\pi_{\ell}-\pi_{\ell_c}|= 1$}} and {{\small $1 \leq rn$}}, so
{{\small $D_{\ell \ell_c}\leq  \beta  + \om \leq \de_1$}} (as ensured by the choice of $\delta_1$ in~\eqref{delta-constraint-1rst-seriation}). \hfill $\square$\smallskip 

To prove  Lemma~\ref{lem-bisection}, we use Lemma~\ref{technical-lem-between-dist} and  follow the same steps as in the simple scenario from Appendix~\ref{proof-lem-bisection}.  Below, we  mainly focus on the elements that  differ from Appendix~\ref{proof-lem-bisection}. Given {{\small $i\in[n]$}}, the  next lemmas give properties of the connected components of {{\small $\Gcal_i$}}. 

\begin{lem}\label{mini-lem-same-side-general-case} If {{\small $\de_1, \de_2$}} are as in~\eqref{delta-constraint-1rst-seriation},  then all nodes of a connected component  are on the same side of {{\small $i$}}. 
\end{lem}
\textit{Elements of proof.} If $k,\ell$ are adjacent in $\Gcal_i$ (that is, connected by a single edge), we have {{\small $D_{k\ell}\leq \de_1$}} and {{\small $D_{i k} \vee D_{i \ell} \geq \de_2$}}. Then, (\ref{cond:C1}-\ref{cond:C2}) yield {{\small $|\pi_{k}-\pi_l|\leq \kappa$}} and {{\small $|\pi_{i}-\pi_{\ell}|\vee |\pi_{i}-\pi_{k}|> \kappa$}}. \hfill $\square$

\begin{lem}\label{min-lem-include-distant-point-general-case}
 If {{\small $\de_1, \de_2$}} are as in~\eqref{delta-constraint-1rst-seriation}, then  all {{\small $k$}} such that {{\small $\pi_k \leq \pi_i - \rho$}} (respectively, {{\small $\pi_k \geq \pi_i + \rho$}}) are in a same connected component.
\end{lem}
\textit{Elements of proof.} Let {{\small $k, l$}} such that  {{\small $\pi_k < \pi_l  \leq \pi_{i}-\rho$}}. There exist {{\small $\ell_0,\ell_1,\ldots, \ell_{|k-l|}$}} such that  {{\small $\pi_{\ell_0} = \pi_k$}} and  {{\small $\pi_{\ell_{|k-l|}}= \pi_l$}} and {{\small $\pi_{\ell_{s+1}} - \pi_{\ell_s} =1$}} for all $s$. Then, (\ref{cond:C3}-\ref{cond:C4}) imply {{\small $D_{i \ell_s}  \geq \delta_{2}$}}  and  {{\small $D_{\ell_s \ell_{s+1}} \leq \delta_{1}$}}. Hence, {{\small $\ell_{s},\ell_{s+1}$}} are connected by an edge.  \hfill $\square$

\begin{lem}\label{mini-lem-1or2-compo-general-case}
If {{\small $\de_1, \de_2, \de_3$}} satisfy~\eqref{delta-constraint-1rst-seriation}, then the number of connected components containing at least one vertex {{\small $k$}} with   {{\small $D_{k i}\geq \de_3$}} 
is either $1$ or $2$.  
\end{lem}
\textit{Elements of proof.} ``Not more than {{\small $2$}}'': If {{\small $|\pi_i-\pi_k|<\rho$}}, then {{\small $D_{ki}< \beta\rho+\om \leq \de_3$}}. So {{\small $D_{k i}\geq \de_3$}} implies that {{\small $|\pi_k-\pi_i|\geq \rho$}}. Thus, if a connected component contains {{\small $k$}} such that   {{\small $D_{k i}\geq \de_3$}}, then it contains {{\small $\{k: \pi_k \leq \pi_i - \rho\}$}} or {{\small $\{k: \pi_k \geq \pi_i +  \rho\}$}} by Lemma~\ref{min-lem-include-distant-point-general-case}. ``At least $1$'': If {{\small $D_{k i}< \de_3$}}, then {{\small $|\pi_k-\pi_i| <(\delta_{3}+\om)/\alpha \leq n/8$}}. Take {{\small $k_0$}} such that  
{{\small  $|\pi_{i}-\pi_{k_0}|  \geq n/8$}}. Then, we have  {{\small $D_{k_0 i}\geq \de_3$}}. \hfill $\square$\smallskip 

\textit{Proof of Lemma~\ref{lem-bisection}.} We repeat the same argument as in Appendix~\ref{proof-lem-bisection}. Using the  lemmas above, with {{\small $\de_1, \de_2, \de_3$}} as in~\eqref{delta-constraint-1rst-seriation}, we  obtain  {{\small $G_i \neq \emptyset$}}, as well as the 1 and the 2 of Lemma~\ref{lem-bisection}. 

If {{\small $G_i'= \emptyset$}}, then  {{\small $D_{ik}< \de_3$}} for at least one {{\small $k \in\pi^{-1}\{1,n\}$}}. Then   {{\small $|\pi_i-\pi_{k} |  <(\delta_{3}+\om)/\alpha \leq n/8$}}. This gives the 3 of Lemma~\ref{lem-bisection}. \hfill $\square$


\subsection{Proof of (\ref{eq-scalar-prod-concentration}).} 
\label{appendix:proof-corssed-term-conentration-dist-elementary}

Recall that {{\small $\|F\|_{\infty} := \max_{i\in[n]} \|F_i\|_2$}}, where {{\small $F_i$}} is the {{\small $i$}}-th row of {{\small $F$}}. 
Fix {{\small $i,j\in[n]$}} with {{\small $i<j$}}.  
By the model {{\small $A = F_\pi + \sigma E$}}, we have
\begin{equation}\label{lem-dist-scalar-prod-decompo}
    \langle A_{i}, A_{j} \rangle - \langle F_{\si_i},F_{\si_j}\rangle \ \ = \ \  \va \langle F_{\si_i},E_j\rangle + \va \langle E_i , F_{\si_j} \rangle  + \va^2 \langle E_i , E_j \rangle \enspace.
\end{equation}

\textit{Control of the linear terms.}
The random variables {{\small $(E_{jk})_{k=1}^n$}} are independent, centered, sub-Gaussian, with variance proxies smaller than {{\small $1$}}.
Hence the term {{\small $\langle F_{\pi_i},E_j\rangle = \sum_{k=1}^n F_{\pi_i k} E_{jk}$}} is sub-Gaussian with parameter at most {{\small $C\|F_{\pi_i}\| \le C\|F\|_\infty$}}.
Therefore
\[
\Big| \langle F_{\pi_i},E_j\rangle \Big| \ \le\ C_1 \, \|F\|_{\infty}\, \sqrt{\log n}
\]
with probability at least {{\small $1-1/n^7$}}, for some numerical constant {{\small $C_1$}}.
The same result holds for {{\small $\langle E_i , F_{\pi_j}\rangle$}}.

\textit{Control of the quadratic terms.}
For fixed {{\small $i\neq j$}}, write
\[
\langle E_i,E_j\rangle \;=\; \sum_{k=1}^n E_{ik}E_{jk}
= \sum_{k\notin\{i,j\}} E_{ik}E_{jk} \;+\; E_{ij}(E_{ii} + E_{jj}).
\]
The variables {{\small $E_{ik}$}} and {{\small $E_{jk}$}} are independent, centered, sub-Gaussian with unit variance proxy, so {{\small $X_k:=E_{ik}E_{jk}$}} is centered sub-exponential.
Moreover, the three entries {{\small $E_{ij},E_{ii},E_{jj}$}} are also independent sub-Gaussian, so {{\small $Y:=E_{ij}(E_{ii} + E_{jj})$}} is centered sub-exponential.
All variables in {{\small $\{X_k : k\notin\{i,j\}\}\cup\{Y\}$}} involve disjoint
entries of the family {{\small $\{E_{ab}:a\le b\}$}}, and are therefore mutually independent.
By Bernstein's inequality applied to this sum of {{\small $(n-1)$}} independent sub-exponential variables, we obtain
\[
\big|\langle E_i,E_j\rangle\big|
\ \le\ C_2 \sqrt{n\log n}
\]
with probability at least {{\small $1-1/n^7$}}.

\textit{Conclusion by union bound.}
Combining the three bounds in \eqref{lem-dist-scalar-prod-decompo} and multiplying by the factors {{\small $\sigma$}} or {{\small $\sigma^2$}}, we obtain, with probability at least {{\small $1-3/n^7$}},
\[
\big| \langle A_{i}, A_{j} \rangle - \langle F_{\pi_i},F_{\pi_j}\rangle \big|
\ \le\ C \Big( \sigma \|F\|_{\infty} \sqrt{\log n} \ + \ \sigma^2 \sqrt{n \log n} \Big)
\ = \ C \Big(\sigma + \frac{\|F\|_{\infty}}{\sqrt{n}}\Big)\, \sigma \, \sqrt{n \log n}\, .
\]
Taking a union bound over all pairs {{\small $1\le i<j\le n$}} (at most {{\small $n^2/2$}} pairs) yields the claimed inequality
\[
\max_{1\le i<j\le n}
\big| \langle A_{i}, A_{j} \rangle - \langle F_{\pi_i},F_{\pi_j}\rangle \big|
\ \le\ C \Big(\sigma + \frac{\|F\|_{\infty}}{\sqrt{n}}\Big)\, \sigma \, \sqrt{n \log n}
\]
with probability at least {{\small $1-1/n^4$}}.


\subsection{Proof of Lemma~\ref{implication:delta=>rho}}
\label{appendix-proof-lem-relations-dist-estimates-&-ordering-dist}

We introduce some high probability events.  For any {{\small $t\in[3]$}}, let
\begin{equation}
\label{event-well-spread-St}
    \cE_{S^t}' = \ac{ \max_{i \in[n]} \ \min_{k\in S^t,\ k\neq i } \ |\pi_i - \pi_k| \ \ \leq \ \ 4 \sqrt{n \log n}} 
\end{equation}
be the event where   the ordering distance between {{\small $S^t$}} and  any  {{\small $i\in[n]$}} is small; and let  
\begin{equation}
\label{event:good--dist-estimate-Dt}
    \cE_{S^t}'' \ = \ \Big{\{}  \ \widehat D(S^t) \ \in  \  \lde(\al, \beta, 2\om_n,  r, S^t)\Big{\}} 
\end{equation}
be the event where the distance estimate {{\small $\wD(S^t)$}} satisfies the subset version of Definition~\ref{cond-LB} (local distance equivalence on {{\small $S^t$}} as defined in Appendix~\ref{appendix-extension-dist-estim}),
 where {{\small $\om_n = C_{\beta \sigma} n^{3/4}(\log n)^{1/4}$}}. Let
\[
   \cE_1 := \bigcap_{t\in[3]}
   \big(\cE_{S^t}' \cap \cE_{S^t}''\big).
\]
By Lemma~\ref{lem-high-proba-event-neighbor+dist-estimates},  
{{\small $\P(\cE_1)\ge 1-2/n^3$}}.
We now show that on {{\small $\cE_1$}},  
properties 1, 2, 3 of Lemma~\ref{implication:delta=>rho} hold. Fix {{\small $i<j$}} with {{\small $H_{ij}=0$}}.

\paragraph*{1. Proximity of {{\small $p_{ij}$}} to {{\small $i$}}.}
On {{\small $\cE_1$}}, there exists {{\small $l\in S^{t_{ij}}$}} with  
{{\small $|\pi_l-\pi_i|\le 4\sqrt{n\log n}$}}.  
Using the upper bound in {{\small $\lde(\alpha,\beta,\omega_n,r)$}},  
\[
   D_{li}
   \le 4\beta\sqrt{n\log n} + \omega_n.
\]
Since {{\small $p_{ij}$}} minimizes {{\small $D_{ip}$}} over {{\small $S^{t_{ij}}$}},  
{{\small $D_{p_{ij}i}\le D_{li}$}}.  
The lower bound in {{\small $\lde(\alpha,\beta,\omega_n,r)$}} gives  
{{\small $\alpha|\pi_{p_{ij}}-\pi_i| - \omega_n\le D_{p_{ij}i}$}}.  
Thus
\begin{equation}
\label{proximity:i-k_ij}
   |\pi_{p_{ij}}-\pi_i|
   \ \le\
   \alpha^{-1}\!\big(
      4\beta\sqrt{n\log n} + 2\omega_n
   \big) \qquad := C_n .
\end{equation}


\paragraph*{2. Large {{\small $\widehat D^{t_{ij}}_{p_{ij}k}$}} implies large ordering distance.}
Set {{\small $\widehat D^{t_{ij}} := \widehat D(S^{t_{ij}})$}}.  
On {{\small $\cE_1$}}, we have {{\small $\widehat D^{t_{ij}}\in\lde(\alpha,\beta,2\omega_n,r,S^{t_{ij}})$}} by~\eqref{event:good--dist-estimate-Dt}.  Fix {{\small $k\in S^{t_{ij}}$}} such that {{\small $\widehat D^{t_{ij}}_{p_{ij}k}\ge\delta_4$}}.

We consider two cases.

\emph{Case 1:} {{\small $|\pi_{p_{ij}}-\pi_k| > rn$}}.
By the calibration of {{\small $\delta_4$}} in~\eqref{constraint-2nd-seriation}, we have {{\small $ 2\rho + C_n  \leq r n$}}. Thus
\[
   |\pi_{p_{ij}}-\pi_k|
   \ >\ rn
   \ \ge\ 2\rho + C_n.
\]
Using~\eqref{proximity:i-k_ij} and the triangle inequality, we obtain
\[
   |\pi_k-\pi_i|
   \ \ge\ 
   |\pi_k-\pi_{p_{ij}}| - |\pi_{p_{ij}}-\pi_i|
   \ \ge\ 
   (2\rho + C_n) - C_n
   \ =\ 2\rho.
\]

\emph{Case 2:} {{\small $|\pi_{p_{ij}}-\pi_k| \le rn$}}.
Here the local distance-equivalence property on {{\small $S^{t_{ij}}$}} applies to the pair {{\small $(p_{ij},k)$}}.  
On {{\small $\cE_1$}}, we have
\[
   \widehat D^{t_{ij}}_{p_{ij}k}
   \ \le\ 
   \beta\,|\pi_{p_{ij}}-\pi_k| + 2\omega_n.
\]
Combining this with {{\small $\widehat D^{t_{ij}}_{p_{ij}k}\ge\delta_4$}} and using again the choice of {{\small $\delta_4$}} in~\eqref{constraint-2nd-seriation}, we get
\[
   |\pi_{p_{ij}}-\pi_k|
   \ \ge\ 
   \beta^{-1}(\delta_4 - 2\omega_n)
   \ \ge\ 
   2\rho + C_n.
\]
Together with~\eqref{proximity:i-k_ij}, this yields
\[
   |\pi_k-\pi_i|
   \ \ge\ 
   |\pi_k-\pi_{p_{ij}}| - |\pi_{p_{ij}}-\pi_i|
   \ \ge\ 
   (2\rho + C_n) - C_n
   \ =\ 2\rho.
\]

In both cases, we have shown that {{\small $|\pi_k-\pi_i|\ge 2\rho$}}.  
Moreover, by the {{\small $\rho$}}-accuracy of {{\small $H$}} (condition~\eqref{prop:input-H}), we know that whenever {{\small $H_{ij}=0$}},
\[
   |\pi_i-\pi_j| \ <\ \rho.
\]
Hence, by the triangle inequality,
\[
   |\pi_k-\pi_j|
   \ \ge\ 
   |\pi_k-\pi_i| - |\pi_i-\pi_j|
   \ \ge\ 
   2\rho - \rho
   \ =\ \rho.
\]
This proves property~2.


\paragraph*{3. Indices far from {{\small $\{i,j\}$}} exceed the threshold {{\small $\delta_4$}}.}
We argue by contrapositive: fix {{\small $k\in S^{t_{ij}}$}} with  {{\small $\widehat D^{t_{ij}}_{p_{ij}k} < \delta_4$}}.  
By the choice of {{\small $\delta_4$}} in~\eqref{constraint-2nd-seriation}, we have {{\small $\de_4 \le rn$}}, hence the local distance-equivalence property applies to the pair {{\small $(p_{ij},k)$}}. Thus,  
 on {{\small $\cE_1$}}, 
\[
   |\pi_{p_{ij}}-\pi_k|
   \ \le\
   \alpha^{-1}(\widehat D^{t_{ij}}_{p_{ij}k} + 2\omega_n)
   \ <\
   \ \alpha^{-1}(\delta_4 + 2\omega_n) .
\]
Combining this with~\eqref{proximity:i-k_ij} and the triangle inequality, we obtain  
\[
   |\pi_k-\pi_i|
   \ \le\
   |\pi_{p_{ij}}-\pi_k| + |\pi_{p_{ij}}-\pi_i|
   \ \le\
    \alpha^{-1}(\delta_4 + 2\omega_n)  + C_n.
\]
By the condition  (\ref{constraint-2nd-seriation}) on {{\small $\delta_4$}}, the right-hand side is smaller than {{\small $c_0 n$}}. 

Therefore, for any {{\small $k\in S^{t_{ij}}$}},
\[
   \bigl(\pi_k < \pi_i\wedge\pi_j - c_0 n
        \quad \text{or}\quad
        \pi_k > \pi_i\vee\pi_j + c_0 n\bigr)
   \ \Longrightarrow\
   \widehat D^{t_{ij}}_{p_{ij}k} \ge \delta_4,
\]
which is exactly property~3.  


\smallskip
Thus, we have proved that, on {{\small $\cE_1$}}, all three properties of Lemma~\ref{implication:delta=>rho} hold for every {{\small $i<j$}} with {{\small $H_{ij}=0$}}.  
Since {{\small $\P(\cE_1)\ge 1-2/n^3$}}, the lemma follows. \hfill$\square$




\subsection{Proof of Lemma~\ref{noise-uniform-control}}
\label{proof-lem-noise-control}

Since  {{\small $\P\ac{\cE_0^c}\leq 2/n^3$}} by Lemma~\ref{implication:delta=>rho}, we have  for any event {{\small  $\mathcal{A}$}},
\begin{align*}
   \P\ac{\cA} \ = \   \P\ac{\cA \cap \cE_0^c} + \P\ac{\cA \cap \cE_0} \ \leq \ \frac{2}{n^3} + \P\ac{\cA \cap \cE_0} \enspace.
\end{align*}
Taking {{\small  $\cA = \cup_{\{i<j: H_{ij}=0\}} \cE^c(\tilde L_{ij}, \tilde R_{ij})$}} and then a union bound over the pairs {{\small  $i<j$}}, we obtain
\begin{align}\label{decompo-any-event-E0}
   \P\ac{\cup_{\{i<j: H_{ij}=0\}} \cE^c(\tilde L_{ij}, \tilde R_{ij})} \ \ \leq \ \ \frac{2}{n^3} \ + \  \frac{n^2}{2} \max_{i<j: H_{ij}=0} \P\ac{\cE^c(\tilde L_{ij}, \tilde R_{ij}) \cap \cE_0} \enspace.
\end{align}
By Lemma~\ref{tilde-Rij:reduce-dependence}, we have   {{\small $\cE^c(\tilde L_{ij}, \tilde R_{ij}) \cap \cE_0 = \cE^c(L_{ij}',  R_{ij}') \cap \cE_0$}}. In  definition~\eqref{def-R'-ij} of  {{\small  $R_{ij}'$}}, we see {{\small  $3$}} sources of randomness:  {{\small  $p_{ij}$}},  and  {{\small  $S^{t_{ij}}$}},  and  {{\small  $\widehat D^{t_{ij}}$}} which is fully determined from {{\small  $A^{t_{ij}}$}}. 
We can view {{\small $R_{ij}'$}} as a measurable function of {{\small $(p_{ij}, S^{t_{ij}}, A^{t_{ij}})$}}: 
{{\small $R_{ij}' = F(p_{ij}, S^{t_{ij}}, A^{t_{ij}})$}} for some deterministic function {{\small $F$}}. 
For each {{\small $p\in[n]$}}, define
\begin{equation}\label{deterministic-functions}
R_{ij}'(p) := F(p, S^{t_{ij}}, A^{t_{ij}}), \qquad 
L_{ij}'(p) := \tilde F(p, S^{t_{ij}}, A^{t_{ij}}).
\end{equation}
Decomposing over the disjoint events {{\small $\{p_{ij}=p\}_{p=1}^n$}} and using a union bound, we obtain
\[
\P\ac{\cE^c(L_{ij}', R_{ij}') \cap \cE_0} 
   \le n \max_{p\in[n]} \P\ac{\cE^c(L_{ij}'(p), R_{ij}'(p))} \enspace.
\]
Using the definition~\eqref{event:gaussian-noise-individual} of the event $\cE(L_{ij}'(p),  R_{ij}'(p))$ and then a union bound, we  obtain
\begin{align}\label{union-bound-noise-control}
    \P\ac{\cE^c(\tilde L_{ij}, \tilde R_{ij}) \cap \cE_0} \leq 2n  \max_{\substack{p\in[n] \\ B \in \{L_{ij}'(p),  R_{ij}'(p)\}}} \ \P\ac{\frac{1}{ \sqrt{2\, |B|}} \Big{|} \sum_{\ell \in B} (E_{i\ell}-E_{j\ell}) \Big{|} \ \  \geq \ t_0} 
\end{align}
for {{\small  $t_0 = \sqrt{12  \log n}$}}. (When {{\small $|B|=0$}}, the sum is {{\small $0$}} and the event {{\small $\cE(L,R)$}} is trivially satisfied; hence we only consider {{\small $|B|\ge1$}} above.)

Fix {{\small $p \in[n]$}}. Conditionally on {{\small $(S^{t_{ij}},A^{t_{ij}})$}},  
the set {{\small $B\in\{L'_{ij}(p),R'_{ij}(p)\}$}} is deterministic (by~\eqref{deterministic-functions}).  
By definition of {{\small $(S^{t_{ij}},A^{t_{ij}})$}} in the seriation procedure, and since {{\small $B\subset S^{t_{ij}}$}} (by (\ref{L'ij-definition}--\ref{def-R'-ij})),  
the variables {{\small $\{E_{i\ell},E_{j\ell}:\ell\in B\}$}} are conditionally independent given {{\small $(S^{t_{ij}},A^{t_{ij}})$}},  and are {{\small  $2|B|$}}
mean-zero sub-Gaussian random variables with variance proxies {{\small $\le1$}}.  
By a standard tail bound~\cite[Corollary~1.7]{rigollet2023high}, we have for all {{\small  $t>0$}}, 
\[
\P_{\,|\,S^{t_{ij}},A^{t_{ij}}}\!\left(
\frac{1}{\sqrt{2|B|}}\Big|\sum_{\ell\in B}(E_{i\ell}-E_{j\ell})\Big|\ge t
\right) \le 2 e^{-t^2/2}.
\]
Since this holds for any fixed {{\small $(S^{t_{ij}},A^{t_{ij}})$}}, the same bound holds unconditionally.

Taking {{\small $t=t_0=\sqrt{12\log n}$}} gives {{\small $2e^{-t^2/2}=2n^{-6}$}}.  
With the prefactor {{\small $2n$}} from the union bound in~\eqref{union-bound-noise-control}, we get  
{{\small $\P\ac{\cE^c(\tilde L_{ij},\tilde R_{ij})\cap\cE_0}\le4/n^5$}}.  
Plugging this into~\eqref{decompo-any-event-E0} yields the final bound~{{\small $4/n^3$}}. 



\subsection{Proof of Lemma~\ref{assump-2-stable-by-sampling}}
\label{section:proof-lemma-E5}

Fix a pair {{\small  $(i,j)$}}, with {{\small $i<j$}}. Define {{\small $x_k^{(ij)} = (F_{ik}-F_{jk}) \mathbf{1}_{\{k < i - c_0 n\}}$}}  and {{\small $\tilde x_k^{(ij)} = (F_{jk}-F_{ik}) \mathbf{1}_{\{k > j + c_0 n\}}$}} for all {{\small $k\in[n]$}}. 
Since {{\small $k<i < j$}} on the support of {{\small $x^{(ij)}$}} and {{\small $F$}} is Robinson, we have 
{{\small $F_{ik} > F_{jk}$}} and thus {{\small $x_k^{(ij)} >0$}}. 
Similarly, on the support of {{\small $\tilde x^{(ij)}$}} we have {{\small $k>j > i$}} and hence {{\small $F_{jk} > F_{ik}$}}, so {{\small $\tilde x_k^{(ij)} > 0$}}.

Since  {{\small $F \in \mathcal{AL}(\alpha, \beta,r,r')$}} (Definition~\ref{defi:PBL}),  we know that at least one of the following inequalities holds:
\begin{equation}
\label{assumption}
\sum_{k=1}^n x_k^{(ij)} \ \geq\  \alpha |i-j|,  \quad \textup{or} \qquad \sum_{k=1}^n \tilde x_k^{(ij)} \ \geq\  \alpha |i-j| .
\end{equation}
Assume that the first inequality  with {{\small $x_k^{(ij)}$}}  holds.

Set {{\small $C_0 =192$}}. By the local {{\small $\ell_2$}} upper bound in {{\small $F \in \mathcal{AL}(\alpha, \beta,r,r')$}}, {{\small $\|F_i-F_j\|\le \beta\,|i-j|/\sqrt{n}$}}.
Hence
\begin{equation}\label{eq:Psi-bern-proof}
 C_0 \ \|F_i-F_j\|\,\log n
\ \le\ C_0 \ \beta\,\frac{|i-j|}{\sqrt{n}}\,\log n \ \le \ \alpha |i-j|,
\end{equation}
where the last inequality holds when {{\small $n$}} is larger than some constant  {{\small $C_{\alpha \beta}$}}, so that {{\small $C_0\beta\,\tfrac{\log n}{\sqrt{n}}\le \alpha$}}.

Using the notation {{\small $X_k^{(ij)} := x_k^{(ij)} \mathbf{1}_{\{k\in S\}}$}} for all {{\small $k\in[n]$}}, 
we apply a Bernstein–type inequality for sampling without replacement (Lemma~\ref{bernstein_sample-no-replacement})  to the population {{\small $\{x_k^{(ij)}:k\in[n]\}$}} with 
population mean {{\small $\mu_{ij} = \tfrac{1}{n}\sum_{k=1}^n x_k^{(ij)}$}} 
and population variance {{\small $\sigma^2_{ij} = \tfrac{1}{n}\sum_{k=1}^n (x_k^{(ij)}-\mu_{ij})^2$}}. Then,
 with probability at least {{\small $1-2\delta$}}, we obtain 
\[
\frac{1}{|S|}\sum_{k\in S} x_k^{(ij)} 
\;\ge\; \mu_{ij} \; - \; 
\sigma_{ij}\,\sqrt{\frac{2\,\log(1/\delta)}{|S|}}
\;-\;
\frac{7}{3}\,\frac{(\max_k x_k^{(ij)})\,\log(1/\delta)}{|S|}\, .
\]
Since {{\small $|S|=n/3$}}, multiplying both sides by {{\small $|S|$}} and taking {{\small $\delta = 1/n^6$}} gives
\[
\sum_{k\in S}x_k^{(ij)}
\ \ge\
\frac{1}{3}\sum_{k=1}^n x_k^{(ij)}
\ -\ 
2 \sigma_{ij} \sqrt{n \log n }\, -  14 \max_k x_k^{(ij)} \log n,
\]
with probability at least {{\small $1 - 2/n^6$.}} 

Using that {{\small $\sigma^2_{ij} \le n^{-1}\|F_i-F_j\|^2$}} and {{\small $\max_k x_k^{(ij)}\le \|F_i-F_j\|$}}, we get
\[
\sum_{k\in S}x_k^{(ij)}
\ \ge\
\frac{1}{3}\sum_{k=1}^n x_k^{(ij)}
\ -\ 
16 \|F_i-F_j\| \log n.
\]
Combining this inequality with the first inequality in~\eqref{assumption}, and \eqref{eq:Psi-bern-proof}, and using {{\small $C_0\ge 192$}}, we obtain that, with probability at least {{\small $1-2/n^6$}},
\[
\sum_{k\in S} x_k^{(ij)}
\ \ge\
\frac{1}{4}\sum_{k=1}^n x_k^{(ij)}
\ \ge\
\frac{\alpha}{4}\,|i-j|.
\]
Repeating the same argument for {{\small $\tilde x^{(ij)}$}} 
and applying a union bound over both cases and all pairs {{\small $(i,j)$}} 
yields the claimed bound, with probability at least {{\small $1 - 2/n^4$}}:
\[
\forall i,j:\ 
\sum_{k < i - c_0 n} (F_{ik}-F_{jk})\mathbf{1}_{\{k\in S\}}
\ \vee\
\sum_{k > j + c_0 n} (F_{jk}-F_{ik})\mathbf{1}_{\{k\in S\}}
\ \ge\
\frac{\alpha}{4}|i-j| .
\]
This completes the proof of Lemma~\ref{assump-2-stable-by-sampling}.



\section{High probability events}

For brevity, we assume  in this appendix that {{\small $n/3$}} is an integer. Let us start with a simple concentration inequality on the number of  sampled points in the subsets of {{\small $[n]$}}.
\begin{lem}
\label{technical-lem:interval-sample}
    Let  a subset {{\small $I \subset [n]$}}, and a uniform sample {{\small $S$}}  of {{\small $[n]$}} such that  {{\small $|S| = n/3$}}. Then 
    \[  |S \cap I|  \ \ \geq \ \ \frac{|I|}{3} \ - \ \sqrt{n \log n} \]
    with probability at least {{\small $1- (2/n^6)$}}. 
\end{lem}
\textit{Proof of Lemma~\ref{technical-lem:interval-sample}.}  The random number {{\small $|S \cap I|$}} follows the hypergeometric distribution  with parameters {{\small $(N=n/3,\ p=|I|/n,\ n)$}}, so that {{\small $\E[|S\cap I|]=Np=|I|/3$}}. We apply Hoeffding inequality~\eqref{eq_Hoeffding_hypergeom} and obtain 
\[\P\ac{\left| |S \cap I| - \frac{|I|}{3}\right| \ \geq \ \sqrt{\frac{n t}{6}}}  \ \ \leq \ \ 2 e^{-t}\]
for all {{\small $t>0$}}. Taking {{\small $t = 6 \log n$}}  completes  the proof of Lemma~\ref{technical-lem:interval-sample}.\hfill $\square$


\begin{lem}
\label{lem-high-proba-event-neighbor+dist-estimates}
    Recall that {{\small $\cE_1 := \cap_{t\in[3]} (\cE'_{S^t} \cap \cE''_{S^t})$}} where {{\small $\cE'_{S^t}$}}  and  {{\small $\cE''_{S^t}$}} are defined in~\eqref{event-well-spread-St} and~\eqref{event:good--dist-estimate-Dt} respectively. If {{\small $n\geq 8$}} and   {{\small $D^* \in \lde(\al, \beta, 0, r)$}}  then {{\small $\P\ac{\cE_1^c} \leq 2/n^3$}}. 
\end{lem}
\textit{Proof of Lemma~\ref{lem-high-proba-event-neighbor+dist-estimates}.} Fix {{\small $t\in[3]$}}.   We analyze  the events   {{\small $\cE'_{S^t}$}}  and {{\small $\cE''_{S^t}$}} separately.

$\circ$ For  {{\small $\cE'_{S^t}$}} in~\eqref{event-well-spread-St}.  We cover {{\small $[n]$}} with disjoint intervals of the form {{\small $I_k= [a_k, b_k)$}} for some integers {{\small $a_k, b_k$}},   with  cardinal numbers {{\small $3\sqrt{n \log n} + 4 \leq |I_k| \leq 4 \sqrt{n \log n}$}}, for {{\small $n\geq 8$}}. 
Define the  sets {{\small $\tilde I_k= \pi^{-1}(I_k)$}} for all {{\small $k$}}. The union of the {{\small $\tilde I_k$}}'s  clearly   covers {{\small $[n]$}}, and each {{\small $\tilde I_k$}} has   the same cardinal number as the corresponding {{\small $I_k$}}. 
Since {{\small $S^t$}} is uniformly distributed over all subsets of {{\small $[n]$}} of size {{\small $|S^t|=n/3$}}, we can apply  Lemma~\ref{technical-lem:interval-sample}, which yields {{\small $|\tilde I_k \cap S^t | \geq 2 $}}   with  probability at least {{\small $1-(2/n^6)$}} (since {{\small $|\tilde I_k|/3-\sqrt{n\log n}\ge 4/3$}} and {{\small $|\tilde I_k\cap S^t|$}} is integer).  Taking   a union bound over all sets {{\small $\tilde I_k$}}, whose total number is less than {{\small $\sqrt{n}$}}, we obtain for {{\small $n\geq 4$}},
\begin{equation}
\label{cE'_S^t-proba-bound}
    \P\ac{\cE_{S^t}'} \ \ \geq \ \ \P \ac{ \min_k \ \ | \tilde I_k  \cap S^t|   \ \geq \  2 } \ \ \geq \ \ 1 - \frac{2\sqrt{n}}{n^6} \ \geq \ 1 - \frac{1}{n^5} \enspace.
\end{equation} 

$\circ$ For  {{\small $\cE''_{S^t}$}} in~\eqref{event:good--dist-estimate-Dt}.
Denote {{\small  $\eta := \eta(S^t)= \max_{i \in S^t} \min_{k\in S^t, k\neq i}  |\pi_i - \pi_k|$}} as in~\eqref{nearest-neighbor-dist-upper-bound}. Then,  on the event {{\small $\cE_{S^t}'$}} in~\eqref{event-well-spread-St} we have  {{\small $\eta \leq 4 \sqrt{n \log n}$}}, hence, as we saw  in Appendix~\ref{appendix-extension-dist-estim}, the error bound {{\small $\om_n(\eta)$}} in~\eqref{distance-estimate-condition} satisfies  {{\small $\om_n(\eta) \leq \om_n$}}, where {{\small $\om_n$}} is defined in~\eqref{om-n-def}. Denoting the event {{\small $\{\wD(S^t) \in \lde(\al, \beta, \om_n + \om_n(\eta), r, S^t)\}$}} by {{\small $\cE''_{S^t, \eta}$}}, we have {{\small $\wD(S^t) \in \lde(\al, \beta, 2 \om_n, r, S^t)$}} on the event {{\small $\cE_{S^t}' \cap \cE''_{S^t, \eta}$}}, which exactly means that  $\cE''_{S^t} \supset \cE_{S^t}' \cap \cE''_{S^t, \eta}$. It directly follows that   {{\small $\cE_{S^t}' \cap \cE''_{S^t, \eta} \subset  \cE_{S^t}' \cap \cE''_{S^t}$}}, and thus
\begin{equation*}
(\cE'_{S^t} \cap \cE''_{S^t})^c \ \ \subset \ \ (\cE'_{S^t} \cap \cE''_{S^t, \eta})^c \ \   = \ \ (\cE'_{S^t})^c \cup (\cE''_{S^t, \eta})^c \enspace.
\end{equation*}
For {{\small $\cE_1 := \cap_{t\in[3]}(\cE'_{S^t} \cap \cE''_{S^t})$}} we obtain
\begin{equation}
\label{events-decompo}
\P\ac{\cE_1^c}   \leq \P\ac{\cup_{t\in[3]}(\cE'_{S^t})^c} + \P\ac{\cup_{t\in[3]} (\cE''_{S^t, \eta})^c} \leq  3/n^5 + \P\ac{\cup_{t\in[3]} (\cE''_{S^t, \eta})^c} 
\end{equation}
where  we used a union bound and~\eqref{cE'_S^t-proba-bound}.
Conditioning on {{\small $(S^1,S^2,S^3)$}} such that   {{\small $D^*(S^t) \in \lde(\al, \beta, \om_n, r, S^t)$}} for all {{\small $t\in[3]$}},   we apply  Lemma~\ref{lem-hypothesis:dist-equiv-estimator}.  This and a union bound over {{\small $t\in[3]$}} give  that {{\small $\cup_{t\in[3]} (\cE''_{S^t, \eta})^c$}} happens with probability at most {{\small $3/n^4$}}.  By Lemma~\ref{D(S)-in-LDE},  the event {{\small $\cap_{t\in[3]} \{D^*(S^t) \in \lde(\al, \beta, \om_n, r, S^t)\}$}} happens with probability at least {{\small $1-1/n^3$}}. Therefore, without conditioning, we have
\[\P\{\cup_{t\in[3]}(\cE''_{S^t, \eta})^c\} \leq 3/n^4 + 1/n^3 \enspace.\] 
Plugging this into~\eqref{events-decompo},  we obtain {{\small $\P\ac{\cE_1^c}   \leq  3/n^5 +  3/n^4 + 1/n^3$}}, which is bounded by {{\small $2/n^3$}} for {{\small $n\ge6$}}. This  completes the proof  of Lemma~\ref{lem-high-proba-event-neighbor+dist-estimates}. \hfill $\square$


\begin{lem}
\label{D(S)-in-LDE}
If $n\geq 6$ and {{\small $D^* \in \lde(\al, \beta, 0, r)$}}, then   {{\small $\P\ac{\cap_{t\in[3]} \{D^*(S^t) \in \lde(\al, \beta, \om_n , r, S^t)\}} \geq 1-1/n^3$}}.
\end{lem}
\textit{Proof of Lemma~\ref{D(S)-in-LDE}.} 
Fix a pair {{\small $(i,j)$}}. 
Denoting {{\small $x_k^{(ij)} :=n  (F_{ik}-F_{jk})^2$}} for all {{\small $k$}}, we have  {{\small $(D^*_{ij})^2 = \sum_{k=1}^n x_k^{(ij)}$}}. 
Fix {{\small $t\in[3]$}}.
Since the (marginal) distribution of {{\small $S^t$}} is the distribution of a sampling without replacement in {{\small $[n]$}}, we can apply  Hoeffding inequality~\eqref{eq_Hoeffding_sample} to the sum {{\small $\tfrac{|S^t|}{n}(D^*_{ij}(S^t))^2 = \sum_{k\in S^t} x_k^{(ij)}$}}.
This gives, for all {{\small $t>0$}},
\[\P\pa{\bigg{|}\frac{1}{|S^t|} \sum_{k\in S^t} x_k^{(ij)}  - \frac{1}{n}\sum_{k=1}^n x_k^{(ij)} \bigg{|} \, \geq \, t} \leq 2\exp\pa{\frac{-2|S^t| t^2}{(\max_k x_k^{(ij)})^2}}  \leq 2\exp\pa{\frac{-2 t^2}{3 n}}\]
since {{\small $x_k^{(ij)} =n  (F_{ik}-F_{jk})^2 \leq n$}} and {{\small $|S^t| = n/3$}} (where we recall that {{\small $n/3$}} is assumed to be an integer). 

Recalling that {{\small $(D^*_{ij}(S^t))^2 = \tfrac{n}{|S^t|} \sum_{k\in S^t} (F_{ik}-F_{jk})^2$}}, 
we multiply both sides of the inequality by {{\small $n$}}, and take {{\small $t = 3 \sqrt{n \log n}$}} to obtain
\[\P\pa{ \Big{|}(D^*_{ij}(S^t))^2 -  (D^*_{ij})^2 \Big{|} \geq  3 n^{3/2} \sqrt{\log n}} \leq \frac{2}{n^6}\enspace.\]
 Taking a union over all pairs {{\small $(i,j)$}} and then using  the inequality {{\small $|a - b|  \leq \sqrt{|a^2 - b^2|}$}} for any {{\small $a,b \ge 0$}}, we get 
\[\P\pa{\forall i,j: \  \Big{|}D^*_{ij}(S^t) -  D^*_{ij} \Big{|} < \sqrt{3} n^{3/4} (\log n)^{1/4}} \geq 1 - \frac{2}{n^4}\enspace.\]
Since {{\small $D^* \in \lde(\al, \beta, 0, r)$}} and {{\small $\sqrt{3} n^{3/4} \pa{\log n}^{1/4} \leq \om_n$}}, this gives {{\small  $\P\ac{D^*(S^t) \in \lde(\al, \beta, \om_n, r, S^t)} \geq 1-2/n^4$}}. 
The probability bound  
holding for each fixed {{\small $t\in[3]$}}, 
a union bound over the three subsamples gives
\[
\P\!\left(\bigcap_{t\in[3]}\{D^*(S^t)\in \lde(\alpha,\beta,\omega_n,r,S^t)\}\right)
\ \ge\ 1-\frac{6}{n^4}
\ \ge\ 1-\frac{1}{n^3}\,,
\]
{{\small since $n\ge 6$}}.
This completes the proof of Lemma~\ref{D(S)-in-LDE}. \hfill $\square$


\section{Proof of Theorem~\ref{1st:lowerBound}}
\label{appendix:lower-bound}

We recall that the lower bound is proved  in the particular case where {{\small $F$}}  is known and equal to the {{\small $F_\al$}} defined above Theorem~\ref{1st:lowerBound}. It is not difficult to check that for {{\small $\al \in (0,1]$}}, the matrix {{\small $F_\al$}} belongs to {{\small $[0,1]^{n\times n}$}} as in model~\eqref{eq:model}, and to
the  bi-Lipschitz matrices class {{\small $\mathcal{BL}(\alpha,  \beta)$}}  for  any {{\small $\beta \geq  \al$}}. We will establish the lower bound {{\small $(\va/\al)\sqrt{\log(n)/n}$}} 
under the condition  {{\small $\al/\va \geq C_0 \sqrt{\log(n)/n}$}} where {{\small $C_0$}} is a numerical constant. This last condition is satisfied  as soon as {{\small $n \geq C_{\al \sigma}$}} for some constant {{\small $C_{\al \sigma}$}} only depending on the constants {{\small $\al$}} and {{\small $\sigma$}}.

Our minimax lower bound is based on Fano's method as stated below.   We denote  the set of  permutations of {{\small $[n]$}} by {{\small $\bPi_n$}}. For two permutations {{\small $\si$}} and {{\small $\si'$}} in {{\small $\bPi_n$}}, we denote the Kullback-Leibler divergence of {{\small $\P_{(F_\al,\si)}$}} and {{\small $\P_{(F_\al,\si')}$}} by {{\small $KL(\P_{(F_\al,\si)} \| \P_{(F_\al,\si')})$}}.  Given the loss {{\small $L_{\max}$}} in~\eqref{eq:max-loss},  a radius {{\small $\ep>0$}} and a subset {{\small $\mathcal{S} \subset  \bPi_n$,}} the packing number {{\small $\mathcal{M}(\ep, \mathcal{S}, L_{\max})$}} is defined as the largest number of points in {{\small $\mathcal{S}$}} that are  at  least {{\small $\ep$}} away from each other with respect to {{\small $L_{\max}$}}. Below, we state a specific version of Fano's lemma. 
\begin{lem}[\cite{yu1997assouad}]
\label{fano:prop}
For any subset $\mathcal{S}\subset  \bPi_n$, define the Kullback-Leibler diameter of $\mathcal{S}$ as 
\begin{equation*}
    d_{KL}(\mathcal{S})= \underset{\si,\si'\in \mathcal{S}}{\textup{sup}} \ \, KL(\P_{(F_\al,\si)} \| \P_{(F_\al,\si')})\enspace .
\end{equation*}Then, for any estimator $\hat \si$ and  any $\ep>0$, we have
\begin{equation*}
    \underset{\si \in \mathcal{S}}{\textup{sup}} \quad  \P_{(F_\al,\si)}\cro{ L_{\max}(\hat \si, \si) \geq \frac{\ep}{2}} \, \geq  \, 1 - \frac{d_{KL}(\mathcal{S}) + \log(2)}{\log \mathcal{M}(\ep, \mathcal{S}, L_{\max})}\enspace  .
\end{equation*}
\end{lem}
In view of the above proposition, we mainly have to choose a suitable subset {{\small $\mathcal{S}$}}, control its Kullback-Leibler diameter, and get a sharp lower bound of its packing number. A difficulty stems from the fact that the loss {{\small $L_{\max}(\hat \si, \si)$}}  is invariant when reversing  the ordering {{\small  $\si$}}. 

Let {{\small $k :=  C_1 (\va/\al ) \sqrt{n\log n}$}}, for a small enough numerical constant {{\small $C_1 \in (0,1]$}} (which will be set later). To ensure that {{\small $k \leq n/4$}}, we enforce the condition {{\small $\al/\va \geq C_0 \sqrt{\log(n)/n}$}}, with {{\small $C_0:= 4 C_1$}}. 
For simplicity, we assume that {{\small $n/4$}} is an integer.  We introduce {{\small $n/4$}} permutations {{\small $\si^{(s)} \in \bPi_n$, $s=1,\ldots,n/4$}}. For each {{\small $s\in [n/4]$}}, let {{\small $\si^{(s)}$}} be such that
 \begin{align*}
   \forall j \in [n]\setminus{\{ s, s+k\}}\, : \ \, \si^{(s)}_j &= j\,, \   \quad   \textup{ and } \qquad  
    \si^{(s)}_s = s+k \,, \quad  \textup{ and }  \qquad 
    \si^{(s)}_{s+k} = s\enspace .
\end{align*}
Each permutation {{\small $\si^{(s)}$}} is therefore equal to the identity {{\small $( j )_{j\in[n]}$}} up to an exchange of the two indices {{\small $ s $}} and {{\small $s+k$}}. This collection of {{\small $n/4$}} permutations is denoted by {{\small $\mathcal{S}:=\{\si^{(1)},\ldots,\si^{(n/4)}\}$}}. For the subset {{\small $\mathcal{S}\subset   \bPi_n$}}, we  readily check that 
\begin{equation*}
\forall s,t \in \left[\frac{n}{4}\right], \, s\neq t\,: \qquad \quad  L_{\max}(\si^{(t)},\si^{(s)}) \,  \geq \,  \frac{ k}{n} \enspace.
\end{equation*} 
This gives a lower bound on  the packing number {{\small $\mathcal{M}(\ep_n, \mathcal{S}, L_{\max})$}} of radius {{\small $\ep_n$}}:  
\[\mathcal{M}(\ep_n, \mathcal{S}, L_{\max}) \geq n/4\enspace, \qquad \textup{ for} \ \, \ep_n :=  k /n \enspace.\]

To upper bound the {{\small KL}} diameter of {{\small $\mathcal{S}$}}, we use the following lemma whose proof is postponed to the end of the section.
\begin{lem}\label{claim:divergenceBernouilli} For any  $n \times n$ matrix $F$, and $\si, \si' \in \bPi_n$,  we have  $KL(\P_{(F, \si)} \| \P_{(F, \si')}) \leq \frac{1}{2 \va^2}\sum_{i,j\in[n]}  (F_{\si_i \si_j}-F_{\si'_i \si'_j})^2$.
\end{lem}
Combining with the definition of {{\small $F_\al$,}} we obtain for any {{\small $\pi , \pi' \in \mathcal{S}$}},
\[KL(\P_{(F_\al, \si)} \| \P_{(F_\al, \si')})\leq C_2 n \frac{(\al \ep_n)^2}{\va^2} = C_2 C_1^2 \log n \enspace ,\] 
for  {{\small $\ep_n =  k /n =  C_1 (\va/\al ) \sqrt{\log(n)/n}$}}, and a numerical constant {{\small $C_2>0$.}} Taking the value  {{\small $C_1 = (2\sqrt{C_2} )^{-1}$}}, we have   
\[d_{KL}(\mathcal{S}) \leq \frac{\log n}{4}\enspace.\]  
Applying Lemma \ref{fano:prop} to the set {{\small $\mathcal{S}$}}, we arrive at
 \begin{equation*}
      \inf_{\hat \si } \ \, \underset{\si \in \mathcal{S}}{\textup{sup}} \quad  \P_{(F_\al,\si)}\cro{L_{\max}(\hat \si, \si) \geq \frac{\ep_n}{2}} \geq 1 - \frac{(\log(n)/4) + \log 2}{\log(n/4)} \geq \frac{1}{2} \enspace ,
\end{equation*}
as soon as {{\small $n$}} is greater than some numerical constant.  The lower bound {{\small $\ep_n/2$}} is of the order of {{\small $(\va/\al)\sqrt{\log(n)/n}$}}.  Theorem~\ref{1st:lowerBound} follows. \hfill $\square$

\subsection{Proof of Lemma \ref{claim:divergenceBernouilli}}
For all {{\small $i,j\in[n]$}}, we denote the   marginal distribution of  {{\small $A_{ij}$}} by {{\small $\P_{(F_{ij}, \si)}$}}.   By definition of the Kullback-Leibler divergence, we have 
\[KL(\P_{(F, \si)} \| \P_{(F, \si')}) = \sum_{i < j} KL(\P_{(F_{ij}, \si)}, \P_{(F_{ij}, \si')})  \leq \sum_{i,j} (F_{\si_i \si_j}-F_{\si'_i \si'_j})^2 / (2 \va^2)\]
where the  equality follows from the independence of the  {{\small $A_{ij}$,  $i<j$,}} and the inequality from Lemma~\ref{KL-subGauss}. The proof of Lemma~\ref{claim:divergenceBernouilli} is complete. 
\begin{lem}[\cite{Duchi-course}]
\label{KL-subGauss}
    Let two  normal distributions {{\small $P=N(\mu_1,\sigma^2)$}} and {{\small $Q=N(\mu_2,\sigma^2)$}},  with respective means {{\small $\mu_1,\mu_2 \in \mathbb{R}$}} and variance  {{\small $\sigma^2$}}. Then, {{\small $KL(P,Q) =  (\mu_1-\mu_2)^2 / (2\va^2)$}}.
\end{lem}


\newpage

\section{Extension to approximate permutations}
\label{section-new_appendix-extension-rates}

This appendix complements Section~\ref{extension-section} of the main paper.  
We begin by motivating the relaxation to approximate permutations 
through an analysis of the homogeneity constraint inherent to exact permutations.  
We then describe the algorithmic extension of {{\small SABRE}} and present its detailed pseudocode.  
Next, we introduce the average-Lipschitz assumption adapted to this setting, 
and provide the formal version of Theorem~\ref{thm-extension} with explicit tuning parameters.  
Finally, Appendix~\ref{appendix-choice-tuning-approx} details the derivation 
and scaling of these parameters.

\subsection{On the homogeneity restriction of exact permutations}
\label{appendix:homogeneity}

Exact permutations implicitly enforce a strong homogeneity constraint on the latent ordering.  
In many applications, however, the latent positions are unevenly populated: 
some regions correspond to densely sampled items (e.g., tightly connected nodes in a network),
while others are sparse.  
We now illustrate this limitation through a simple calculation.

In the exact permutation setting, consecutive rows of a bi-Lipschitz 
matrix {{\small $F \in \mathcal{BL}(\alpha,\beta)$}} are nearly identical.  More precisely, for any {{\small $i$}},
\[
  |F_{ik} - F_{(i+1)k}| \ \asymp_{\alpha,\beta} \ n^{-1},
\]
so that the squared Euclidean distance between consecutive profiles satisfies
\[
   \min_{k\neq i} \| F_{\pi_i} - F_{\pi_k}\|^2 \ \asymp_{\alpha,\beta} \ n^{-1}.
\]
Thus, consecutive objects cannot be well separated, which imposes a
strong homogeneity constraint.

In contrast, when {{\small $\pi \in \ap(\zeta)$}}, an object {{\small $i$}} may be 
separated from all others by at least {{\small $\zeta$}} in ordering distance, so that
\[
   \min_{k\neq i} \| F_{\pi_i} - F_{\pi_k}\|^2 \ \gtrsim_{\alpha,\beta} \ \zeta^2/n .
\]
Hence, the squared distance between consecutive profiles may diverge,
for instance when {{\small $\zeta \ge \sqrt{n\log n}$}}.  
This simple observation motivates the approximate-permutation framework, 
which removes the homogeneity constraint while preserving the same estimation rate 
(Theorem~\ref{thm-extension}).


\subsection{Algorithmic extension and pseudocode}
\label{appendix-extension-SABRE}

Algorithmically, {{\small SABRE}} extends naturally to the setting of approximate permutations.  
The only modification concerns the data-splitting scheme: instead of the balanced tripartition 
{{\small $(S^1,S^2,S^3)$}} used in the original algorithm, we now employ leave-one-out subsets 
{{\small $S^{ij} = [n]\setminus\{i,j\}$}} for all pairs {{\small $(i,j)$}}.

When latent positions are unevenly spaced, local neighborhoods may contain very few points -- 
sometimes only one per interval of length {{\small $\zeta$}}.  
In such regions, the balanced tripartition {{\small $(S^1,S^2,S^3)$}} can strongly amplify the sparsity:  
within some subsets {{\small $S^t$}}, the distance between nearest neighbors may greatly exceed 
{{\small $\zeta$}}, leading to biased estimates {{\small $\widehat D_{ij}$}} 
in the distance-estimation stage and to unreliable local tests in the refined-seriation stage, 
where the proxy sets {{\small $\tilde L_{ij}$}} (constructed from {{\small $L_{ij}\cap S^t$}}) 
may contain too few points for stable comparison.

To avoid this degeneracy, the extension replaces the tripartition by leave-one-out subsets 
{{\small $S^{ij}=[n]\setminus\{i,j\}$}}, ensuring that each pair {{\small $(i,j)$}} is evaluated 
using all other points as references.  
This modification substantially increases the computational cost 
({{\small $O(n^5)$}} instead of {{\small $O(n^3)$}}) 
but provides greater statistical robustness and stability.

\paragraph*{Pseudocode.}
For completeness, the detailed procedure is given below. We recall that  {{\small $0_{n\times n}$}} denotes the {{\small $n\times n$}} zero matrix.

\begin{algorithm}[H]
    \caption{{{\small \texttt{Extension of Re-evaluate Comparisons}}}} \label{algo:complete-seriation-extension}
    \begin{algorithmic}[1]
        \REQUIRE $(H,D,A,\sigma,\delta_4)$  
        \ENSURE $\widehat H \in \{-1,0,1\}^{n\times n}$
        \STATE Initialize $\widetilde H=0_{n\times n}$
        \FOR{$i<j$ with $H_{ij}=0$}
         \STATE {{\small $S^{ij} = [n]\setminus\{i,j\}$}}
         \STATE Compute {{\small $\widehat{D}^{ij} = \texttt{Estimate} \; \texttt{Distance}(A,S^{ij})$}}
        \STATE  $L_{ij}=\{k : H_{ik}=H_{jk}=1\}$,\quad
               $R_{ij}=\{k : H_{ik}=H_{jk}=-1\}$
        \STATE  $p_{ij} \in \argmin_{p\in S^{ij}} D_{ip}$
        \STATE  $\tilde L_{ij}=\{k\in L_{ij}\cap S^{ij}: 
              \widehat D^{ij}_{p_{ij}k}\ge\delta_4\}$
        \STATE $\tilde R_{ij}=\{k\in R_{ij}\cap S^{ij}: 
              \widehat D^{ij}_{p_{ij}k}\ge\delta_4\}$
        \STATE $\widetilde H_{ij} =$
               \texttt{Evaluate} \texttt{Comparison}$(A_i,A_j,
               \tilde L_{ij},\tilde R_{ij},\sigma)$
        \STATE $\widetilde H_{ji} = -\widetilde H_{ij}$
    \ENDFOR
       \STATE $\widehat H = H + \widetilde H$
    \end{algorithmic}
\end{algorithm}




\subsection{average-Lipschitz structure under approximate permutations}
\label{appendix-extension-approx-assumptions}

We now extend the regularity assumption 
to the approximate-permutation setting of Section~\ref{extension-section}, 
where the latent index vector satisfies 
{{\small $\pi \in \ap(\zeta)$}} for some spacing parameter {{\small $\zeta > 0$}}.
In this setting, the structural conditions are expressed 
on the latent coordinates {{\small $\pi_i$}} instead of the canonical indices {{\small $i$}} of the signal matrix~{{\small $F$}}.

\begin{definition}[average-Lipschitz matrices relative to approximate permutations]
\label{def:AL-approx-perm}
Let {{\small $\alpha, \beta, r, r' > 0$}}.  
A symmetric matrix {{\small $F\in[0,1]^{n\times n}$}} is said to be
\textit{average-Lipschitz relative to the latent ordering~$\pi$} if 
for all {{\small $i,j\in[n]$}} with {{\small $\pi_i < \pi_j$}} the following conditions hold:

\begin{enumerate}[leftmargin=1.3em,label=(\roman*)]

\item \textbf{Local conditions.}
If {{\small $\pi_j-\pi_i \le n r$}}, then
  \begin{enumerate}[leftmargin=1.3em,label=--]
  \item \textit{Average $\ell_2$ upper bound:} {{\small $\|F_{\pi_i}-F_{\pi_j}\|  \ \le\  \beta\,\frac{|\pi_i-\pi_j|}{\sqrt{n}}$}}.
  \item \textit{Average $\ell_1$ lower bound:}
 \[
\Bigg(
     \sum_{k:\, \pi_k < \pi_i - c_0 n} (F_{\pi_i \pi_k} - F_{\pi_j \pi_k})
     \ \vee\
     \sum_{k:\, \pi_k > \pi_j + c_0 n} (F_{\pi_j \pi_k} - F_{\pi_i \pi_k})
     \Bigg)
     \ \ge\ 
     \alpha\,|\pi_i-\pi_j|.\]
where {{\small $c_0 = 1/32$}}.
  \end{enumerate}

\item \textbf{Non-collapse at large distance.}  
If {{\small $\pi_j-\pi_i > n r$}}, then {{\small $\|F_{\pi_i}-F_{\pi_j}\|\ >\  r'\sqrt{n}$}}.
\end{enumerate}

We denote this class by {{\small $\mathcal{AL}_\pi(\alpha,\beta,r,r')$}}.
\end{definition}

When {{\small $\pi$}} is an exact permutation (i.e., {{\small $\zeta=0$}}), 
the class {{\small $\mathcal{AL}_\pi(\alpha,\beta,r,r')$}} coincides with
{{\small $\mathcal{AL}(\alpha,\beta,r,r')$}} in Definition~\ref{defi:PBL}.
More generally, the present formulation simply interprets 
the same Lipschitz regularity along the latent coordinates {{\small $\pi_i$}}.


\subsection{Formal statement and tuning parameters}
\label{appendix-extension-approx-rates}

We now provide the formal version of Theorem~\ref{thm-extension} stated 
in the main paper, together with explicit tuning parameters.

\begin{thm}[Formal version of Theorem~\ref{thm-extension}]
\label{thm-extension-formal}
For any {{\small $(\alpha,\beta,r,r', \sigma)$}} and any sequence 
{{\small $\bar\zeta = (\zeta_n)_{n\ge 1}$}} satisfying {{\small $\zeta_n/n \to 0$}}, 
there exists a constant {{\small $C_{\alpha \beta r r' \sigma \bar\zeta}$}} depending only on 
{{\small $(\alpha,\beta,r,r',\sigma,\bar\zeta)$}} such that the following holds 
for all {{\small $n \ge C_{\alpha \beta r r'\sigma \bar\zeta}$}}.

Suppose 
{{\small $\pi \in \ap(\zeta_n)$}} (Definition~\ref{spatial-sparsity})
and {{\small $F \in \mathcal{AL}_\pi(\alpha,\beta,r,r')$}} (Definition~\ref{def:AL-approx-perm}).

Run the extended {{\small SABRE}} algorithm of Appendix~\ref{appendix-extension-SABRE} with tuning parameters
\[
\delta_1 \;=\; n^{3/4}\log n \;+\; \sqrt{(2\zeta_n+1)\,n}\,
\log\!\Big(\frac{n}{2\zeta_n+1}\Big),
\qquad
\delta_{k+1} \;=\; \delta_k \log\!\Big(\frac{n}{2\zeta_n+1}\Big), \quad k\in[3].
\]
Then, with probability at least {{\small $1-1/n^2$}}, the estimate {{\small $\hat\pi$}} satisfies
\[
L_{\max}(\hat\pi,\pi)
\ \leq\ 
C \; \frac{\sigma}{\alpha}\,\sqrt{\frac{\log n}{n}} ,
\]
for some numerical constant {{\small $C>0$}}.
\end{thm}

\medskip
Thus, in the relaxed setting of approximate permutations, 
the polynomial-time algorithm of Appendix~\ref{appendix-extension-SABRE} 
achieves the same rate as {{\small \textsc{SABRE}}} (Theorem~\ref{thm-matrix-setting}), 
up to constants. The tuning parameters follow the same scaling rules as in the exact permutation analysis,
with adjustments accounting for the spacing {{\small $\zeta_n$}} as detailed in Appendix~\ref{appendix-choice-tuning-approx} below.


\subsection{Choice of tuning parameters}
\label{appendix-choice-tuning-approx}

\paragraph*{High-level reasoning.}
The choice of thresholds {{\small $(\delta_1,\delta_2,\delta_3,\delta_4)$}} follows the same principles as in the exact-permutation setting 
(where {{\small $F\in\mathcal{AL}(\alpha,\beta,r,r')$}}),
but must be adapted under 
{{\small $F\in\mathcal{AL}_\pi(\alpha,\beta,r,r')$}}
when the latent positions are unevenly spaced.
In such cases, local neighborhoods can contain very few points,
which introduces a systematic bias in the empirical distance estimates
and makes the initial seriation stage more fragile.

To mitigate  this effect, the first threshold {{\small $\delta_1$}} 
must be increased to account for the bias 
induced by the spacing parameter {{\small $\zeta_n$}}.  
The subsequent thresholds {{\small $(\delta_k)_{k\ge2}$}} are then defined recursively, following the same logarithmic scaling as in the exact-permutation case.
This choice ensures that the graphs {{\small $\mathcal G_i$}} constructed during the first seriation stage remain sufficiently connected, so that the overall structure of {{\small \textsc{SABRE}}} and its statistical guarantees are preserved.

\smallskip
\paragraph*{Detailed derivation.}
The thresholds must satisfy a set of deterministic inequalities 
with the same structure as 
(\ref{delta-constraint-1rst-seriation}--\ref{constraint-2nd-seriation}),
but adjusted to the approximate permutation setting. When estimating distances in the exact permutation case, the nearest-neighbor distance between consecutive latent positions
is exactly one, which yields an empirical distance error of order 
{{\small $\omega_n \asymp_{\beta,\sigma} n^{3/4}(\log n)^{1/4}$}}. 
When {{\small $\pi \in \ap(\zeta_n)$}}, the nearest-neighbor spacing may grow up to
\[
\eta := \max_{i\in[n]} \min_{k\neq i} |\pi_i - \pi_k| \ \le\ 2\zeta_n + 1 ,
\]
which introduces an additional bias in the nearest-neighbor approximation used to build
the empirical distance matrix~{{\small $\widehat D$}}.  
Repeating the same analysis as in the {{\small \textsc{SABRE}}} proof with this spacing~{{\small $\eta$}} 
shows that the distance estimator now satisfies
\[
\tilde\omega_n 
\ \lesssim\ 
\sqrt{\beta (2\zeta_n+1)\,n} 
\ +\ 
\sqrt{(\sigma + 1)\sigma}\,n^{3/4}(\log n)^{1/4} .
\]
This defines the modified estimation error~{{\small $\tilde\omega_n$}} 
that replaces~{{\small $\omega_n$}} in the tuning inequalities.

The connectivity condition ensuring that consecutive items in the latent order
are linked in the graph built by {{\small\textsc{Aggregate} \textsc{Bisections}}}
must therefore hold with the new spacing parameter~{{\small $\eta$}}.  
This leads to the generalized requirement
\[
\tilde\omega_n + \beta(2\zeta_n+1) \ \le\ \delta_1 ,
\]
which extends the original condition~{{\small $\omega_n + \beta \le \delta_1$}}
from the exact permutation setting.  
Under these modified constraints (and similar ones obtained by revisiting the proof 
of Theorem~\ref{thm-matrix-setting-general}), the same arguments apply.

Accordingly, the first threshold {{\small $\delta_1$}} combines two dominant contributions:  
the term {{\small $n^{3/4}\log n$}} controls stochastic fluctuations of the empirical distances,
as in the exact-permutation analysis; 
while the additional term 
{{\small $\sqrt{(2\zeta_n+1)n}\,\log\!\big(\tfrac{n}{2\zeta_n+1}\big)$}} 
accounts for the bias induced by the enlarged nearest-neighbor spacing {{\small $(2\zeta_n+1)$}}.  
This yields
\[
\delta_1 \;=\; n^{3/4}\log n 
\;+\; 
\sqrt{(2\zeta_n+1)\,n}\,
\log\!\Big(\frac{n}{2\zeta_n+1}\Big),
\]
and, recursively, the subsequent thresholds 
{{\small $(\delta_k)_{k\le4}$}} 
are defined as 
{{\small $\delta_{k+1}=\delta_k\log\!\big(\tfrac{n}{2\zeta_n+1}\big)$}},
following the same recursive rule as in the original tuning scheme.





\section{Concentration inequalities}

\begin{lem}[Hypergeometric distribution]
For {{\small $p \in [0,1]$}} and {{\small $1\leq N \leq n$}}, let {{\small $X$}} be  a hypergeometric random variable with parameters {{\small $(N,p,n)$.}} 
Then, for all {{\small $t>0$,}}
\begin{equation}
\label{eq_Hoeffding_hypergeom}
\P\ac{ |X - Np| \ \geq \ \sqrt{\frac{Nt}{2}} } \ \  \leq \ \ 2 e^{-t} \enspace.
\end{equation}
\end{lem}

\begin{lem}[Hoeffding-type bound for sampling without replacement] 
Let {{\small $\Xcal =\{x_1,\ldots,x_n\}$}} where {{\small $x_i >0$}} for all {{\small $i$}}. If {{\small $X_1,\ldots, X_N$}} is a random sample drawn without replacement from {{\small $\Xcal$}}, then  for all {{\small $t>0$}}, 
\begin{equation}
\label{eq_Hoeffding_sample}
  \P\pa{\bigg{|}\frac{1}{N} \sum_{i=1}^{N} X_i - \frac{1}{n}\sum_{i=1}^n x_i \bigg{|}\, \geq \, t} \ \ \leq \  \ 2 \exp\pa{- \frac{2N t^2}{\max_{i} x_i^2}}\enspace.
\end{equation}
\end{lem}

\smallskip 

\begin{lem}[Bernstein–type bound for sampling without replacement {\cite[Cor.~3.6]{bardenet2015concentration}}]
\label{bernstein_sample-no-replacement}
Let {{\small $\Xcal =\{x_1,\ldots,x_n\}$}} where {{\small $x_i >0$}} for all {{\small $i$}}, with population  mean $\mu := \tfrac{1}{n}\sum_{i=1}^n x_i$ and population variance $\sigma^2 := \tfrac{1}{n}\sum_{i=1}^n (x_i-\mu)^2$. 
Draw $X_1,\dots,X_N$ without replacement from this population, with $1\le N\le n$. 
Then for any $\delta\in(0,1)$, with probability at least $1-2\delta$,
\[
\frac{1}{N}\sum_{i=1}^N X_i  
\;\ge\; \mu \; - \;  
\sigma\,\sqrt{\frac{2\,\log(1/\delta)}{N}}
\;-\;
\frac{7}{3}\,\frac{(\max_{i} x_i)\,\log(1/\delta)}{N}\, .
\]
\end{lem}

\end{appendix}



\end{document}